\documentclass[matprg]{svjour3}

\usepackage{url}
\usepackage{caption}
\usepackage{mathtools}
\usepackage{bbm}
\usepackage{enumitem}
\usepackage{esvect}
\usepackage{centernot}
\usepackage{float} 
\usepackage{fullpage}
\usepackage{amsfonts}
\usepackage{latexsym}
\usepackage{amssymb}
\usepackage{color}
\usepackage{hyperref}
\usepackage{graphicx}
\graphicspath{{./figures/}}
\usepackage{epstopdf}
\usepackage{algorithm,algpseudocode}
\usepackage{verbatim}
\usepackage{mathtools}
\usepackage{enumitem}
\usepackage{esvect}
\usepackage{centernot}
\usepackage{float} 
\usepackage{siunitx}
\usepackage{mathrsfs}
\usepackage[flushleft]{threeparttable}
\usepackage{textcomp}
\usepackage{booktabs}
\usepackage{tikz}
\usetikzlibrary{positioning,angles,quotes}
\usepackage{empheq}
\usepackage{amsmath}
\usepackage{amsfonts}
\usepackage{amssymb}
\usepackage{makecell}

\usepackage{tikz}
\usepackage{pgfplots}
\usetikzlibrary{bayesnet}
\usetikzlibrary{arrows}

\numberwithin{equation}{section}

\allowdisplaybreaks
\def\vgap{\vspace*{.1in}}
\newtheorem{assumption}{Assumption}

\newcommand{\E}{\mathbb{E}}

\newcommand{\R}{\mathbb{R}}

\newcommand{\La}{\mathcal{L}}

\def \grad {\nabla}
\DeclareMathOperator*{\argmax}{arg\,max}
\DeclareMathOperator*{\argmin}{arg\,min}

\newcommand{\inner}[2]{\langle {#1,#2} \rangle}
\NewDocumentCommand\norm{ m O{}}{\lVert #1 \rVert_{#2}}

\def\endproof{{\ \hfill\hbox{%
      \vrule width1.0ex height1.0ex
    }\parfillskip 0pt}\par}

\newcommand{\tsum}{\textstyle \sum}

\newcommand{\ep}{\epsilon}

\makeatletter
\newcommand{\algmargin}{\the\ALG@thistlm}
\makeatother
\newlength{\whilewidth}
\algdef{SE}[parWHILE]{parWhile}{EndparWhile}[1]
  {\parbox[t]{\dimexpr\linewidth-\algmargin}{%
     \hangindent\whilewidth\strut\algorithmicwhile\ #1\ \algorithmicdo\strut}}{\algorithmicend\ \algorithmicwhile}%
\algnewcommand{\parState}[1]{\State%
  \parbox[t]{\dimexpr\linewidth-\algmargin}{\strut #1\strut}}

\def \bigO {\mathcal{O}}

\usepackage{pifont}

\newcommand{\probP}{\text{I\kern-0.15em P}}
\newcommand{\ts}{\textsuperscript}
\newcommand{\addt}[1]{#1^t}
\newcommand{\addtt}[1]{#1^{t-1}}
\newcommand{\addttt}[1]{#1^{t-2}}
\newcommand{\addp}[1]{#1^{t+1}}
\newcommand{\addT}[1]{#1^N}
\newcommand{\addTT}[1]{#1^{N-1}}

\newcommand{\addO}[1]{#1^1}

\newcommand{\subaddt}[1]{#1_t}
\newcommand{\subaddtt}[1]{#1_{t-1}}

\newcommand{\subaddT}[1]{#1_N}

\NewDocumentCommand\fancySubAdd{m m m}{%
    \IfNoValueTF{#1}{%
        #2^{#3}
    }{%
        #2[#1]^{#3}
    }%
}

\NewDocumentCommand\fancyAddT{m m}{\fancySubAdd{#1}{#2}{N}}
\NewDocumentCommand\fancyAddTT{m m}{\fancySubAdd{#1}{#2}{N-1}}
\NewDocumentCommand\fancyAddt{m m}{\fancySubAdd{#1}{#2}{t}}
\NewDocumentCommand\fancyAddtt{m m}{\fancySubAdd{#1}{#2}{t-1}}
\NewDocumentCommand\fancyAddp{m m}{\fancySubAdd{#1}{#2}{t+1}}
\NewDocumentCommand\fancyAddttt{m m}{\fancySubAdd{#1}{#2}{t-2}}
\NewDocumentCommand\fancyAddO{m m}{\fancySubAdd{#1}{#2}{0}}

\def \xbar {\bar{x}}
\def \xbarT {\addT{\xbar}}

\def \zbar {\bar{z}}

\def \pibar {\bar{\pi}}
\NewDocumentCommand\myLa{O{x} O{\pi} O{}}{\La_{#3}(#2;#1)}
\NewDocumentCommand\nmyLa{O{x} O{\pi} O{}}{\La_{#3}(#1;#2)}
\NewDocumentCommand\myLalam{O{x} O{\pi} O{}}{\La^\lambda_{#3}(#2;#1)}
\def \xstar {x^*}

\def \zstar {z^*}

\def \Mg {M_g}

\def \DX {D_X}

\NewDocumentCommand\xt{O{t}}{x^{#1}}
\def \xtt {\addtt{x}}
\def \xttt {\addttt{x}}
\def \xp {\addp{x}}
\def \xT  {\addT{x}}

\def \z {z}

\def \zt {\addt{\z}}

\def \xtil {\tilde{x}}
\NewDocumentCommand\xtlt{O{t}}{\tilde{x}^{#1}}
\NewDocumentCommand\xtilt{O{t}}{\xtil^{#1}}

\NewDocumentCommand\wt{O{t}}{\omega_{#1}}

\def \wtt {\subaddtt{w}}

\def \wT {\subaddT{w}}

\NewDocumentCommand\thetat{O{t}}{\theta_{#1}}

\def \taut {\subaddt{\tau}}
\def \tauT {\subaddT{\tau}}

\NewDocumentCommand\etat{O{t}}{{\eta_{#1}}}

\def \etatt {\subaddtt{\eta}}

\def \etaT {\subaddT{\eta}}
\def \sp {\text{span}}

\NewDocumentCommand\Mit{O{1} O{t}}{\mathcal{M}_{#1}^{#2}}
\NewDocumentCommand\Xt{O{t}}{\mathcal{X}^{#1}}

\NewDocumentCommand\etai{O{i}}{\eta_{#1}}
\NewDocumentCommand\etait{o}{\fancyAddt{#1}{\etai}}
\NewDocumentCommand\etaitt{o}{\fancyAddtt{#1}{\etai}}
\NewDocumentCommand\etaittt{o}{\fancyAddttt{#1}{\etai}}
\NewDocumentCommand\etaip{o}{\fancyAddp{#1}{\etai}}
\NewDocumentCommand\etaiT{o}{\fancyAddTT{#1}{\etai}}
\NewDocumentCommand\etaiO{o}{\fancyAddO{#1}{\etai}}

\NewDocumentCommand\taui{O{i}}{\tau_{#1}}
\NewDocumentCommand\tauit{O{i} O{t}}{{\tau}_{#1, #2}}
\NewDocumentCommand\tauip{O{i}}{{\tau}_{#1, t+1}}
\NewDocumentCommand\tauiT{O{i}}{{\tau}_{#1, N-1}}
\NewDocumentCommand\tauitt{O{i}}{{\tau}_{#1, t-1}}
\NewDocumentCommand\tauiO{O{i}}{{\tau}_{#1, 0}}

\NewDocumentCommand\tauisig{O{i}}{\tau_{#1,\sigma}}
\NewDocumentCommand\tauisigt{O{i}}{\addt{\tau}_{#1,\sigma}}
\NewDocumentCommand\tauisigp{O{i}}{\addp{\tau}_{#1,\sigma}}
\NewDocumentCommand\tauisigT{O{i}}{\addTT{\tau}_{#1,\sigma}}
\NewDocumentCommand\tauisigtt{O{i}}{\addtt{\tau}_{#1,\sigma}}
\NewDocumentCommand\tauisigO{O{i}}{\addO{\tau}_{#1,\sigma}}

\NewDocumentCommand\tauipi{O{i}}{\tau_{#1,\pi}}
\NewDocumentCommand\tauipit{O{i}}{\addt{\tau}_{#1,\pi}}
\NewDocumentCommand\tauipip{O{i}}{\addp{\tau}_{#1,\pi}}
\NewDocumentCommand\tauipiT{O{i}}{\addTT{\tau}_{#1,\pi}}
\NewDocumentCommand\tauipitt{O{i}}{\addtt{\tau}_{#1,\pi}}
\NewDocumentCommand\tauipiO{O{i}}{\addO{\tau}_{#1,\pi}}

\NewDocumentCommand\hpii{O{i}}{\hat{\pi}_{#1}}
\NewDocumentCommand\hpiit{O{i}O{t}}{\hat{\pi}^{#2}_{#1}}

\NewDocumentCommand\hxi{O{i}}{\hat{\xi}_{#1}}

\NewDocumentCommand\Pii{O{i}}{\Pi_{#1}}
\NewDocumentCommand\Piit{O{t} O{1} }{{\Pi}^{#1}_{#2}}
\NewDocumentCommand\Piitl{O{i}}{\tilde{\Pi}_{#1}}
\def \La {\mathcal{L}}
\def \hLa {\hat{\mathcal{L}}}

\def \tilLa {{\La}}
\NewDocumentCommand\mytLa{O{x} O{\pi}}{\tilLa(#1; #2)}

\NewDocumentCommand\myfi{O{i}}{f_{{#1}}}
\NewDocumentCommand\myfipr{O{i}}{f_{#1}'}
\NewDocumentCommand\fii{O{i}}{f_{#1}}
\NewDocumentCommand\fistar{O{i}}{f^*_{#1}}

\NewDocumentCommand\yi{O{i}}{y_{#1}}
\NewDocumentCommand\yit{O{i}O{t}}{{y}^{#2}_{#1}}
\NewDocumentCommand\yip{O{i}}{\addp{y}_{#1}}
\NewDocumentCommand\yitil{O{i}}{\tilde{y}_{#1}}
\NewDocumentCommand\yitilt{O{i} O{t}}{\tilde{y}^{#2}_{#1}}
\NewDocumentCommand\uyi{O{i}}{\underline{y}_{#1}}
\NewDocumentCommand\uyit{O{i} O{t}}{\underline{y}^{#2}_{#1}}
\NewDocumentCommand\uyitE{O{i} O{t}}{\underline{y}^{#2}_{#1, \E}}
\NewDocumentCommand\uyitt{o}{\fancyAddtt{#1}{\uyi}}
\NewDocumentCommand\uyip{o}{\fancyAddp{#1}{\uyi}}
\NewDocumentCommand\uyiO{o}{\fancyAddO{#1}{\uyi}}

\NewDocumentCommand\byi{O{i}}{\bar{y}_{#1}}
\NewDocumentCommand\byit{O{i}}{\addt{\bar{y}}_{#1}}
\NewDocumentCommand\byiT{O{i}}{\addT{\bar{y}}_{#1}}
\NewDocumentCommand\byip{O{i}}{\addp{\bar{y}}_{#1}}
\NewDocumentCommand\yibar{O{i}}{\bar{y}_{#1}}
\NewDocumentCommand\Mpi{O{}}{M_{\pi_{#1}}}
\newcommand{\dom}[1]{\text{dom}#1}

\NewDocumentCommand\pii{O{i}}{\pi_{#1}}
\NewDocumentCommand\piit{O{i} O{t}}{\pii[#1]^{#2}}
\NewDocumentCommand\piitt{o}{\fancyAddtt{#1}{\pii}}
\NewDocumentCommand\piittt{o}{\fancyAddttt{#1}{\pii}}
\NewDocumentCommand\piip{o}{\fancyAddp{#1}{\pii}}
\NewDocumentCommand\piiT{o}{\fancyAddT{#1}{\pii}}
\NewDocumentCommand\piiO{o}{\fancyAddO{#1}{\pii}}

\NewDocumentCommand\piistar{O{i}}{\pi^*_{#1}}

\NewDocumentCommand\Mpii{O{i}}{M_{\Pi_{#1}}}

\def \pik {\pii[1:]}

\def \pikt {\piit[1:]}

\def \upi {\underline{\pi}}
\NewDocumentCommand\upii{O{i}}{\upi_{#1}}
\NewDocumentCommand\upiit{o}{\fancyAddt{#1}{\upii}}
\NewDocumentCommand\upiitt{o}{\fancyAddtt{#1}{\upii}}
\NewDocumentCommand\upiittt{o}{\fancyAddttt{#1}{\upii}}
\NewDocumentCommand\upiip{o}{\fancyAddp{#1}{\upii}}
\NewDocumentCommand\upiiT{o}{\fancyAddT{#1}{\upii}}

\NewDocumentCommand\piibar{O{i}}{\bar{\pi}_{#1}}
\NewDocumentCommand\piibart{o}{\fancyAddt{#1}{\piibar}}
\NewDocumentCommand\piibartt{o}{\fancyAddtt{#1}{\piibar}}
\NewDocumentCommand\piibarttt{o}{\fancyAddttt{#1}{\piibar}}
\NewDocumentCommand\piibarp{o}{\fancyAddp{#1}{\piibar}}
\NewDocumentCommand\piibarT{o}{\fancyAddT{#1}{\piibar}}

\NewDocumentCommand\piitil{O{i}}{\tilde{\pi}_{#1}}
\NewDocumentCommand\piitilp{o}{\fancyAddp{#1}{\piitil}}

\NewDocumentCommand\pixij{O{i} O{}}{\pi_{#1}(\xi_{#1}^{#2})}
\NewDocumentCommand\pixijt{O{i} O{j}}{\addt{\pi}_{#1}(\xi_{#1}^{#2})}
\NewDocumentCommand\pihxit{O{i} O{t}}{\pi^{#2}_{#1}(\hat{\xi}_{#1})}
\NewDocumentCommand\pixijp{O{i} O{}}{\addp{\pi}_{#1}(\xi_{#1}^{#2})}
\NewDocumentCommand\pixijtt{O{i} O{}}{\addtt{\pi}_{#1}(\xi_{#1}^{#2})}
\NewDocumentCommand\pixijO{O{i} O{}}{\addO{\pi}_{#1}(\xi_{#1}^{#2})}
\NewDocumentCommand\gradfijt{O{i} O{j} O{t}}{\grad f_{#1}(\uyit[#1][#3], \xi_{#1}^{#2})}
\NewDocumentCommand\gradfijht{O{i} O{j} O{t-1}}{\grad f_{#1}(\uyit[#1][#3], \hat{\xi}_{#1}^{#2})}
\NewDocumentCommand\gradfit{O{i} O{t}}{\grad f_{#1}(\uyit[#1][#2])}
\NewDocumentCommand\Laijtt{O{i} O{j} O{t-1} O{t}}{\La_{#1}(x^{#3}; \piit[\ito[#1]][#4](\xi_{\ito[#1]}^{#2}))}
\NewDocumentCommand\hLaijtt{O{i} O{j} O{t-1} O{t}}{\hat{\La}_{#1}(x^{#3}; \piit[\ito[#1]][#4](\xi_{\ito[#1]}^{#2}))}
\NewDocumentCommand\Laitt{O{i}  O{t-1} O{t}}{\La_{#1}(x^{#2}; \piit[\ito[#1]][#3])}

\def \nsumt {\tsum_{t=1}^N}
\def \sumt {\tsum_{t=1}^{N}}

\def \sumwt {\tsum_{t=1}^{N} \wt}

\def \sumk {\tsum_{i=1}^k}

\NewDocumentCommand\Lfi{O{i}}{L_{{#1}}}
\NewDocumentCommand\Li{O{i}}{L_{{#1}}}

\NewDocumentCommand\Vi{O{1}}{V_{#1}}
\NewDocumentCommand\vi{O{1}}{v_{#1}}
\NewDocumentCommand\vit{O{t}O{1}}{v_{#2}^{#1}}

\NewDocumentCommand\Vii{O{i}}{V_{#1}}
\NewDocumentCommand\vii{O{i}}{v_{#1}}
\NewDocumentCommand\viit{O{i}O{t}}{v_{#1}^{#2}}
\NewDocumentCommand\viitt{O{i}O{t-1}}{v_{#1}^{#2}}

\NewDocumentCommand\vixijt{O{i}O{j}O{t}}{v^{#3}_{#1}(\xi_{#1}^{#2})}
\NewDocumentCommand\vixijtt{O{i}O{j}O{t-1}}{v^{#3}_{#1}(\xi_{#1}^{#2})}
\NewDocumentCommand\vvit{O{i:}O{t}}{v_{#1}^{#2}}
\NewDocumentCommand\vitt{O{t-1}O{1}}{v_{#2}^{#1}}
\NewDocumentCommand\gamit{O{1}O{t}}{\gamma_{#1, #2}}
\NewDocumentCommand\gamiit{O{i}O{t}}{\gamma_{#1, #2}}
\NewDocumentCommand\vibar{O{1}}{\bar{v}_{#1}}
\NewDocumentCommand\vistar{O{1}}{v^*_{#1}}
\def \ybar {\bar{y}}
\def \Lg {L_g}
\def \gstar {g^*}
\def \Dgstar {D_{g^*}}
\def \pibar {\bar{\pi}}
\def \uw {u_w}
\def \guw {g_{\uw}}
\def \transpose {^\top}
\def \guwstar {\guw^*}

\NewDocumentCommand\Dfi{O{i}}{D_{f_{#1}}}
\NewDocumentCommand\Dfistar{O{i}}{D_{f^*_{#1}}}
\NewDocumentCommand\dDfistar{O{i}}{\mathcal{D}_{f^*_{#1}}}
\def \frN {\mathfrak{N}}
\NewDocumentCommand\frNi{O{1:i}}{\frN_{#1}}
\def \frNo {\mathfrak{N}_o}
\def \frS {\mathfrak{S}}
\def \frP {\mathfrak{P}}

\def \Dx {\mathcal{D}_X}
\NewDocumentCommand\DVi{O{1}}{D_{V_{#1}}}

\NewDocumentCommand\Dpii{O{i}}{\mathcal{D}_{\Pi_{#1}}}
\NewDocumentCommand\Dpiio{O{i}}{\mathcal{D}_{\Pi_{#1}}}

\NewDocumentCommand\sig{O{}}{\sigma_{#1}}
\NewDocumentCommand\sigtl{O{}}{\tilde{\sigma}_{#1}}
\NewDocumentCommand\sigfi{O{i}}{\sigma_{f_#1}}
\NewDocumentCommand\sigpii{O{i}}{\sigma_{\pi_{#1}}}
\NewDocumentCommand\sigLai{O{i}}{\sigma_{\La_{#1}}}

\NewDocumentCommand\pji{O{j}O{i}}{\pi_{#1, #2}}
\NewDocumentCommand\pjit{O{j}O{i}}{\addt{\pi}_{#1, #2}}
\NewDocumentCommand\pjitt{O{j}O{i}}{\addtt{\pi}_{#1, #2}}
\NewDocumentCommand\pjip{O{j}O{i}}{\addp{\pi}_{#1, #2}}
\NewDocumentCommand\pjiT{O{j}O{i}}{\addT{\pi}_{#1, #2}}

\NewDocumentCommand\pjitilp{O{i}}{\tilde{\pi}^{t+1}_{#1}}
\NewDocumentCommand\pjitil{O{i}}{\tilde{\pi}_{#1}}

\NewDocumentCommand\Mfi{O{i-1}}{\mathcal{M}_{f_{#1}}}
\NewDocumentCommand\Mi{O{i}}{M_{#1}}
\NewDocumentCommand\Mitl{O{}}{{M}_{#1}}
\NewDocumentCommand\Mitil{O{}}{\tilde{M}_{#1}}
\NewDocumentCommand\Mih{O{}}{\hat{M}_{#1}}
\NewDocumentCommand\Mibar{O{}}{\bar{M}_{#1}}

\def \Ltl {\tilde{L}}

\NewDocumentCommand\Lfiitil{O{i}}{\tilde{L}_{f_{#1}}} 
\NewDocumentCommand\Mpiitil{O{i}} {\tilde{M}_{\Pii[#1]}}

\NewDocumentCommand\SigPiitil{O{i}} {\tilde{\sigma}_{\pii[#1]}}

\NewDocumentCommand\SigLaii{O{i}}{\sigma_{\La_{#1}}}

\NewDocumentCommand\Ait{O{i}} {\mathcal{A}_{#1}^t}
\NewDocumentCommand\Tit{O{i}} {\mathcal{T}_{#1}^t}
\NewDocumentCommand\Cit{O{i}} {\mathcal{C}_{#1}^t}
\NewDocumentCommand\Delxt{O{i}} {\Delta^t_{x,#1}}

\NewDocumentCommand\Bixt{O{i}} {\Delxt[#1]}
\NewDocumentCommand\Delpit{O{i}} {\Delta^t_{\pi_{#1}}}
\NewDocumentCommand\Bipit{O{i}} {\Delpit[#1]}
\NewDocumentCommand\Bideltat{O{i}} {\delit[#1]}
\NewDocumentCommand\Bit{O{i}} {\mathcal{B}_{#1}^t}
\NewDocumentCommand\Delit{O{i}}{\Delta^t_{#1}}

\NewDocumentCommand\delit{O{i}}{\delta_{#1}^t}

\NewDocumentCommand\normsq{m}{\norm{#1}^2}

\def \piQ {\pi_q}

\NewDocumentCommand\piiP{O{i}}{\pi_{p_{#1}}}
\NewDocumentCommand\piiPt{O{i}}{\addt{\piiP[#1]}}
\NewDocumentCommand\piiPp{O{i}}{\addp{\piiP[#1]}}
\NewDocumentCommand\piiPtt{O{i}}{\addtt{\piiP[#1]}}
\NewDocumentCommand\piiPttt{O{i}}{\addttt{\piiP[#1]}}
\NewDocumentCommand\piiPT{O{i}}{\addT{\piiP[#1]}}

\NewDocumentCommand\piiQ{O{i}}{\pi_{q_{#1}}}
\NewDocumentCommand\piiQt{O{i}}{\addt{\piiQ[#1]}}
\NewDocumentCommand\piiQp{O{i}}{\addp{\piiQ[#1]}}
\NewDocumentCommand\piiQtt{O{i}}{\addtt{\piiQ[#1]}}
\NewDocumentCommand\piiQttt{O{i}}{\addttt{\piiQ[#1]}}
\NewDocumentCommand\piiQT{O{i}}{\addT{\piiQ[#1]}}

\NewDocumentCommand\piQxii{O{i}}{\piQ(\xi_q^{#1})}
\NewDocumentCommand\piQxiit{O{i}}{\addt{\piQ}(\xi_q^{#1})}
\NewDocumentCommand\piQxiip{O{i}}{\addp{\piQ}(\xi_q^{#1})}
\NewDocumentCommand\piQxiitt{O{i}}{\addtt{\piQ}(\xi_q^{#1})}
\NewDocumentCommand\piQxiittt{O{i}}{\addttt{\piQ}(\xi_q^{#1})}
\NewDocumentCommand\piQxiiT{O{i}}{\addT{\piQ}(\xi_q^{#1})}

\NewDocumentCommand\MPi{O{}}{M_{p_{#1}}}
\NewDocumentCommand\MQi{O{}}{M_{q_{#1}}}
\NewDocumentCommand\MLaii{O{i+1}}{M_{\La_{#1:}}}

\NewDocumentCommand\Hpt{O{t}}{H_+^{#1}}

\NewDocumentCommand\Hmt{O{t}}{H_-^{#1}}

\def \piitil {\tilde{\pi}_i}

\NewDocumentCommand\Yi{O{i}}{Y_{#1}}
\NewDocumentCommand\Yit{O{t}O{1}}{\mathcal{Y}^{#1}_{#2}}
\NewDocumentCommand\ito{O{i}}{#1 :}

\NewDocumentCommand\myfixi{O{i}O{x}O{j}}{f_{#1}(#2, \xi_{#1}^{#3})}
\NewDocumentCommand\myfixipr{O{i}O{x}O{j}}{f'_{#1}(#2, \xi_{#1}^{#3})}
\NewDocumentCommand\fistarxi{O{i}O{\pii}O{}}{f^*_{#1}(#2, \xi_{#1}^{#3})}

\def \Var {\text{Var}}
\NewDocumentCommand\Oi{O{i}}{\mathcal{O}_{#1}}
\NewDocumentCommand\SOi{O{i}}{\mathcal{SO}_{#1}}
\NewDocumentCommand\DSOi{O{i}}{\mathcal{DSO}_{#1}}
\def \ip {i+1}

\NewDocumentCommand\ifrom{O{i-1}}{1:{#1}}

\NewDocumentCommand\partp{O{l}}{p_{#1}}
\NewDocumentCommand\partq{O{l}}{q_{#1}}
\NewDocumentCommand\STCi{O{i}}{\text{STC}_{#1}}

\NewDocumentCommand\myhi{O{i}}{h_{#1}}

\NewDocumentCommand\myfiga{O{}}{f^\gamma_{#1}}

\NewDocumentCommand\upperi{O{i}} {^{(#1)}}

\NewDocumentCommand\Ui{O{i}}{U_{#1}}

\def \targmax {{\textstyle\argmax}}
\def \targmin {{\textstyle\argmin}}
\def \tmax {{\textstyle\max}}
\def \tmin {{\textstyle\min}}

\title{Optimal Algorithms for Convex Nested Stochastic Composite Optimization
\thanks{
This research was partially supported by the ARO grant W911NF-18-1-0223 and the NSF grant 1953199.}}
\author{
Zhe Zhang\thanks{H. Milton Stewart School of Industrial and Systems
    Engineering, Georgia Institute of Technology, Atlanta, GA, 30332.
    (email: {\tt jimmy\_zhang@gatech.edu}).}
    \and Guanghui Lan 
    \thanks{H. Milton Stewart School of Industrial and Systems
    Engineering, Georgia Institute of Technology, Atlanta, GA, 30332.
    (email: {\tt george.lan@isye.gatech.edu}).}
}

\begin{document}
\maketitle
\begin{abstract}Recently, convex nested stochastic composite optimization (NSCO) has received considerable attention for its applications in reinforcement learning and risk-averse optimization. The current NSCO algorithms have worse stochastic oracle complexities, by orders of magnitude, than those for simpler stochastic composite optimization problems (e.g., sum of smooth and nonsmooth functions) without the nested structure. Moreover, they require all outer-layer functions to be smooth, which is not satisfied by some important applications.
These discrepancies prompt us to ask: ``does the nested composition make stochastic optimization more difficult in terms of the order of oracle complexity?"  In this paper, we answer the question by developing order-tight algorithms for the convex NSCO problem constructed from an arbitrary composition of smooth, structured non-smooth and general non-smooth layer functions. When all outer-layer functions are smooth, we propose a stochastic sequential dual (SSD) method to achieve an oracle complexity of $\mathcal{O}(1/\epsilon^2)$ ($\mathcal{O}(1/\epsilon)$) when the problem is non-strongly (strongly) convex. When there exists some structured non-smooth or general non-smooth outer-layer function, we propose a nonsmooth stochastic sequential dual (nSSD) method to achieve an oracle complexity of $\mathcal{O}(1/\epsilon^2)$. We provide a lower complexity bound to show the latter $\mathcal{O}(1/\epsilon^2)$ complexity to be unimprovable even under a strongly convex setting. All these complexity results seem to be new in the literature and they indicate that the convex NSCO problem has the same order of oracle complexity as those without the nested composition in all but the strongly convex and outer-non-smooth problem.  

\end{abstract}


\section{Introduction} \label{sec_intro}
\subsection{Motivation}
Composite optimization has attracted
considerable interest for its applications in compressed sensing, image processing and machine learning.
Many algorithmic studies have been focused on composite optimization
of the form $\min_{x \in X}f(x) + g(x)$, where $f$
is a smooth convex function and $g$ is a nonsmooth function with certain special structures. Optimal first-order methods have been developed in \cite{nesterov2007gradient,tseng2008accelerated,beck2009fast,Lan14-1,lanOuyang2016GradientSliding} for solving these problems 
under different assumptions about $g$. In the stochastic setting, Lan~\cite{Lan08,Lan08-1} presented
an accelerated stochastic approximation method that can achieve the
optimal iteration/sampling complexity when one only has access to stochastic (sub)gradients of the objective function (see, e.g.,~\cite{ghadimi2012optimal,GhaLan13-1,chen2014optimal} for extensions).

The study of composite optimization has later expanded to more complex nested composition problems. Specifically, Lewis and Wright \cite{lewis2016proximal}  developed a globally convergent algorithm for solving $\min_{x \in X}f(g(x))$, where the outer layer function $f$ can be non-smooth, non-convex and extended-valued. Lan~\cite{Lan15Bundle} also studied the complexity of these problems when $f$ is relatively simple.
Wang et al. \cite{wang2017accelerating} are the first to study
nested stochastic composite optimization (NSCO) problems when 
$f$ and $g$ are given as expectation functions. Since NSCO finds wide applications in reinforcement learning \cite{wang2017accelerating}, meta-learning \cite{chen2020solving}, and risk-averse optimization \cite{ruszczynski2020stochastic}, it becomes an important topic in stochastic optimization.

A key challenge in NSCO is the lack of unbiased gradient estimators for the nested function. This issue can be illustrated with a simple two-layer problem, $\min_{x \in X} \{f(x) := f_1 (f_2(x))\}.
$
By the chain rule, the (sub)gradient of $f$ is given  by 
$f'(x) = f_1'(f_2(x))f'_2(x)\footnote{We use the Jacobian matrix $\myfi'(y) \in  R^{n\times m}$ to represent the (sub)-gradient of $f_i: \R^m \rightarrow \R^n$. It helps to simplify the notation\textbf{} when deriving the gradient of a nested composite function.}$.
Now assume $f_2$ is accessible only through a stochastic first-order oracle which returns unbiased estimators $(f_2(x, \xi), f'_2(x, \xi)),
$ then the point at which to evaluate $\myfi[1]'$, $f_2(x)$, is not available. If the sub-gradient is evaluated instead at some  $\myfi[2](x, \xi)$, the nested estimator is biased except for an affine $\myfi[1]$, i.e.
\begin{equation}\label{eq:bias_illus}
\E[\myfi[1]'( \myfi[2](x, \xi))] \neq \grad \myfi[1]' (\myfi[2](x)).\end{equation}
So it is not possible to obtain an unbiased estimator for $f_1'(f_2(x))$ nor $f_1'(f_2(x))\myfi[2]'(x)$.

The current approach to resolve the aforementioned challenge in convex NSCO is gradient approximation, i.e., control the bias of the left-hand side of \eqref{eq:bias_illus} by evaluating $\myfi[1]'$ at some "close enough" estimate to $\myfi[2](x)$. 
For example, if $\myfi[1]$ is assumed to be $L_1$-Lipschitz smooth, the empirical average from $N=\bigO(1/\ep^2)$ samples, $\ybar(\xi):=\sum_{i=1}^{N}\myfi[2](x, \xi_i)/N$, satisfies $\E[\norm{\ybar(\xi) - \myfi[2](x)}] \leq \bigO(\ep)$ such that 
$\E[\norm{\myfi[1]'(\ybar(\xi)) -  \myfi[1]' (\myfi[2](x))}]\leq \bigO(\Lfi[1]\ep).$
So the stochastic gradient descent method with $f_1'(\ybar(\xi))\myfi[2]'(x)$ as the gradient-proxy can  find an $\ep$-optimal solution, i.e., $\E[f(x^N) - f(\xstar)] \leq \ep$, within $\bigO(1/\ep^4)$ calls to the stochastic oracle. Indeed, this simple method achieves the same oracle complexity as the SCGD algorithm in \cite{mengdi2017stochastic}. However, the $\ybar(\xi)$-estimation is more subtle in the SCGD algorithm since it is computed as the faster part of a two-time-scale scheme. 

Comparing this result to that for solving the simple one-layer stochastic optimization problem without the nested structure, we notice two discrepancies. First, the gradient-approximation type methods require the outer-layer function to have a  Lipschitz-continuous gradient so that the desired gradient can be estimated by estimating its argument. Second, the oracle complexity are worse than the $\bigO(1/\ep^2)$ complexity for the one-layer problem by orders of magnitude. These observations motivate us to ask the following research question:
\begin{empheq}[box=\fbox]{align*}
&\text{Does the nested composition make stochastic composite optimization more difficult }\\
&\text{in terms of the order of oracle complexity?}
\end{empheq}

On the one hand, the answer appears to be negative for some non-convex NSCO problems. In the one-layer problem, we know it takes $\bigO(1/\ep^2)$ \cite{ghadimi2013stochastic} queries to the stochastic oracle to find an $\ep$-stationary solution $\xbar$, i.e., $\E[\normsq{\grad f(\xbar)}]\leq \ep$. Recently, by using a specially-designed potential function, Ghadimi et al. \cite{ghadimi2020single} developed an $\bigO(1/\ep^2)$ algorithm for the two-layer problem and extended it in \cite{balasubramanian2022stochastic} to solve the multi-layer problem with the same oracle complexity of $\bigO(1/\ep^2)$. Additionally, under some stronger smoothness assumptions, a variance reduction algorithm proposed in \cite{zhang2019multi} can improve the oracle complexity further to $\bigO(1/\ep^{1.5})$.

On the other hand, the answers is still unclear for convex NSCO problems. To find an $\ep$-optimal solution $\xbar$ for the one-layer problem, i.e., $\E[f(\xbar) - f(\xstar)]\leq \ep$, we know from \cite{nemirovsky1983problem} that the order-optimal oracle complexities are $\bigO(1/\ep^2)$ if $f$ is non-smooth and $\bigO(1/\ep)$ if $f$ is also strongly convex. However, the results for NSCO in the literature fail to match them. To the best of our knowledge, finite time convergence bounds for convex NSCO problems are available only when outer-layer functions are all smooth (see \cite{mengdi2017stochastic,Yang2019Multilevel,wang2017accelerating}). These works are based on gradient approximation; they use moving averages to track function values of the inner layers and apply multi-timescale schemes to ensure their faster convergence. However, as shown in Table \ref{tb:rates} and \ref{tb:krates}, these complexities appear to be sub-optimal. For the two-layer problem, if the innermost layer function is non-smooth, the oracle complexity is $\bigO(1/\ep^{4})$. With an  additional smoothness assumption for $f_2$, the oracle complexity can be improved to $\bigO(1/\ep^{3.5})$  and to $\bigO(1/\ep^{1.25})$ if the problem is also strongly convex. For the multi-layer problem, the complexities are even worse as they  depend exponentially on the number of layers $k$. Additionally, the assumption of smooth outer-layer functions is violated by some important applications of convex NSCO, for example, the risk averse two-stage stochastic program in Section 5.
\subsection{Our Contributions}\label{subsec:intro_our_c}

\begin{table}[tb]
\centering
\makebox[0pt][c]{\parbox{\textwidth}{%
    \begin{minipage}[b]{0.50\hsize}\centering
        \begin{threeparttable}
        \caption{Two-Layer Oracle Complexity}
        \label{tb:rates}
        \centering

        \begin{tabular}{|l l|  c |c|}

        \toprule
         Problem & Type  &
         In the Literature & SSD/SSDp    \\
        \midrule
         Convex  & \makecell{Outer Nonsmooth \\ Outer Smooth \\ All Smooth} &\makecell{N.A.  \\ $\bigO(1/\ep^{4})$ \cite{mengdi2017stochastic} \\ $\bigO(1/\ep^{3.5})$ \cite{wang2017accelerating}}
          & $O(1/\ep^2)$ \\

          \midrule
        \makecell{Strongly\\ Convex}  & \makecell{Outer Nonsmooth \\ Outer Smooth \\ All Smooth} & \makecell{N.A. \\ $\bigO(1/\ep^{1.5})$ \cite{mengdi2017stochastic} \\$\bigO(1/\ep^{1.25})$ \cite{wang2017accelerating} }
         & \makecell{$\bigO(1/\ep^2)$\\ $\bigO(1/\ep^{\ })$\\ $\bigO(1/\ep^{\ })$} \\

        \bottomrule
        \end{tabular}



        \end{threeparttable}
    \end{minipage}
    \hfill
    \begin{minipage}[b]{0.35\hsize}\centering
        \begin{threeparttable}
        \caption{k-Layer  }
        \label{tb:krates}
        \centering

        \begin{tabular}{| c | c|}

        \toprule
         
         In the Literature & SSD/SSDp   \\

        \midrule          
          \makecell{N.A. \\ $\bigO(1/\ep^{2k})$\cite{Yang2019Multilevel} \\ $\bigO(1/\ep^{(7+k)/4})$\cite{Yang2019Multilevel}} & $\bigO(1/\ep^2)$ \\
          \midrule        
          \makecell{N.A. \\ N.A. \\ $\bigO(1/\ep^{(3+k)/4})$\cite{Yang2019Multilevel}} & \makecell{$\bigO(1/\ep^2)$\\ $\bigO(1/\ep^{\ })$\\ $\bigO(1/\ep^{\ })$} \\
        \bottomrule
        \end{tabular}
        \end{threeparttable}
    \end{minipage}
 }}
\end{table}


In this paper, we answer the question by developing efficient methods to achieve order-optimal oracle complexities under some mild assumptions.
 We study the two-layer nested problem with a (possibly strongly) convex  regularization term $u(x)$ given by
\begin{equation}\label{eq:2_prob}
\min_{x \in X} \{f(x) := f_1 \circ f_2 (x) + u(x)\},
\end{equation}
where $X$ is a compact and convex set with a finite radius, i.e., $\max_{x, \xbar \in X} \norm{x - \xbar} \leq \Dx < \infty.$ 
We impose the following \textit{compositional convexity} assumption throughout the paper.
\begin{assumption}\label{ass:cp_cv}
A nested function $f_1 \circ f_2  (x)$ in \eqref{eq:2_prob} is said to satisfy the \textit{compositional convexity} assumption if 
\begin{itemize}
\item every layer function $f_i: \R^{n_i}\rightarrow \R^{n_{i-1}} $ is proper closed and convex. 
\item $\myfi[1]$ is  component-wise non-decreasing if $f_2$ is not affine.
\end{itemize}
\end{assumption}
The convexity assumption implies a bi-conjugate reformulation important for our development: 
\begin{equation}\label{eq:bi-con}
f_i(y_i) = \max_{\pii \in \Pii} \pii y_i - \fistar(\pii),
\end{equation}
 where $\fistar$ is the Fenchel conjugate of $f_i$ and $\Pii=\dom(\fistar):=\{\pii \in \R^{n_{i-1}\times n_i}: \fistar[i,j](\pii[i,j]) < \infty\ \forall j \in [n_{i-1}] \}$.
The above monotonicity and the layerwise convexity assumption together form the classic sufficient condition for a nested function to be convex (e.g. see \cite{nesterov2003introductory}). Compared to \cite{wang2017accelerating,mengdi2017stochastic,Yang2019Multilevel},  Assumption \ref{ass:cp_cv} is stronger than their assumption of the nested function $f$ being convex.  However, we do not require the outer-layer function $\myfi[1]$ to be smooth. \\

Rather than the gradient approximation approach in \cite{mengdi2017stochastic,wang2017accelerating,Yang2019Multilevel}, we take a reformulation approach.
 We obtain a saddle-point problem by replacing every layer function with the maximum over some linearizations indexed by (sub)gradients (e.g. see \eqref{eq:bi-con}), i.e.,  $\myfi[1]'$ and $\myfi[2]'$  play the roles of dual variables $\pii[1]$ and $\pii[2]$. This relaxes
the tight coupling between $x$ and the dual variables. 
Instead of requiring $\piit[1]\piit[2]$ to be the true gradient, we only need the iterative algorithm to generate $(\xt; \piit[1], \piit[2])$ to converge to some stationary point.
 So rather than $\piit[1]=\myfi[1]'(\myfi[2](\xt))$, we can select easier $\piit[1]$'s with readily available unbiased estimators such that the update of $\xt$ can  be performed with unbiased arguments. 
The availability of unbiased arguments is the key to our order-optimal methods. Specifically, our development can be summarized into four steps.

First, we consider a simple two-layer problem with $\myfi[1]$ and $\myfi[2]$ being smooth. Linearizing both layer functions using the bi-conjugate in \eqref{eq:bi-con}, we arrive at a ``$\min_x\max_{\pii[1]}\max_{\pii[2]}$" saddle-point problem similar to that in \cite{zhang2019efficient}. We extend the deterministic sequential dual (SD) method \cite{zhang2019efficient} to a stochastic sequential dual (SSD) method, each iteration of which consists of computing prox-mappings for $\pii[2]$ and $\pii[1]$ and then computing the prox-mapping for $x$.  
Choosing the prox-functions in a similar fashion as \cite{LanZhou18RPDG}, the prox-mappings for $\pii[1]$ and $\pii[2]$ simplify to gradient evaluations so that the SSD method can be implemented in a primal form.
  The SSD method achieves the oracle complexities of 
$\bigO(1/\ep^2)$ under the non-strongly convex setting  and $\bigO(1/\ep)$ under the strongly convex setting.  As illustrated in Table \ref{tb:rates}, they  improve over the best oracle complexities in the literature by orders of magnitude and 
 have the same orders as those for solving the simpler one-layer problem. 
  Moreover, in the deterministic setting, the gradient oracle complexities of the SSD method matches that of Nesterov accelerated gradient method \cite{nesterov2003introductory}. 
 We also consider a closely related problem where the layer function $\myfi[1]$ is structured non-smooth \cite{Nestrov2004Smooth}.
By selecting $\normsq{\cdot}$ as the prox-function when performing the $\pii[1]$-prox mapping, the SSD method achieves an oracle complexity of $\bigO(1/\ep^2)$. This complexity bound is clearly optimal for the non-strongly convex case, since it has the same order as the one for the simpler one-layer problem. Somewhat surprisingly, we demonstrate this bound is also tight for the strongly convex case, by presenting a new and matching lower complexity bound for this case.

Second, we consider a more complicated two-layer problem where both $\myfi[1]$ and $\myfi[2]$ are general non-smooth. 
Such a problem is particularly challenging due to the non-smooth outer-layer function $\myfi[1]$. From the (sub)gradient approximation perspective, estimating $\myfi[1]'(\myfi[2](x))$ is not possible because $\myfi[1]'$ may not be continuous with respect to its argument. From the SSD perspective, the efficient computation of the $\pii[1]$-prox mapping with any strongly convex prox-function is difficult. 
This motivates us to introduce a tri-conjugate reformulation to $\myfi[1]$, $\myfi[1](\yi[1])= \min_{\vi}\max_{\pii[1]} \inner{\pii[1]}{\yi[1] - \vi[1]} + \myfi[1](\vi),$ where $\vi$ is an auxiliary primal variable. In the tri-conjugate reformulation, the prox-mappings for both $\vi[1]$ and $\pii[1]$ can be computed efficiently with $\normsq{\cdot}$ as the prox-function. 
 So we propose a non-smooth stochastic sequential dual (nSSD) method, whose iteration consists of computing prox-mappings for $\pii[2]$ and $\pii[1]$ and then computing prox-mappings for $\vi[1]$ and $x$. The nSSD method achieves an oracle complexity of $\bigO(1/\ep^2)$. To the best of our knowledge, it is the first method to achieve any finite time convergence guarantee for the NSCO problem with a non-smooth outer-layer function. This complexity bound is optimal for the non-strongly convex case since it is in the same order of magnitude as the one for the simpler one-layer problem. Moreover, we present a lower complexity bound to show that it is not improvable for the strongly convex case as well.

Third, we extend the two-layer SSD and nSSD methods to the multi-layer setting where the $k$-layer problem is constructed from an arbitrary nested composition of smooth, structured non-smooth and general non-smooth layer functions. 
In order to construct unbiased arguments for the prox-mappings, we propose a novel repeated sampling scheme to ensure conditional independence among  estimators for the sequentially updated dual variables. 
As shown in Table \ref{tb:krates}, when all the outer-layer functions are smooth, the proposed multi-layer SSD method achieves the order-optimal oracle complexities of $\bigO(1/\ep^2)$ in the non-strongly convex case and of $\bigO(1/\ep)$ in the strongly convex case, improving over the exponential dependence of the order of oracle complexity on $k$ in \cite{Yang2019Multilevel}. When there is either a structured non-smooth or a general non-smooth outer-layer function, the multi-layer SSD/nSSD method achieves an order-optimal oracle complexity of $\bigO(1/\ep^2)$. Such a complexity result appears to be the first finite-time  convergence guarantee for the multi-layer NSCO problem with non-smooth outer-layer functions. Moreover, 
the  oracle complexity of $\bigO(1/\ep^2)$ can be attained with parameter-free stepsizes, which could be helpful in practice.

Fourth, we illustrate the SSD and nSSD methods by applying them to two interesting examples; minimizing the mean-upper-semideviation risk of order $1$ for a two-stage stochastic program and minimizing the maximum loss associated with a system of stochastic composite functions. A direct application of our methods leads to order-optimal oracle complexities.  We show that the constant dependence of  our methods can be further improved if some intra-layer problem structures are exploited. 

The rest of the paper is organized as follows. Section 2 introduces the SSD method for the smooth and structured-nonsmooth two-layer problem.  Section 3 introduces the nSSD method for the general non-smooth two-layer problem. Section 4 extends these methods to the multi-layer setting and Section 5 provides two concrete applications. Some  concluding remarks are made in Section 6.

\subsection{Notations\& Assumptions}\label{subsc:notation}

The following notations and assumptions are used throughout the paper. 

\begin{itemize}

\item The feasible region $X$ is convex and compact with $\Dx := \max_{x_1, x_2 \in X} \norm{x_1 - x_2}[2] < \infty$. We assume the solution set $X^*:= \argmin_{x \in X} f(x)$ to be non-empty and use $x^* \in X^*$ to denote an arbitrary optimal solution.

\item Every layer function $\myfi: \R^{n_i} \rightarrow \R^{n_{i-1}}$ is closed, convex and proper. We use the notation $\grad \myfi(\yi) \in \R^{n_{i-1}} \times \R^{n_i}$  ($\myfi'(\yi)  \in \R^{n_{i-1}} \times \R^{n_i}$) to denote the Jocabian (sub-gradient) matrix. 

\item There exists a stochastic first-order oracle $\SOi$ associated with every layer function $\myfi$.  When queried at some $\yi \in \Yi$, the $\SOi$ returns a pair of unbiased estimators $(\myfixi[i][\yi][], \myfixipr[i][\yi][])$ for $(\myfi(\yi), \myfipr(\yi))$.
Results returned by different queries to $\SOi$ are independent, and all $\SOi$'s are independent. 

 \item The Fenchel conjugate of a convex function $g(x)$ is defined as  $g^*(\pi) := \max_{x \in \R^n} \inner{x}{\pi} - g(x)$. The Bregman distance function (or prox-function) associated with a convex function $g$ is defined as $D_g(x; y) = g(x) - g(y) - \inner{g'(y)}{x -y}.$ For a m-dimensional vector valued function $g(x)$, its dual variable $\pi:=[\pii[1]; \pii[2];\ldots; \pii[m]]$ is an $m\times n$ matrix, and $g^*$ and $D_{g^*}$ are $m$-dimensional vector functions. Specifically, 
 the $j$\ts{th} entries of $g^*$ and $D_{g^*}$ are defined according to 
 $g^*_j(\pii[j]):= \max_{x \in \R^n} \inner{x}{\pii[j]} - g_j(x)$ and $D_{g^*,j}(\pi; \pibar) = g_j^*(\pii[j]) - g^*_j(\piibar[j]) - \inner{{g^*_j}'(\piibar[j])}{\pii[j] -\piibar[j]}.$

\item We use the term prox-mapping of $h$  to refer to the following type of computation:
\begin{equation}\label{def:prox}
y^t \leftarrow \targmin_{y \in Y} \inner{y}{g} + h(y) + \eta V(y; \bar{y}),\end{equation}
where $Y \subset \R^n$ is closed and convex, $h$ is a convex function, and $V(\cdot; \cdot)$ is some Bregman distance function. When $h$ is an $m$-dimensional vector function and $y \subset \R^{m \times n}$ is a matrix,  we use the term prox-mapping of $h$ to refer to computing the prox-mapping for every row of $y$, i.e., $y^t:= [y^t_1; \ldots; y^t_m]$ with
$$y_i^t\leftarrow \targmin_{y_i \in Y_i} \inner{y_i}{g} + h_i(y_i) + \eta V_i(y_i; \bar{y}_i) \ \forall i\in[m].$$
In both cases, we call $g$ the argument, $\bar y$ the prox center, $V$ the prox-function and $\eta$ the stepsize parameter. Moreover, we also use the term prox-mapping for $y$ to denote the prox-mapping of $h$ if the associated function $h(y)$ is clear from context, e.g., prox-mapping for $x$ refers to the prox-mapping of $u(x).$ 

\item $\norm{\cdot}$ denotes the $l_2$ (operator) norm unless specified otherwise. 
\item We use the notation $[\cdot]^+$ to denote projection onto the positive orthant, i.e.,  $[x]^+:=\max\{x, 0\}$, and the notation $[\cdot]_Y^+$ to denote the projection onto set $Y$.
\end{itemize}

    

\section{Smooth and Structured Non-smooth Two Layer Problems}

In this section, we present the SSD method for the two-layer problem in \eqref{eq:2_prob}.
 We assume $\myfi[2]$, the inner layer function, to be smooth and $\myfi[1]$, the outer layer function, to be either smooth or structured non-smooth. 
Specifically, for a function $g(y)$ defined on $Y$, we call it smooth if its gradient $\grad g(y)$ is Lipschitz continuous, i.e.,
$\norm{\grad g(y) - \grad g(\ybar)} \leq L \norm{ y- \ybar} \ \forall y, \ybar \in Y.$
We call it structured non-smooth if there is some convex closed set $\Pi$ and convex closed and proper function $g^*$ such that 
\begin{equation}\label{def:str-ns}
g(y) = \tmax_{\pi \in \Pi} \inner{\pi}{y} - g^*(\pi), \forall y \in Y. 
\end{equation}
We assume the prox-mapping of $g^*$ with $\norm{\cdot}[2]^2$ as the prox-function to be efficiently computable (see Subsection \ref{subsc:notation}). Notice such a definition differs from the one proposed by Nesterov \cite{Nestrov2004Smooth} where the inner product contains a general linear operator $A$, i.e., $\inner{\pi}{A y}$. However, our definition is not restrictive in NSCO because $A y$ can be regarded as the output from an inner linear layer function $A (\cdot)$. 

The coming subsections are organized as follows. Subsection \ref{sb:2-smo-alg} introduces the SSD method, followed by its convergence results in Subsection \ref{sb:2-smo-re}. Next, Subsection \ref{sb:2-smo-l} presents a lower complexity bound for the structured non-smooth problem and Subsection \ref{sb:2-smo-pf} presents the detailed convergence analysis.  

\subsection{The SSD Method}\label{sb:2-smo-alg}
As suggested in Subsection \ref{subsec:intro_our_c}, the development of the SSD method is inspired by a $\min-\max-\max$ reformulation of \eqref{eq:2_prob} given by 
\begin{equation}\label{eq:2_sad}
\min_{x \in X} \max_{\pii[1] \in \Pii[1]}\max_{\pii[2] \in \Pii[2]} \{\mytLa[x][\pii[1], \pii[2]]:= \nmyLa[x][\pii[1], \pii[2]][1] + u(x)\},
\end{equation}
where $\Pii[1]$ and $\Pii[2]$ are respective  domains of $\fistar[1]$ and $\fistar[2]$ (see \eqref{eq:bi-con}) and the compositional Lagrangian functions are defined according to  
\begin{align}\label{2eq:com_La}\begin{split}
\nmyLa[x][\pii[2]][2] = \pii[2] x - \fistar[2](\pii[2]), \text{ and }
\nmyLa[x][\pii[1], \pii[2]][1] = \pii[1] \nmyLa[x][\pii[2]][2] - \fistar[1](\pii[1]).
\end{split}
\end{align}
For simplicity, we will use the notation $z:=(x; \pii[1], \pii[2])$ and $Z:= X \times \Pii[1] \times \Pii[2]$ for the rest of the section. Since
$\La_1(x; \pii[1], \pii[2])$ can be interpreted as the nested composition of a lower linear approximation to $\myfi[2]$, specified by $\pii[2]$ and $\fistar[2](\pii[2])$, and   a lower linear approximation  to $\myfi[1]$, specified by $\pii[1]$ and $\fistar[1](\pii[1])$,  a certain duality relationship holds between $\La$ (c.f. \eqref{eq:2_sad}) and the original problem \eqref{eq:2_prob}.
\begin{lemma}\label{lm:duality}
Let $f$ and $\tilLa$ be defined in \eqref{eq:2_prob} and \eqref{eq:2_sad}. Then 
\begin{enumerate}
    \item[a)] Weak duality: $f(x) \geq \mytLa[x][\pii[1], \pii[2]]\ \ \forall (\pii[1], \pii[2]) \in \Pii[1] \times \Pii[2],\, \forall x \in X$. 
    \item[b)] Strong duality: for a given $x\in X$, $f(x) = \mytLa[x][\hpii[1], \hpii[2]]$ for any $\hpii[2] \in \partial \myfi[2](x)$ and $\hpii[1] \in \partial \myfi[1](\myfi[2](x))$.
    \item[c)] There exists a pair $(\piistar[1], \piistar[2])$ such that $\zstar:=(\xstar; \piistar[1], \piistar[2])$ is a saddle point, i.e., 
    $$\La(\xstar; \pii[1], \pii[2])\leq \La(\xstar; \piistar[1], \piistar[2]) \leq \La(x; \piistar[1], \piistar[2])\ \forall (x, \pii[1], \pii[2])\in Z.$$
    \item[d)] For any  $(x; \pii[1], \pii[2]) \in Z$, an upper bound on the optimality gap of $x$ is given by:
    \begin{equation}
    f(x) - f(\xstar) \leq \max_{\piibar[1] \in \Pii[1], \piibar[2] \in \Pii[2]} \La(x; \piibar[1], \piibar[2]) - \La(\xstar; \pii[1], \pii[2])
    \end{equation}
\end{enumerate}
\end{lemma}
\begin{proof}
First, for Part b), let $x \in X$ be given. It follows from $\hpii[1]\in \partial \myfi[1](\myfi[2](x))$ and $\hpii[2]\in \partial \myfi[2](x)$ that
\begin{align*}
 f_2(x) = \nmyLa[x][\hpii[2]][2], \text{ and }
 f_1(f_2(x)) = \nmyLa[x][\hpii[1], \hpii[2]][1].
\end{align*}
Thus $f(x) = \nmyLa[x][\hpii[1], \hpii[2]][1] + u(x)$.
Regarding Part a), for any $x \in X$, the following decomposition is valid for any $(\pii[1], \pii[2])\in \Pii[1] \times \Pii[2]$:
\begin{align*}
\mytLa[x][\hpii[1], \hpii[2]]- \mytLa[x][\pii[1], \pii[2]] &= \nmyLa[x][\hpii[1], \hpii[2]] - \nmyLa[x][\pii[1],\hpii[2]]+ \nmyLa[x][\pii[1], \hpii[2]] - \nmyLa[x][\pii[1], \pii[2]]\\
&= \underbrace{\nmyLa[x][\hpii[1], \hpii[2]][1] - \nmyLa[x][\pii[1], \hpii[2]][1]}_{A} +  \underbrace{\pii[1](\nmyLa[x][\hpii[2]][2] - \nmyLa[x][\pii[2]][2])}_{B},
\end{align*}
where $A \geq 0$ because $\hpii[1] \in \argmax_{\pii[1] \in  \Pii[1]} \{\pii[1] \myfi[2](x) - \fistar[1](\pii[1]) \equiv \myLa[x][\pii[1], \hpii[2]][1]\}$. If $\myfi[2]$ is affine, $\Pii[2]$ is a singleton set such that $B=0$. Otherwise  a non-negative $\Pii[1]$ (Assumption \ref{ass:cp_cv}) and $\myLa[x][\hpii[2]][2] - \myLa[x][\pii[2]][2]\ge 0$ imply $B\geq0$. Therefore $f(x) = \nmyLa[x][\hpii[1], \hpii[2]] \geq \nmyLa[x][\pii[1], \pii[2]]$. 

As for Part c), the first-order optimality condition implies that there exist some $\piistar[1] \in \partial \myfi[1] (\myfi[2](\xstar))$, $\piistar[2] \in \partial \myfi[2](\xstar)$ and $u' \in \partial u(\xstar)$ such that $\piistar[1]\piistar[2](x - \xstar) + \inner{u'}{x - \xstar} \geq 0\ \forall x\in X.$ Thus, with $D_{u} (x;\xstar):= u(x) - u(\xstar) - \inner{u'}{x-\xstar}\geq 0$, we get
$$\La(x; \piistar[1], \piistar[2]) - \La(\xstar; \piistar[1], \piistar[2]) = \piistar[1]\piistar[2](x - \xstar)+ \inner{u'}{x - \xstar} + D_{u} (x;\xstar) \geq 0.$$
The relation $\La(\xstar; \pii[1], \pii[2])\leq \La(\xstar; \piistar[1], \piistar[2])$ is a direct consequence of the strong duality. 

 Part d) is  a direct consequence of Parts a) and b).
\end{proof}
\vgap

Now we introduce a $Q$-gap function which is often-used for saddle point problems (e.g. see \cite{LanBook} and \cite{zhang2019efficient}).  
For a point $\zt:=(\xt; \piit[1], \piit[2]) \in Z$, the $Q$-gap function, defined with respect to some reference point $z  \in Z$, is given by 
\begin{equation}\label{def:gap_func}
Q(\zt, z):= \myLa[\pii[1],\pii[2]][\xt] - \myLa[\piit[1], \piit[2]][x].
\end{equation}
The $Q$-gap function plays a central role in our development for two reasons. First, it provides an upper bound to the optimality gap (see Lemma \ref{lm:duality}.d)), so minimizing the $Q$-gap leads to an optimal solution for the original problem. 
Second, the $Q$-gap function admits a decomposition conducive to algorithm design.
Specifically, the development of the SSD method is motivated by the following decomposition:
\begin{equation} \label{decompose_gap}
Q(\zt, z) = Q_2(\zt, z) +  Q_1(\zt, z) +  Q_0(\zt, z),
\end{equation}
\vspace{-8mm}
\begin{align}
Q_2(\zt, z) &:= \mytLa[\xt][\pii[1], \pii[2]] - \mytLa[\xt][\pii[1], \piit[2]] 
= \pii[1][\pii[2]\xt - \fistar[2](\pii[2]) ] {\boxed{-\pii[1][\piit[2] \xt - \fistar[2](\piit[2])]}}, \label{eq:def_Q2} \\
Q_1(\zt, z) &:= \mytLa[\xt][\pii[1], \piit[2]] - \mytLa[\xt][\piit[1], \piit[2]] 
= \pii[1] \myLa[\xt][\piit[2]][2] - \fistar[1](\pii[1])  {\boxed{-[ \piit[1]\myLa[\xt][\piit[2]][2] - \fistar[1](\piit[1])]}},\label{eq:def_Q1}\\
Q_0(\zt, z) &:= \mytLa[\xt][\piit[1], \piit[2]] - \mytLa[x][\piit[1], \piit[2]] 
= {\boxed{ \piit[1]\piit[2] \xt + u(\xt)}} - (\piit[1]\piit[2]x + u(x)), \label{eq:def_Q0}
\end{align}
where $Q_2, Q_1,$ and $ Q_0$ relate to the optimality of $\piit[2], \piit[1]$ and $\xt$, respectively. The conceptual sequential dual (SD) method, originally proposed in \cite{zhang2019efficient}, performs  prox-mappings for $\pii[2]$, $\pii[1]$ and  $x$ in order to reduce $Q_2$, $Q_1$, and $Q_0$, i.e., the boxed terms in \eqref{eq:def_Q2}, \eqref{eq:def_Q1} and \eqref{eq:def_Q0}. With $(\xt[0]; \piit[1][0], \piit[2][0]) \in Z $, the $t$\ts{th} iteration of the SD method is given by  
\begin{align}\label{ls:SD_idea}
\begin{split}
&\piit[2] \leftarrow \argmax_{\pii[2] \in \Pii[2]} \pii[2]\tilde{x}^t - \fistar[2](\pii[2]) - \tauit[2] U_2(\pii[2];\piitt[2]), \text{ where } \tilde{x}^t:=  x^{t-1} + \thetat (\xtt - \xt[t-2]);\\
&\piit[1] \leftarrow \argmax_{\pii[1] \in \Pii[1]} \pii[1]\tilde{y}_1^t - \fistar[1](\pii[1]) - \tauit[1] U_1(\pii[1];\piitt[1]), \text{ where}\ \tilde{y}_1^t:= \La_2(\xtt; \piit[2]) + \thetat \piitt[2](\xtt - \xt[t-2]);\\
&\xt \leftarrow \targmin_{x \in X}\  \yitilt[0]x + u(x)+\tfrac{\etat}{2} \normsq{x - \xtt},\ \text{where } \yitilt[0]:= \piit[1]\piit[2].
\end{split}
\end{align}
In the above listing, $U_2$ and $U_1$ denote  general prox-functions and scalars  $\tauit[2]$, $\tauit[1]$ and $\etat$ represent non-negative stepsizes. Since $\xt$ is yet available, the arguments to the $\piit[2]$ and $\piit[1]$ prox-mappings are  extrapolated points, i.e., $\tilde{x}^t$ and $\tilde{y}_1^t$ to predict  $\xt$ and $\La_2(\xt; \piit[2])$.

The deterministic SD method has been shown to achieve optimal complexities under various settings \cite{zhang2019efficient,lan2022optimal}. 
To adapt it to the NSCO setting, we propose to replace the deterministic arguments, $\yitilt[1]$ and $\yitilt[0],$  with some stochastic estimators, $\yitilt[1](\xi)$ and $\yitilt[0](\xi)$. This leads to a stochastic version of the SD method, namely, the SSD method. Initialized to $(\xt[0]; \piit[1][0], \piit[2][0])\in Z$ and $\xt[-1] = \xt[0]$, its $t$\ts{th} iteration  is given by 
 \begin{align}\label{ls:2-SSD-gen}
\begin{split}
&\piit[2] \leftarrow \targmax_{\pii[2] \in \Pii[2]}\ \pii[2]\tilde{x}^t - \fistar[2](\pii[2]) - \tauit[2] U_2(\pii[2];\piitt[2]), \text{ where } \tilde{x}^t:=  x^{t-1} + \thetat (\xtt - \xt[t-2]);\\
&\piit[1] \leftarrow \targmax_{\pii[1] \in \Pii[1]}\ \pii[1]\yitilt[1](\xi) - \fistar[1](\pii[1]) - \tauit[1] U_1(\pii[1];\piitt[1]), \text{ where}\ \tilde{y}_1^t:= \La_2(\xtt; \piit[2]) + \thetat \piitt[2](\xtt - \xt[t-2]);\\
&\xt \leftarrow \targmin_{x \in X}\ \yitilt[0](\xi)x + u(x) +\tfrac{\etat}{2} \normsq{x - \xtt}, \text{where } \yitilt[0]:= \piit[1]\piit[2].
\end{split}
\end{align}

Next we provide concrete implementation to the above prox-mappings and construction of the unbiased estimators, i.e.,  $\yitilt[1](\xi)$ and $\yitilt[0](\xi)$ for $\yitilt[1]$ and $\yitilt[0]$.
First, let us consider the two-layer smooth problem where $\myfi[1]$ is smooth. A key challenge appears to be that the conjugate functions, $\fistar$'s, are not explicitly known. However, if $\Dfistar$, the Bregman distance function generated by $\fistar$, is selected as the prox-function $U_i$, there exists some primal equivalences to the prox-mappings of $\fistar$'s such that the steps in \eqref{ls:2-SSD-gen} simplify to gradient evaluations. 
Towards that end, an important result is a  duality relation \cite{beck2017first} between a  closed convex and proper function $g$ and its Fenchel conjugate $g^*$,
\begin{equation}\label{rel:conj_duality}
\pii[] \in \targmax_{\piibar[] \in \dom(g^*)} \piibar[] \uyi[] - g^*(\piibar[]) \Longleftrightarrow \pii[] \in \partial g(\uyi[]) \Longleftrightarrow g^*(\pii[]) = \pii[]\uyi[] -  g(\uyi[]).
\end{equation}
The first relation in \eqref{rel:conj_duality} implies  the equivalence of the prox-mapping of $\fistar$ to a gradient evaluation at some averaged point. Specifically, the following lemma is an extension of a similar result in \cite{LanZhou18RPDG,LanBook} to vector-valued functions and its proof is given in the appendix.
\begin{lemma} \label{lm:proximal_gradient_lemma}
Given a convex, closed and proper vector function $g$, if $\piitt[]$ is associated with some primal point $\uyitt[]$, i.e., $\piitt[] = g'(\uyitt[]) \in \partial g(\uyitt[])$, then the prox-mapping for $\piit[]$, i.e., 
   $\piit[] \in \targmin_{\pi \in \dom(g^*)} -\pii[] \yitilt[] + g^*(\pii[]) + \taut D_{g^*}(\pi; \piitt[]),$ is equivalent to 
\begin{equation}\label{eq:implicit-prox}
\piit[] = g'(\uyit[]) \in \partial g(\uyit[]) \text{ with }\ \uyit[] := (\yitilt[] + \taut \uyitt[])/ (1 + \taut).
\end{equation}
\end{lemma}
Thus, initialized to $\piit[1][0] = \grad \myfi[1](\uyit[1][0])$ and $\piit[2][0] = \grad \myfi[2](\uyit[2][0])$, each iteration of the SSD algorithm in its primal form is given by 
\begin{align}\label{ls:SD_idea_primal}
\begin{split}
&\piit[2]\leftarrow \grad \myfi[2](\uyit[2]),  \text{ where }  \uyit[2] := (\xtilt + \tauit[2] \uyit[2][t-1])/(1 + \tauit[2]), \ \tilde{x}^t:=  x^{t-1} + \thetat (\xtt - \xt[t-2]);\\
&\tilde{y}_1^t:= \myfi[2](\uyit[2]) +  \grad \myfi[2](\uyit[2][t])(\xtt - \uyit[2]) + \thetat \grad \myfi[2](\uyit[2][t-1])(\xtt - \xt[t-2]);\\
&\piit[1] \leftarrow \grad \myfi[1](\uyit[1]),  \text{ where } \uyit[1] := (\tilde{y}_1^t(\xi) + \tauit[1] \uyit[2][t-1])/(1 + \tauit[1]);\\
&\xt \leftarrow \targmin_{x \in X} \yitilt[0](\xi)x + u(x)+\tfrac{\etat}{2} \normsq{x - \xtt} \text{ where } \yitilt[0]:= \grad \myfi[1](\uyit[1]) \grad \myfi[2](\uyit[2]).
\end{split}
\end{align}
Notice that the second line in \eqref{ls:SD_idea_primal} utilizes the relation 
$\myfi[2](\uyit[2]) +  \grad \myfi[2](\uyit[2][t])(\xtt - \uyit[2]) = \La_2(\xtt; \piit[2]),$ which holds because  the second equivalence in \eqref{rel:conj_duality} imply  
$\fistar[2](\piit[2]) = \grad \myfi[2](\uyit[2][t]) \uyit[2] - \myfi[2](\uyit[2])$ if $\piit[2] = \grad \myfi[2](\uyit[2][t])$ (first line of \eqref{ls:SD_idea_primal}).
\vspace{2mm}

The construction for $\yitilt[1](\xi)$ and $\yitilt[2](\xi)$ are provided in the concrete primal-form SSD method shown in Algorithm \ref{alg:ssd_2}. There are two observations worth mentioning.  
First, due to the primal equivalence to  $\yitilt[1]$ in Line 2 of \eqref{ls:SD_idea_primal}, an unbiased estimator to it, $\yitilt[1](\xi)$, can be constructed in Line 4 of Algorithm \ref{alg:ssd_2} with stochastic estimators from $\SOi[2]$, i.e., $f_2(\uyit[2], \xi_2^1)$, $ \grad \myfi[2](\uyit[2], \xi_2^1)$,$ \myfi[2](\uyitt[2], \hat{\xi}_2)$. It is interesting to note the estimator for $\fistar[2](\piit[2])$ (as a part of $\La_2(\xtt; \piit[2])$ in line 2 of \eqref{ls:2-SSD-gen}) is obtained without explicitly knowing either $\fistar[2]$ nor  $\piit[2]$.
 Second, the dependence between the sequentially generated iterates  requires us to make three independent queries to $\SOi[2]$ (e.g. Line 3 of Algorithm \ref{alg:ssd_2}).  For example, consider $\piit[1]\piit[2]$, the argument to the $\xt$-prox mapping in \eqref{ls:SD_idea_primal}. As shown in Figure \ref{fig:2_sampling}, since the estimator $\grad \myfi[2](\piit[2], \xi_2^1)$ is used to generate $\piit[1]$ in Line 4 of Algorithm \ref{alg:ssd_2}, $\E[\grad \myfi[2](\uyit[2], \xi_2^1)|\piit[1]] \neq  \grad \myfi[2](\uyit[2]).$ Line 5 needs a new estimator $\grad \myfi[2](\uyit[2], \xi_2^0)$ independent of $\piit[1]$ (conditioned on $\piit[2]$) to ensure 
  $$\E[\grad \myfi[1](\uyit[1], \xi_1^0) \grad \myfi[2](\uyit[2], \xi_2^0)|\piit[1], \piit[2]] = \grad \myfi[1](\uyit[1]) \grad \myfi[2](\uyit[2]) = \piit[1]\piit[2].$$
 An additional estimator $\myfi[2](\uyitt[2], \hat{\xi}_2)$ is also required in the next iteration for a similar reason. To highlight their independence, we use the notation $\xi_i^j$ to emphasize that an independent $j$\ts{th} estimator is drawn from $\SOi$ and that it is used as part of the argument for the prox-mapping to reduce $Q_j$ and use the notation $\hat{\xi}_i$ to emphasize it being used for momentum extrapolation in the next iteration. 

\begin{figure}
\centering
\tikz{
  \node[obs] (pitwo) {$\piit[2]$};
  \node[rectangle, draw,right=of pitwo, yshift=-1cm] (ftwo) {$\myfixi[2][\uyit[2]][1], \grad \myfi[2](\uyit[2], \xi_2^1)$};
  \node[rectangle, draw, right=of pitwo, yshift=1cm] (pixizero) {$\grad \myfi[2](\uyit[2], \xi_2^0)$};
  \node[obs, right=of ftwo] (piipone) {$\piit[1]$};
  \node[rectangle, draw, right=of piipone] (pixione) {$\grad \myfi[1](\uyit[1], \xi_1^0)$};
  \plate [inner sep=.25cm,yshift=.2cm] {plate1} {(pitwo)(ftwo)(pixizero)(piipone)(pixione)} {};
  \edge {pitwo} {ftwo, pixizero};
  \edge {ftwo} {piipone};
  \edge {piipone} {pixione};
}
\caption{Illustration of stochastic dependency: $\grad \myfi[2](\uyit[2], \xi_2^0)$ is independent of  $\grad \myfi[1](\uyit[1], \xi_1^0)$ conditioned on $\piit[2]= \grad \myfi[2](\uyit[2])$.
}\label{fig:2_sampling}
\end{figure}

\begin{algorithm}[htb]
\caption{Stochastic Sequential Dual (SSD) Method for the Smooth Two Layer Problem}
\label{alg:ssd_2}
\begin{algorithmic}[1]
\Require $x^0 \in X$.
\parState{Set $\uyit[2][0] = x^{-1} \leftarrow x^0$.
Call $\SOi[2]$ to obtain $\grad \myfi[2](\uyit[2][0], \hat{\xi}_2)$ and set $\uyit[1][0] \leftarrow \myfi[2](x^0, \hxi[2])$.}
\For{$t = 1,2,3 ... N$}
\vspace{1mm}
\parState{Set $\tilde{x}^t\leftarrow x^{t-1} + \thetat (\xtt - \xt[t-2])$. Set $\uyit[2] \leftarrow (\tauit[2] \uyitt[2] + \xtilt) /(1+ \tauit[2]).$ \vspace{1mm}\\
Call $\SOi[2]$ to obtain independent estimates  $ \{f_2(\uyit[2], \xi_2^1), \grad \myfi[2](\uyit[2], \xi_2^1), \grad \myfi[2](\uyit[2], \xi_2^0),  \grad \myfi[2](\uyit[2], \hat{\xi}_2)\}$. \vspace{1mm}
}
\parState{Set $\yitilt[1](\xi) \leftarrow  f_2(\uyit[2], \xi_2^1) + \grad \myfi[2](\uyit[2], \xi_2^1)(\xtt-\uyit[2]) + \thetat \grad \myfi[2](\uyitt[2], \hat{\xi}_2) (\xtt - \xttt).$ \vspace{1mm}\\
Set $\uyit[1]\leftarrow (\tauit[1] \uyit[1][t-1](\xi) + \yitilt[1]) / (1 + \tauit[1])$.
Call $\SOi[1]$ to obtain $\grad \myfi[1](\uyit[1], \xi_1^0)$.\vspace{1mm}
}
\parState{Set $\yitilt[0](\xi) \leftarrow \grad \myfi[1](\uyit[1], \xi_1^0)  \grad \myfi[2](\uyit[2], \xi_2^0)$. Compute $\xt \leftarrow \argmin_{x \in X} \yitilt[0](\xi) x + u(x) + \tfrac{\etat}{2} \norm{x -\xtt}^2$, .\vspace{1mm}
}
\EndFor
\State Return $\xbarT = \nsumt \wt \xt / (\sumwt)$.
\end{algorithmic}
\end{algorithm}

Now we present some modifications required to handle a structured non-smooth $\myfi[1]$ (c.f. \eqref{def:str-ns}), i.e., 
\begin{equation}\label{def:st-ns-f1}
 \myfi[1](y_1):= \max_{\pii[1] \in \Pii[1]} \inner{\pii[1]}{y_1} - \fistar[1](\pii[1]), \ \forall y_1 \in \R^{n_1}.\end{equation} 
 Recall that our definition of structured non-smoothness implies that the prox-mapping of $\fistar[1]$, with $\normsq{\cdot}$ as the prox-function, can be efficiently computed. So the concrete implementation of \eqref{ls:2-SSD-gen}, the SSD method for the two-layer problem with a structured non-smooth $\myfi[1]$, is the same as Algorithm \ref{alg:ssd_2} except for initializing $\piit[1][0]$ to some point in $\Pii[1]$ and replacing Line 4 and 5 by 
\begin{align}\label{ls:2-ssd-strns}\begin{split}
 &\text{Set } \yitilt[1](\xi) \leftarrow  f_2(\uyit[2], \xi_2^1) + \grad \myfi[2](\uyit[2], \xi_2^1)(\xtt-\uyit[2]) + \thetat \grad \myfi[2](\uyitt[2], \hat{\xi}_2) (\xtt - \xttt).\\
&\text{Set } \piit[1]\leftarrow \targmax_{\pii[1] \in \Pii[1]} \inner{\pii[1]}{\yitilt[1](\xi)} - \fistar[1](\pi) - \tauit[1]\normsq{\pii[1] - \piitt[1]}/2.\\
&\text{Set } \yitilt[0](\xi) \leftarrow \piit[1]  \grad \myfi[2](\uyit[2], \xi_2^0) \text{ and compute }  \xt \leftarrow \targmin_{x \in X} \ \yitilt[0](\xi) x + u(x) + \tfrac{\etat}{2} \norm{x -\xtt}^2.
\end{split}
 \end{align}

\subsection{Convergence Results}\label{sb:2-smo-re}
We present in this subsection the convergence guarantees of the proposed SSD method. The proofs are deferred to Subsection \ref{sb:2-smo-pf}.

First, we need to specify a few problem parameters. We use the following sets, $Y_1$ and $Y_2$, to define the effective domains for $\myfi[1]$ and $\myfi[2]$. 
\begin{equation}\label{def:domain}
Y_2 := X, Y_1:=\text{Conv}\left(\{\La_2(\pii[2]; x)\ |\ \pii[2] = \grad \myfi[2](\uyi[2]),\ \uyi[2], x \in Y_2\}\right).\end{equation}
Notice that $Y_1$ is bounded set  because $X$ is assumed to be bounded.
For the layer function $\myfi[2]$, the definitions of the smoothness constant $\Lfi[2]$, the Lipschitz-continuity constant $\Mi[2]$ (over $Y_2$) and variance constants, $\sig[2]$ and $\sig[\myfi[2]]$, are listed below. 
\begin{align}\label{cst:f2}
\begin{split}
&\norm{\grad \myfi[2](y_2) - \grad \myfi[2](\yibar[2])} \leq \Lfi[2] \norm{y_2 - \yibar[2]}, \forall y_2, \yibar[2] \in \R^n.\\
 &\norm{\grad \myfi[2](y_2) } \leq M_2,\  \forall y_2 \in Y_2.  \\
&\E[\normsq{\grad \myfi[2](\uyi[2]) - \grad \myfi[2](\uyi[2], \xi_2)}] \leq \sig[2]^2,\ \E[\normsq{ \myfi[2](\uyi[2]) - \myfi[2](\uyi[2], \xi_2)}] \leq \sig[\myfi[2]]^2 \leq\sig[2]^2 \Dx^2, \forall \uyi[2] \in Y_2.
 \end{split}\end{align}
 Notice that we assume the variance of  $\myfi[2](x, \xi)$ to  be bounded also by $\sig[2]^2 \Dx^2$ to simplify the notation.
For the layer function $\myfi[1]$, we require the Lipschitz-smoothness $\Lfi[1]$ and the variance constant $\sig[1]$ to be defined over $\R^{m_1}$ because the stochastic estimate $\yitilt[1](\xi)$ could deviate from $Y_1$. However, the Lipschitz-continuity constant $\Mi[1]$ is defined only over the bounded effective domain $\Yi[1]$ (c.f. \eqref{def:domain}). This ensures that it  always remains finite. Specifically, their definitions are given by
\begin{align}\label{cst:f1}\begin{split}
&\norm{\grad \myfi[1](y_1) - \grad \myfi[1](\yibar[1])} \leq \Lfi[1]\norm{y_1 - \yibar[1]},\  \E[\normsq{\grad \myfi[1](\uyi[1]) - \grad \myfi[1](\uyi[1], \xi_1)}] \leq \sig[1]^2,\ \forall y_1, \yibar[1] \in \R^{n_1}. \ \  \\
 &\norm{\grad \myfi[1](\uyi[1]) } \leq M_1,\  \forall \uyi[1] \in \Yi[1].
\end{split}\end{align}.

Additionally, it is also useful to define uniform upper bounds for variances associated with the stochastic arguments to the prox-mappings in Algorithm \ref{alg:ssd_2}.
\begin{equation}\label{cst:agg_var}
\sigtl[1] := \{\max_{t \geq 0} \E[\normsq{\yitilt[1](\xi) - \E[\yitilt[1](\xi)] }]\}^{1/2},
\ \sigtl[x] := \{\max_{t \geq 0} \E[\normsq{\yitilt[0](\xi) - \E[\yitilt[0](\xi)] }] ]\}^{1/2}.\end{equation}
Now we are ready to state the convergence rate for the non-strongly convex problem. 
\begin{theorem}\label{thm:2-sm}
Let a smooth two-layer function $f$ be given, with its problem parameters defined in \eqref{cst:f2}, \eqref{cst:f1} and \eqref{cst:agg_var}. If $\{\xt\}$ is generated by Algorithm \ref{alg:ssd_2} with
\begin{equation}\label{stp:2-smo-non}
\wt = t,\ \thetat = (t-1)/t,\ \tauit[1]=\tauit[2]=(t-1)/2,\end{equation}
then $\sigtl[1]^2 \leq 3 \sig[2]^2 \Dx^2 + 2 \sig[\myfi[2]]^2 \leq 5 \sig[2]^2 \Dx^2$ and $\sigtl[x]^2 \leq \sig[1]^2 \Mi[2]^2 + 2\Mi[1]^2\sig[2]^2 + \sig[2]^2(\sig[1]^2 + 10 \Lfi[1]^2 \sig[2]^2 \Dx^2).$ If, in addition, $\etat$ is chosen as 
\begin{equation}\label{stp:2-smo-non-p}
\etat = \max \{\tfrac{2}{t+1} (\Mibar[1]\Lfi[2] + \Lfi[1]\Mi[2]^2 ), \tfrac{\sqrt{t} \sigtl[x]}{\DX}\},
\end{equation}
for some $\Mibar[1] \geq \norm{\grad \myfi[1](\myfi[2]( \xbarT))},$ the ergodic average solution $\xbarT$ satisfies 
$$\E[f(\xbarT) - f(\xstar)] \leq \tfrac{\Mibar[1]\Lfi[2] + \Lfi[1]\Mi[2]^2}{N(N+1)} \normsq{\xt[0] - \xstar} + \tfrac{4 \Lfi[1] \sigtl[1]^2}{N} + \tfrac{2 \Mi[1] \sigtl[1]}{\sqrt{N}} + \tfrac{4\sigtl[x]\Dx}{\sqrt{N}}, \forall N \geq 2.$$
\end{theorem}

We make two remarks regarding the result. First, in the deterministic case, i.e., $\sig[1] = \sig[2] =0$, the SSD method has an oracle complexity of $\bigO( (\Mibar[1]\Lfi[2] + \Lfi[1]\Mi[2]^2)^{1/2} \norm{\xt[0] - \xstar}/\sqrt{\ep}).$ Since $\myfi[1] \circ \myfi[2]$ restricted to $X$ has a smoothness constant of $\Lfi[1]\Mi[2]^2 + \Mi[1]\Lfi[2]$  and $\Mibar[1] \leq \Mi[1]$,\vspace{1mm} the oracle complexity of the SSD method is of the same order as the optimal  oracle complexity of $\bigO( (\Mi[1]\Lfi[2] + \Lfi[1]\Mi[2]^2)^{1/2} \norm{\xt[0] - \xstar}/\sqrt{\ep})$ for solving a smooth $f$, e.g. obtained by using the Nesterov's accelerated gradient method. In practice, if there exists some prior knowledge about $\norm{\grad \myfi[1](\myfi[2]( \xbarT))}$, say the output solution $\xbarT$ is inside some ball around $\xstar$, then we can select $\Mibar[1]$ to be significantly smaller than $\Mi[1]$ such that the SSD method outperforms the Nesterov's accelerated gradient method. This improvement is made possible by exploiting the nested problem structure. Second, the stepsize parameters of prox-mappings for the dual variables in \eqref{stp:2-smo-non} are parameter independent.  
The specific of choice of $\etat$ in \eqref{stp:2-smo-non-p} is only required to achieve the desired constant dependence. In practice, any $\etat$ that scales $\Theta(\sqrt{t})$  leads to the order-optimal  stochastic oracle complexity of $\bigO(1/\ep^2)$.

The next theorem presents the convergence rate for the strongly convex problem. 
\begin{theorem}\label{thm:2-sm-str}
Let a smooth two-layer function $f$ (c.f. \eqref{eq:2_prob}) be given and let the regularization term $u(x)$ have a positive strong convexity modulus of $\alpha$. Let problem parameters be defined in \eqref{cst:f2}, \eqref{cst:f1} and \eqref{cst:agg_var}. If $\{\xt\}$ is generated by Algorithm \ref{alg:ssd_2} with
\begin{equation*}
\wt = t,\ \thetat = (t-1)/t,\ \tauit[1]=\tauit[2]=(t-1)/2,\end{equation*}
we have $\sigtl[1]^2 \leq 3 \sig[2]^2 \Dx^2 + 2 \sig[\myfi[2]]^2 \leq 5 \sig[2]^2 \Dx^2$ and $\sigtl[x]^2 \leq \sig[1]^2 \Mi[2]^2 + 2\Mi[1]^2\sig[2]^2 + \sig[2]^2(\sig[1]^2 + 10 \Lfi[1]^2 \sig[2]^2 \Dx^2).$ If, in addition, $\etat$ is chosen as 
$$\etat := \max \{\tfrac{2}{t+1} (\Mibar[1]\Lfi[2] + \Lfi[1]\Mi[2]^2 ), \tfrac{(t-1)\alpha}{2}\},$$
for some $\Mibar[1] \geq \norm{\grad \myfi[1](\myfi[2]( \xbarT))},$ the ergodic average solution $\xbarT$ and the last iterate $\xT$ satisfy 
\begin{align}
\E[f(\xbarT) - f(\xstar)] &\leq (\tfrac{\Lfi[1]\Mi[2]^2}{\alpha} + 1)[\tfrac{\log(N + 1)\tilde L}{N^2} \normsq{\xt[0] - \xstar} + \tfrac{4}{N}( \Lfi[1] \sigtl[1]^2 + \tfrac{\sigtl[x]^2}{\alpha})], \label{eq:2-str-sm-conv}\\
\E[\normsq{\xt[N] - \xstar}] &\leq \tfrac{\tilde{L}}{\alpha N(N+1)} \normsq{\xt[0] -\xstar} + \tfrac{4}{\alpha N}( \Lfi[1] \sigtl[1]^2 + \tfrac{\sigtl[x]^2}{\alpha}),
\end{align}
where $\tilde{L}:=\Mibar[1]\Lfi[2] + \Lfi[1]\Mi[2]^2$ denotes the overall smoothness constant.
\end{theorem}

A few remarks are in order. First, a simpler stepsize choice of $\etat = (t-1)\alpha/2$ still leads to an order-optimal stochastic oracle complexity of $\bigO(1/\ep),$ improving over the $\bigO(1/\ep^2)$ complexity  of the non-strongly problem. Second, by applying a certain restarting technique to Algorithm \ref{alg:ssd_2} (see Section 4.2.3 of \cite{LanBook}), the stochastic oracle complexity for finding an $\ep$-close $\xt[N]$, i.e., $\E[\normsq{\xt[N] - \xstar}] \leq \ep$, can be improved further to 
$$\bigO[\sqrt{\tilde{L}/\alpha} \log(\normsq{x^0 -\xstar}/\ep) + \tfrac{4}{\alpha \ep}( \Lfi[1] \sigtl[1]^2 + \tfrac{\sigtl[x]^2}{\alpha})].$$
Notice that the deterministic part, $\bigO[\sqrt{\tilde{L}/\alpha} \log(\normsq{x^0 -\xstar}/\ep)]$, is not improvable and the stochastic part has optimal dependence with respect to both $\alpha$ and $\epsilon$. Third, the stochastic oracle complexity for finding an $\ep$-optimal solution, i.e., $\E[f(\xbarT) - f(\xstar)] \leq \ep$, is also $\bigO(1/\ep\alpha^2)$. Its dependence on $\alpha$ is worse than the $\bigO(1/\ep\alpha)$ complexity  for solving a one-layer stochastic optimization problem. It results from the analysis technique of combining the $Q$-gap function and the distance $\norm{\xbarT - \xstar}^2$ together to derive a function-value bound. Such a worse dependence on $\alpha$ is also observed in \cite{lan2018rgem} where a similar technique is used for analyzing a randomized algorithm. However, it is unclear if the dependence on $\alpha$ is improvable for the NSCO problem.

Now we move on to consider the case with a structured non-smooth $\myfi[1]$. We use $\Mitl[1]$ to denote its Lipschitz continuity constant, i.e.,
\begin{equation}\label{cst:str-ns-f1}
\max_{\pii[1] \in \Pii[1]} \norm{\pii[1]} \leq \Mitl[1].\end{equation}
Notice that the above definition implies that $\myfi[1]$ is $\Mitl[1]$-Lipschitz continuous over $\R^{n_1}$, rather than over the bounded effective domain $\Yi[1]$ in \eqref{cst:f1}. The convergence guarantee for the non-strongly case is presented in the next theorem. 
\begin{theorem}\label{thm:2-strns}
  Let a two-layer function $f$ (c.f. \eqref{eq:2_prob}) comprised of  a structured non-smooth $\myfi[1]$ and  a smooth $\myfi[2]$ be given. Let the problem parameters be defined in \eqref{cst:f2}, \eqref{cst:agg_var} and \eqref{cst:str-ns-f1}. If $\{\xt\}$ is generated according to \eqref{ls:2-ssd-strns} with
\begin{equation*}
\wt = t,\ \thetat = (t-1)/t,\ \tauit[2]=(t-1)/2,\end{equation*}
we have $\sigtl[1]^2 \leq 3 \sig[2]^2 \Dx^2 + 2 \sig[\myfi[2]]^2 \leq 5 \sig[2]^2 \Dx^2$ and $\sigtl[x]^2 \leq \Mi[1]^2\sig[2]^2 .$ If, in addition, the remaining prox-mapping stepsizes are chosen as 
\begin{equation}\label{stp:2-strns}
\etat := \max \{\tfrac{2\Mi[1]\Lfi[2]}{t+1}  + \tfrac{2 \Mi[1] \Mi[2]}{\Dx}, \tfrac{\sqrt{t} \sigtl[x]}{\Dx}\}, \tauit[1] := \max\{\tfrac{\Mi[2]\Dx}{2 \Mi[1]}, \tfrac{\sqrt{t} \sigtl[1]}{2\Mi[1]}\},\end{equation}
 the ergodic average solution $\xbarT$ satisfies  
\begin{align}\label{eq:2-strns-conv}
\E[f(\xbarT) - f(\xstar)] \leq \tfrac{\Mitl[1] L_2}{N^2} \normsq{\xt[0] - \xstar} + \tfrac{4\Mitl[1] \Mi[2]  \Dx}{N} + \tfrac{4 \Mitl[1] \sigtl[1]}{\sqrt{N}} + \tfrac{2 \Dx \sigtl[x]}{\sqrt{N}}.
\end{align}
\end{theorem}  

We make three remarks regarding the result. First, in the deterministic case with $\sig[2] =0$, the oracle complexity simplifies to 
$\bigO((\Mitl[1]\Lfi[2])^{1/2} \norm{\xt[0] - \xstar}/ \sqrt{\ep} + \Dx \Mitl[1]\Mi[2]/\ep),$
where the first and second term can be attributed to the smooth  $\myfi[2]$ and the structure non-smooth  $\myfi[1]$, respectively. Second, the specific choices of $\etat$ and $\tauit[1]$ in \eqref{stp:2-strns} is required only for the desired constant dependence. In the stochastic case, any $\etat=\Theta(\sqrt{t})$ and $\tauit[1]=\Theta(\sqrt{t})$ would lead to the order-optimal oracle complexity of $\bigO(1/\ep^2)$. Third, the $\bigO(1/\ep^2)$ complexity is not improvable even under the strongly convex setting, which we will discuss in the next subsection.

\subsection{Lower Complexity Bound}\label{sb:2-smo-l}
We present a lower complexity result for the strongly convex problem with a structured non-smooth outer-layer function. 
Specifically, we develop a lower bound on the number of queries to $\SOi[2]$ required to obtain an (expected) $\ep$-optimal solution for a class of first-order methods.
For simplicity, we assume $X$ to be a ball centered on $0$, $B(0, r)$, and $u(x):= \alpha \normsq{x}/2$. 

We take an approach similar to Nesterov \cite{nesterov2003introductory} by proposing an abstract computation scheme which updates on affine subspaces reachable by $x$, $y_1$ and $\pii[1]$.  
The abstract scheme consists of the following steps. In the beginning, provided with $\xt[0] \in X$, $\yit[1][0] \in \R^{n_1}$ and $\piit[1][0] \in \Pii[1]$, the affine sub-spaces are initialized to $\Xt[0]:= \sp(\xt[0])$, $\Yit[0]:=\{\yit[1][0]\}$, and $\Piit[0] = \sp(\piit[1][0])$. In the $t$\ts{th} iteration, the updates are given by (the ``+''  in the definitions of $\Xt[t]$, $\Yit[t]$ and $\Pii[1]^t$ represents the Minkowski sum. )
\emph{
\begin{align}\label{alg:abs-strns}
\begin{split}
&\text{Query}\ \SOi[2] \text{ to obtain } (\myfi[2](y_2^t, \xi_2^t), \myfi[2]'(y_2^t, \xi_2^t)) \text{ at some } y_2^t \in \Xt[t-1].\\
&\Yit[t] := \Yit[t-1]+\sp\{\myfi[2](y_2^t, \xi_2^t) - y_1^0\} + \{ \myfi[2]'(y_2^j, \xi_2^j) x: j\leq t, x \in \Xt[t-1]\} .\\
&\Piit[t] := \Piit[t-1] + \sp(\piit[1]) \text{ where } \piit[1] = \targmax_{\pii[1] \in \Pii[1]} \inner{\pii[1]}{y_1^t} - \fistar[1](\pii[1]) - \tauit[1] \normsq{\pii[1] - \piitt[1]}/2, \\
&\quad \quad \text{ for some } y_1^t \in \Yit[t], \tauit[1]\geq 0, \piitt[1] \in \Piit[t-1].\\
& \Xt[t] :=  \Xt[t-1]+\{ \myfi[2]'(y_2^j, \xi_2^j)^\top \pii[1]^\top: j \leq t, \pii[1] \in \Pii[1]^t\} . 
\end{split}
\end{align}}
After $N$ iterations, the scheme can output any $x^N \in \Xt[N]$.

The abstract scheme is quite general. Since $X$ is a zero-centered ball, $\Xt[t]$ includes all possible prox-mappings from every prox-center in $\Xt[t-1]$, with every non-negative stepsize $\etat$ and with every possible argument in $\{\pii[1] \myfi[2]'(y_2^j, \xi_2^j)\}_{j \leq t, \pii[1]\in \Pii[1]^t}$. In particular, since 
$\Yit[t]$ contains every possible convex combination of evaluated function values $\{\myfi[2](y_2^j, \xi_2^j)\}$ and the prox-mapping for $\piit[1]$ simplifies to gradient evaluation if $\tauit[1] =0$, the argument set $\{\pii[1] \myfi[2]'(y_2^j, \xi_2^j)\}_{j \leq t, \pii[1]\in \Pii[1]^t}$ covers the pseudo-gradient in gradient approximation type algorithm, i.e., the SCGD algorithm in \cite{mengdi2017stochastic} is a special case. Moreover, since $\Yit[t]$ also contains all possible momentum-extrapolation term in $\{ \myfi[2]'(y_2^j, \xi_2^j) x: j\leq t, x \in \Xt[t-1]\}$, our SSD method \eqref{ls:2-ssd-strns} is also a special case. It is worth noting that one iteration of the SSD method corresponds to three iterations of the abstract scheme since it requires three independent queries to $\SOi[2]$ at $\uyit[2]$. The next theorem states the lower complexity result and its proof is deferred to the Appendix.

\begin{theorem}\label{thm:l-strns}
Given problem parameters $\ep >0$, $\Mitl[1] > 0$, $\bar{\alpha}\leq \Mi[1]^2/(4\ep)$, and $\sig[\myfi[2]]\geq 4\ep / \Mi[1]$,  there exists a  nested two-layer problem \eqref{eq:2_prob} consisted of a structured non-smooth $\myfi[1]$ and a stochastic linear $\myfi[2]$ such that $\myfi[1]$ is $\Mitl[1]$-Lipschitz continuous (c.f. \eqref{cst:str-ns-f1}), the variance of $\myfi[2](x, \xi)$ is bounded by $\sig[\myfi[2]]$ (c.f. \eqref{cst:f1}), and $u(x) = \bar{\alpha} \normsq{x}/2$. If the abstract scheme in \eqref{alg:abs-strns} is initialized to $\xt[0]=0$, $\piit[1][0]=0$ and some $\yit[1][0]$ and its output solution $x^N$ satisfies $\E[f(x^N) - f(\xstar)]\leq \ep$, then $N\geq \Omega(\Mitl[1]^2\sig[\myfi[2]]^2/\ep^2).$ 
\end{theorem}

The preceding lower bound shows the $\bigO(1/\ep^2)$ oracle complexity in \eqref{eq:2-strns-conv} to be order-optimal even under the strongly convex setting, implying that the structured non-smooth outer-layer function makes a composite nested problem intrinsically more difficult than the one-layer problem. 

\subsection{Convergence Proofs}\label{sb:2-smo-pf}
 We present in this subsection the detailed convergence proofs for the results in Subsection \ref{sb:2-smo-re}. We begin by providing a convergence bound of the $Q$-gap function which is valid for Theorem \ref{thm:2-sm}, \ref{thm:2-sm-str} and \ref{thm:2-strns}.

 \begin{proposition}\label{pr:2-sm}
Consider a two-layer problem of the form \eqref{eq:2_prob}, with solution iterates $\{\zt:=(\xt; \piit[1], \piit[2])\}$ generated according to \eqref{ls:2-SSD-gen} with $\Ui[2] = \Dfistar[2]$ and $\piit[2][0] = \grad \myfi[2](\xt[0])$. Suppose $\Ui[1]$ is $\beta$-strongly convex and $\fistar[1]$ is $\mu$-strongly convex with respect to $\Ui[1]$, i.e., 
 $$\Ui[1](\pii[1]; \piibar[1]) \geq \beta \normsq{\pii[1] - \piibar[1]}/2,\ \Dfistar[1](\pii[1]; \piibar[1]) \geq \mu \Ui[1](\pii[1]; \piibar )  \ \forall \pii[1], \piibar[1] \in \Pii[1].$$
 Suppose the Lipschitz-continuity and smoothness constants of $\myfi[2]$ are defined in \eqref{cst:f2} and the stochastic arguments $\yitilt[1](\xi)$ and $\yitilt[0](\xi)$ satisfy the variance bounds in \eqref{cst:agg_var}. 
 Let $z:=(\xstar; \pii[1], \pii[2] = \grad \myfi[2](y_2))$, for some $ y_2 \in X$, denote a reference point which could potentially depends on $\{\zt\}$ and assume the following requirements are satisfied for all $t \geq 1$ with some non-negative weight $\wt$:
 \begin{align}
  \begin{split}
  &\wt[t] = \thetat[t+1]\wt[t+1],\ \wt(\tauit[2] + 1) \geq \wt[t+1] \tauit[2][t+1], \\
&\etat \geq \tfrac{\thetat[t+1] \norm{\pii[1]} \Li[2]}{\tauit[2][t+1]} + \tfrac{\thetat[t+1]  \Mitil[2]^2}{ \tauit[1][t+1] \beta},\  
 \etat[N] \geq \tfrac{ \norm{\pii[1]} \Li[2]}{\tauit[2][N] + 1} + \tfrac{  \Mitil[2]^2}{ (\tauit[1][N] + \mu)  \beta},
  \end{split}
  \end{align} 
  where the constant $\Mitil[2]\geq \norm{\piit[2]}\ \forall t$. 
  Then following $Q$-gap bound is valid with $\wt[0]=0$
  \begin{align}\label{eq:2-sm-Q-conv}
  \begin{split}
  \E[\sumwt Q(\zt; z) &+ {\wt[N](\etat[N] + \alpha)} \normsq{\xt[N] - \xstar}/2] \leq \E[\sumt \{\wt\etat - \wt[t-1](\etat[t-1] + \alpha)\} \normsq{\xt[t-1] - \xstar}/2]\\
  &+ \E[\sumt \{\wt\tauit[1] - \wt[t-1](\tauit[1][t-1] + \mu)\} \Ui[1](\pii[1];\piit[1][t-1])] + \wt[1]\tauit[2][1] \Dfistar[2](\pii[2]; \piit[2][0])\\
  &+ \sumt \wt \sigtl[x]^2 /(\etat + \alpha) + \sumwt \sigtl[1]^2 /\beta (\tauit[1] + \mu) + \E[\sumt \wt \inner{\pii[1]}{\yitilt[1]- \yitilt[1](\xi)}].
  \end{split}
  \end{align}

 \end{proposition}
 \begin{proof}
 First, let us  use the definition of  prox-mapping for $\pii[1]$ in \eqref{ls:2-SSD-gen} to derive a convergence bound for $Q_1$ (c.f. \eqref{eq:def_Q1}).
 The $\mu-$strong convexity of $\fistar[1]$ with respect to $\Ui[1]$ implies a three point inequality (e.g. Lemma 3.8 of \cite{LanBook})
 \begin{equation}\label{pf:3pointpi1}
 - (\piit[1] -\pii[1])   \yitilt[1](\xi) + \fistar[1](\pii[1]) - \fistar[1](\piit[1]) + (\tauit[1] + \mu) \Ui[1](\pii[1];\piit[1]) + \tauit[1] \Ui[1](\piit[1];\piitt[1]) - \tauit[1] \Ui[1](\pii[1];\piitt[1]) \leq 0.\end{equation}

Let's focus on $- (\piit[1] -\pii[1])   \yitilt[1](\xi).$ Setting $\delta_1^t:=  \yitilt[1] - \yitilt[1](\xi)$, we get 
\begin{align*}
\E[- (\piit[1] -\pii[1])   \yitilt[1](\xi)] = \E[-(\piit[1] - \pii[1])\myLa[\xt][\piit[2]][2] + 
(\piit[1] - \pii[1])\delta_1^t + (\piit[1] - \pii[1])(\myLa[\xt][\piit[2]][2] - \yitilt[1])], 
\end{align*}
whereby 
\begin{align*}(\piit[1] - \pii[1])(\piit[2]\xt - \yitilt[1]) &= (\piit[1] - \pii[1])\piit[2] (\xt -\xtt) \vspace{-5mm}
- \thetat (\piitt[1] - \pii[1])\piitt[2] (\xtt -x^{t-2}) \\
&+ \thetat (\piitt[1] - \piit[1])\piitt[2] (\xtt -x^{t-2}),\\
\E[(\piit[1] - \pii[1])\delta_1^t ] 
&= \E[\piit[1]\delta_1^t ] -\E [\pii[1]\delta_1^t].
\end{align*}
Since $\E[\delta_1^t | \yitilt[1]] = 0$ and $\fistar[1](\pii[1]) + \tauit[1]\Ui(\pii[1]; \piit[1][t-1])$ has a strong convexity modulus of $(\mu + \tauit[1])\beta$, Lemma \ref{lm:prox-noise-bd} from the Appendix implies that 
$$-\E[\piit[1]\delta_1^t] \leq \sigtl[1]^2/\beta (\mu + \tauit[1]).$$
Moreover, 
\begin{align*}
&\E[ - \thetat (\piitt[1] - \piit[1])\piitt[2] (\xtt -x^{t-2}) - \tauit[1] \Ui[1](\piit[1]; \piitt[1]) \leq  \thetat^2 L_1\Mitil[2]^2 \normsq{\xtt - \xttt}/\beta\tauit[1], \\
&\E[- (\piiT[1] - \pii[1])\pii[2]^{N} (\xT -x^{N-1}) - (\tauT + \alpha) \Ui[1](\pii[1]^{N}; \pii[1]^{N-1})]  \leq  L_1 \Mitil[2]^2 \normsq{\xT - x^{N-1}}/\beta(\tauT + \mu),\\
&\E[ - \thetat (\piiO[1] - \pii[1]^1)\piiO[2] (x^{0} -x^{-1})] = 0. 
\end{align*}
Thus, the $\wt$-weighted sum of \eqref{pf:3pointpi1} satisfies  
\begin{equation}\label{pfpr1:Q1_bound}
\begin{split}
\E[\nsumt \wt Q_1(\zt; z)] \leq& \sumt [\wt\tauit[1] - \wt[t-1](\tauit[1][t-1] + \mu)]\Ui[1](\pii[1]; \piit[1][t-1]) + \sumwt \sigtl[1]^2 /[\beta (\tauit[1] + \mu)] \\
&+ \E [\nsumt \pii[1]\delta_1^t] +  \tsum_{t=2}^N \tfrac{\wtt \thetat \Mitil[2]^2}{\beta\tauit[1]}  \normsq{\xtt - \xttt} + \tfrac{\wT\Mitil[2]^2}{\beta(\tauit[1][N] + \mu)} \normsq{\xT - x^{N-1}}.
\end{split}
\end{equation}

Next, we consider the prox-mapping for $\pii[2]$. The definition of $\piit[2]$ in \eqref{ls:2-SSD-gen}, together with $\fistar[2]$ being $1$-strongly convex with respect to $\Dfistar[2]$, implies the following three point inequality
  \begin{equation*}
 - (\piit[2] -\pii[2])   \xtilt[t] + \fistar[2](\pii[2]) - \fistar[2](\piit[2]) + (\tauit[2] + 1) \Dfistar[2](\pii[2];\piit[2]) + \tauit[2] \Dfistar[2](\piit[2];\piitt[2]) - \tauit[2] \Dfistar[2](\pii[2];\piitt[2]) \leq 0.\end{equation*}
 Multiplying its rows with $\pii[1]$ leads to, 
   \begin{equation}\label{pf:3pointpi2}
 - \pii[1]((\piit[2] -\pii[2])   \xtil[t] + \fistar[2](\pii[2]) - \fistar[2](\piit[2])) + (\tauit[2] + 1) \pii[1]\Dfistar[2](\pii[2];\piit[2]) + \tauit[2] \pii[1]\Dfistar[2](\piit[2];\piitt[2]) - \tauit[2] \pii[1]\Dfistar[2](\pii[2];\piitt[2]) \leq 0.\end{equation}
 Since Lemma \ref{lm:proximal_gradient_lemma} implies that $\piit[2] = \grad \myfi[2] (\yit[2])$ for some $\yit[2]$ for all $t$, it follows from Lemma \ref{lm:smooth_breg} in the appendix that 
 $$\norm{\pii[1]}\pii[1]\Dfistar[2](\piit[2];\piitt[2]) \geq \normsq{\pii[1] (\piit[2] - \piitt[2])}/(2\Lfi[2]). $$
 Then we can obtain a $Q_2$ bound using a similar argument as above  
 \begin{equation}\label{pfpr1:Q2_bound}
 \begin{split}
\E[\nsumt \wt Q_2(\zt; z)] \leq& \wt[1]\tauit[2][11] \pii[1] \Dfistar[2](\pii[2]; \piit[2][0])  \\
& +  \tsum_{t=2}^N \tfrac{\wtt \thetat L_2 \norm{\pii[1]}}{\tauit[1]}   \normsq{\xtt - \xttt} + \tfrac{\wT L_2 \norm{\pii[1]}}{\tauT + \alpha} \normsq{\xT - x^{N-1}}.
\end{split}
\end{equation}

Last, we can derive from prox-mapping for $\xt$ in \eqref{ls:2-SSD-gen} in a similar fashion that 
\begin{equation}\label{pfpr1:Q0_bound}
 \begin{split}
 \E[\nsumt \wt Q_0(\zt; z) +  \tfrac{\wT(\etaT + \alpha)}{2} &\normsq{\xT - \xstar} ] + \E[\sumt \tfrac{\wt \etat}{2} \normsq{\xt - \xtt}]\leq 
  \sumt \wt \sigtl[x]^2 /(\etat + \alpha) \\
  &  \E[\tsum_{t=1}^{N} \{\wt \etat - \wtt (\etatt + \alpha)\}\normsq{x^{t-1} - \xstar}].
  \end{split}
 \end{equation} 
The desired inequality in \eqref{eq:2-sm-Q-conv} then follows from adding up \eqref{pfpr1:Q1_bound}, \eqref{pfpr1:Q2_bound} and \eqref{pfpr1:Q0_bound}.
 \end{proof}


We specialize Proposition \ref{pr:2-sm} to prove Theorem \ref{thm:2-sm}.\\
\textbf{Proof of Theorem \ref{thm:2-sm}}
Clearly, the proposed Algorithm \ref{alg:ssd_2} is a special case of \eqref{ls:2-SSD-gen} with $\Ui[1] = \Dfistar[1]$, $\piit[1] = \grad \myfi[1](\uyit[1])$ and $\piit[2] = \grad \myfi[2] (\uyit[2]).$ 

We begin by deriving a few basic properties satisfied by the stepsize choices in \eqref{stp:2-smo-non}. First, we provide some more useful characterizations of $\uyit[1]$ and $\uyit[2]$. Since $\thetat= (t-1)/t$ and $\tauit[2]=(t-1)/2$, we have $\uyit[1][0] = \xt[0] \in X$ and that
$$\uyit[2] = 2(\tsum_{l=1}^{t-1} l \xt[l] + t \xt[t-1])/[t(t+1)] \in X\ \forall t \geq 1.$$
Thus $\norm{\piit[2]} = \norm{\grad \myfi[2](\uyit[2])} \leq \Mi[2], \forall t.$ Similarly, setting $\hLa_2(x; \piit[2](\xi)):=\myfi[2](\uyit[2], \xi_2^1) - \grad \myfi[2](\uyit[2], \xi_2^1) \uyit[2] +\grad \myfi[2](\uyit[2], \hat{\xi_2})x + [\grad \myfi[2](\uyit[2], \xi_2^1) - \grad \myfi[2](\uyit[2], \hat{\xi_2})]\xtt$, we have
\begin{equation}\label{eq:uyit-decom}
\uyit[1] = 2[\tsum_{l=1}^{t-1} l \hLa_2(\xt[l]; \piit[2][l](\xi)) + t \hLa_2(\xtt; \piit[2](\xi)]/[t(t+1)], \forall t \geq 1.
\end{equation}
We now provide bounds for the variance constants in \eqref{cst:agg_var}. The following  inequality is valid for all $t \geq 1$:
\begin{align*}
\E&[\normsq{\{f_2(\uyit[2], \xi_2^1) - f_2(\uyit[2])\} + \{\grad \myfi[2](\uyit[2], \xi_2^1) - \grad \myfi[2](\uyit[2])\}(\xtt-\uyit[2]) + \{\grad \myfi[2](\uyit[2], \hat{\xi}_2)- \grad \myfi[2](\uyit[2])\} (\xtt - \xttt)}]\\
&\leq 2 \sig[\myfi[2]]^2 + 2 \sig[2]^2 \Dx^2 + \sig[2]^2 \Dx^2 \leq 2 \sig[\myfi[2]]^2 + 3 \sig[2]^2 \Dx^2.
\end{align*}
Thus $\sigtl[1]^2 \leq 2 \sig[\myfi[2]]^2 + 3 \sig[2]^2 \Dx^2$.

Next, we provide a bound on $\sigtl[x]^2$. Define  
$\uyitE[1]:= 2[\tsum_{l=1}^{t-1} l \La_2(\xt[l]; \piit[2][l] + \La_2(\xt; \piit[2])]/[t(t+1)] \in \Yi[1].$
By the Jensen's inequality, we have
$\E[\normsq{\uyit[1] - \uyitE[1]}] \leq \max_{l \leq t} \E[\normsq{\La_2(\xt[l], \piit[2][l](\xi))- \La_2(\xt; \piit[2])]}]\leq 5 \sig[2]^2 \Dx^2.$ Thus the $\Lfi[1]$-smoothness of $\myfi[1]$ implies that 
$$\E[\normsq{\grad \myfi[1](\uyit[1])}] \leq  2\E[\normsq{\grad \myfi[1](\uyitE[1])}] + 2\Lfi[1]^2\E[\normsq{\uyit[1] - \uyitE[1]}] \leq 2 \Mi[1]^2 +10 \Lfi[1]^2\Dx^2\sig[2]^2.  $$
Therefore, the conditional independence of $\grad \myfi[1](\uyit[1], \xi_1^0)$ and $\grad \myfi[2](\uyit[2], \xi_2^0)$ implies that 
\begin{align*}
\E[\normsq{\uyit[0](\xi) - \uyit[0]}] &= \E[\normsq{\grad \myfi[1](\uyit[1], \xi_1^0)\grad \myfi[2](\uyit[2], \xi_2^0) - \grad \myfi[1](\uyit[1])\grad \myfi[2](\uyit[2])}]\\
&\leq \sig[1]\Mi[2]^2 + \sig[1]^2\sig[2]^2 + \sig[2]^2\E[\norm{\grad\myfi[1](\uyit[1])}] \leq \sig[1]^2 \Mi[2]^2 + 2\Mi[1]^2\sig[2]^2 + \sig[2]^2(\sig[1]^2 + 10 \Lfi[1]^2 \sig[2]^2 \Dx^2).
\end{align*}
Therefore the desired upper bound on $\sigtl[x]^2$ is also valid.

Now we return to derive the function value convergence bound. Clearly, the requirements in Proposition \ref{pr:2-sm} is satisfied with $\beta = 1/\Lfi[1]$, $\mu =1$, and $\Mitl[2]= \Mi[2]$, so it follows from \eqref{eq:2-sm-Q-conv} that 
\begin{align}\label{pfthm1:Q}
  \begin{split}
  \E[\sumwt Q(\zt; z)]   &\leq \sumt \wt[1]\etat[1]  \normsq{\xt[0] - \xstar}/2
  + \sumt \wt \sigtl[x]^2 /\etat  + 2N \sigtl[1]^2  + \E\sumt \wt \inner{\pii[1]}{\yitilt[1]- \yitilt[1](\xi)}.
  \end{split}
  \end{align}
  In particular,  the Cauchy-Schwartz inequality implies that 
  \begin{equation}\label{eq:cauchy-sm}
  \E[\sumt \wt \inner{\pii[1]}{\yitilt[1]- \yitilt[1](\xi)}] = \E^{1/2}[\norm{\pii[1]}^2] \E^{1/2}[\normsq{\sumwt\{\yitilt[1]- \yitilt[1](\xi)\}}] \leq  N^{3/2}\sigtl[1]^2 \E^{1/2}[\norm{\pii[1]}^2]\ \forall N \geq 2.\end{equation}
  Moreover, if $\hpiit[1][N] = \grad \myfi[1](\myfi[2](\xbarT))$ and $\hpiit[2][N] = \grad \myfi[2](\xbarT)$, the strong duality and the weak duality relations in Lemma \ref{lm:duality}, and the convexity of $\La(x; \pii[1], \pii[2])$ with respect to $x$ imply
  \begin{align}\label{eq:2-sm-Q-convert}
  \begin{split}
  (\sumwt) (f(\xbarT) - f(\xstar)) &= \sumt \wt \La(\sumwt\xt/[\sumwt]; \hpiit[1][N], \hpiit[2][N]) - \sumt \wt f(\xstar) \\
  &\leq \sumt \wt[t] [\La(\xt; \hpiit[1][N], \hpiit[2][N]) - \La(\xstar; \piit[1], \piit[2])]\\
  &\leq \sumwt Q(\zt; (\xstar; \hpiit[1][N], \hpiit[2][N])).
  \end{split}
  \end{align}
  In view of the preceding two inequalities and that $\norm{\hpiit[1]} \leq \Mi[1]$, the desired function-value convergence rate can be derived from setting $z$ in \eqref{pfthm1:Q} to $(\xstar;\hpiit[2], \hpiit[2])$ and dividing the both sides by $(\sumwt)$.
\endproof

\vgap

\textbf{Proof of Theorem \ref{thm:2-sm-str}}
First, the bounds on $\sigtl[x]$ and $\sigtl[1]$ follows from the same argument as that of Theorem \ref{thm:2-sm}. 

Now we consider the convergence of $\normsq{x^N-\xstar}$. Let $\zstar:=(\xstar; \piistar[1], \piistar[2])$ denote the saddle point to \eqref{eq:2_sad}. Because our stepsizes satisfy the requirements in Proposition \ref{pr:2-sm} with $\mu=1$ and $\beta=1/\Lfi[1]$, the consequent bound in \eqref{eq:2-sm-Q-conv} holds for $z=\zstar$. Taking in account $Q(\zt; \zstar) \geq 0\ \forall t$ and that $\E[\nsumt \wt  \inner{\piistar[1]}{\yitilt[1]- \yitilt[1](\xi)}]=0$, we get the 
 \begin{equation*}
\E[\normsq{\xt[N] - \xstar}] \leq \tfrac{\tilde{L}}{\alpha N(N+1)} \normsq{\xt[0] -\xstar} + \tfrac{4}{\alpha N}( \Lfi[1] \sigtl[1]^2 + \tfrac{\sigtl[x]^2}{\alpha}), \forall N \geq 1.
\end{equation*}
In particular, the Jensen's inequality implies 
\begin{equation}\label{pfeq:smostr21}
\E[\normsq{\xbarT - \xstar}] \leq \tfrac{2\tilde{L}\log(N + 1)}{\alpha N^2} \normsq{\xt[0] -\xstar} + \tfrac{8}{\alpha N}( \Lfi[1] \sigtl[1]^2 + \tfrac{\sigtl[x]^2}{\alpha}).
\end{equation}
Next, let $\hpii[2]\in \partial \myfi[2](\xbarT)$ and $\hpii[1] = \grad \myfi[1](\myfi[2](\xbarT)).$ Applying Proposition \ref{pr:2-sm} again but with a reference point $(\xstar; \piistar[1], \hpii[2])$ for some $\hpii[1]\in \partial \myfi[2](\xbarT)$, we get
\begin{equation}\label{pfeq:smostr22}
\E [\nmyLa[\xbarT][\piistar[1], \hpii[2]]- f(\xstar)] \leq \tfrac{\tilde{L}}{\alpha N(N+1)} \normsq{\xt[0] -\xstar} + \tfrac{4}{\alpha N}( \Lfi[1] \sigtl[1]^2 + \tfrac{\sigtl[x]^2}{\alpha}).
\end{equation}
Moreover, since $\hpii[1] = \grad \myfi[1](\myfi[2](\xbarT)) \Leftrightarrow \myLa[\xbarT][\hpii[2]][2] = \myfi[2](\xbarT) \in  \partial \fistar[1] (\hpii[1]),$ the smoothness of $\myfi[1]$ implies
\begin{align*}
\E[f(\xbarT)& - \La(\xbarT; \piistar[1], \hpii[2])] = \E[\La(\xbarT;\piistar[1], \hpii[2])-\La(\xbarT;\piistar[1], \hpii[2])] = \E[\Dfistar[1](\piistar[1];\hpii[1])]\\
&= \E[D_{f_1}(\myfi[2](\xbarT);\myfi[2](\xstar))] \leq \tfrac{L_1}{2} \E[\normsq{\myfi[2](\xbarT)-\myfi[2](\xstar)}] \leq \tfrac{L_1 \Mi[2]^2}{2} \E[\normsq{\xbarT-\xstar}].
\end{align*}
Thus we can conclude the desired functional value convergence bound in \eqref{eq:2-str-sm-conv} by combining \eqref{pfeq:smostr21} and \eqref{pfeq:smostr22}.
\endproof
\vgap

\textbf{Proof of Theorem \ref{thm:2-strns}}
Clearly, the proposed Algorithm in \eqref{ls:2-ssd-strns} is a special case of \eqref{ls:2-SSD-gen} with $\Ui[1](\cdot) = \normsq{\cdot}$ and $\piit[2] = \grad \myfi[2] (\uyit[2]).$ 
The bound on $\sigtl[1]$ can be derived similarly to that of Theorem \eqref{thm:2-sm}. For $\sigtl[x]$, we have
\begin{align*}
\E[\normsq{\uyit[0](\xi) - \uyit[0]}] &= \E[\normsq{\piit[1]\{\grad \myfi[2](\uyit[2], \xi_2^0) - \grad \myfi[2](\uyit[2])\}}]\leq \Mi[1]^2 \sig[2]^2.
\end{align*}

It is easy to check the requirements in Proposition \ref{pr:2-sm} are satisfied with $\beta=1$ and $\mu=0$. Then function value convergence can also be deduced from an argument similar to that of Theorem \ref{thm:2-sm}.
\endproof

\section{General Nonsmooth Two-layer Problem}\label{sec:2-ns}

In this section, we study the general nonsmooth two-layer problem \eqref{eq:2_prob} where both $\myfi[1]$ and $\myfi[2]$ are stochastic and nonsmooth. As suggested in the introduction, the nonsmooth outer-layer function $\myfi[1]$ poses a critical challenge to the SSD method.  It is not possible to implement the prox-mapping for $\pii[1]$ (see Line 2 in \eqref{ls:2-SSD-gen}) with a strongly convex prox-function. On the one hand, since $\Dfistar[1]$ is not necessarily strongly convex, using it as the prox-function cannot provide the desired stabilization against its stochastic input, $\myfi[2](x)$. On the other hand, since $\myfi[1]$ is general non-smooth, the prox-mapping of $\fistar[1]$ with $\normsq{\cdot}$ as the prox-function may not be easily computable. 
We resolve the challenge by proposing a novel tri-conjugate  reformulation and developing a non-smooth stochastic sequential dual (nSSD) method. Specifically, Subsection \ref{sb:ssdp} introduces the nSSD method,
Subsection \ref{sb:2-ns-lower} presents a lower complexity result when the problem is strongly convex, and Subsection \ref{sb:2-ns-pf} concludes the section with the detailed convergence analysis.

 \subsection{The nSSD Method and Convergence Guarantee}\label{sb:ssdp}

 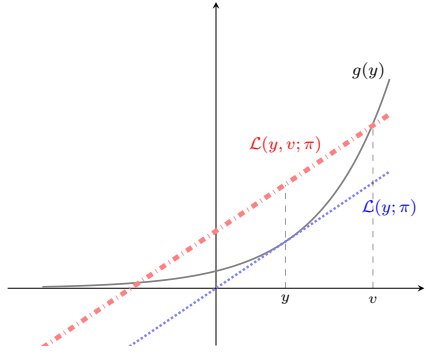
\begin{figure}[bt]
\centering
  \begin{tikzpicture}[scale=0.8]

   \begin{axis}[
        xmin=-3, xmax=3,
        ymin=-1, ymax=5,
        xtick distance=10, ytick distance=10, axis x line=middle, axis y line=middle ]

      \addplot [domain=-2.5:2.5, samples=100, thick, color=black!50]
        {0.3* exp(x) };

      \addplot [domain=-2.5:2.5, samples=100, densely dotted, line width=0.4mm, color=blue!50]
        {0.3 * exp(1) * x};

      \addplot [domain=-2.5:2.5, samples=100, dash dot, line width=0.8mm, color=red!50]
        {0.3 * exp(1) * x + 1};

      \draw [dashed, opacity=0.4] (axis cs:{2.26,0}) -- (axis cs:{2.26,2.8});
      \draw [dashed, opacity=0.4] (axis cs:{1,0}) -- (axis cs:{1,1.8});

      \node[color=blue, font=\footnotesize] at (axis cs: 2.5,1.4) {$\mathcal{L}(y; \pi)$};
      \node[color=red, font=\footnotesize] at (axis cs: 1, 2.5) {$\mathcal{L}(y, v; \pi)$};
      \node[color=black, font=\footnotesize] at (axis cs: 2.2,3.8) {$g(y)$};
      \node[color=black, font=\scriptsize] at (axis cs: 2.26,-0.2) {$v$};
      \node[color=black, font=\scriptsize] at (axis cs: 1,-0.2) {$y$};
    \end{axis}

  \end{tikzpicture}
  \caption{A comparison between the bi-conjugate reformulation, $\La_g(y; \pi)$, and the tri-conjugate reformulation, $\La_g(y,v;\pi)$ when $\pi \in \partial g(y)$.}\label{fig:lin-reform}
  \end{figure}

 First, we present a nested linearization reformulation for \eqref{eq:2_prob}. We begin by introducing a $\min-\max$ tri-conjugate reformulation to a convex non-smooth function $g$ defined on  $Y$:
\begin{equation}\label{def:lin-reform}
g(y) =  \max_{\pi \in \Pi} \min_{v \in Y} \{\La_g(y, v; \pi):= \inner{\pi}{y -v} + g(v)\},\end{equation}
where $\Pi:= \{s \in \partial g(y): y \in Y\}.$ There exist two interpretations to \eqref{def:lin-reform}. Treating $\pi$ as the Lagrange multiplier to the constraint $v=y$, \eqref{def:lin-reform} represents the Lagrangian dual  to the constrained optimization problem, 
$g(y) = \min_{v \in V} \{g(v) \text{ s.t. } v = y\}.$
Treating $v$ as the dual variable to $\pi$, \eqref{def:lin-reform} represents a certain tri-conjugate reformulation, i.e., an additional conjugate to the bi-conjugate: 
$$g(y) = \max_{\pi \in \Pi} \inner{\pi}{y}\quad  \underbrace{-\quad \quad  g^*(\pi).}_{:=\min_{v \in V} -\inner{\pi}{v} + g(v)}$$
A comparison between the tri-conjugate and the bi-conjugate reformulations is illustrated in Figure \ref{fig:lin-reform}.
Because of their close relationship, we use the common notation $\La_g$ for both, but emphasize the tri-conjugate reformulation by the auxiliary primal variable $v$. 
The key advantage of $\La_g(y, v; \pi)$ for us is the implementability of the prox-mappings for both $v$ and $\pi$ with $\normsq{\cdot}$ as the prox-function, which would be crucial for the nSSD method. 

Returning  to the two-layer problem in \eqref{eq:2_prob}, we consider a $\min-\min-\max-\max$ reformulation given by 
\begin{equation}\label{def:2-ns-reform}
\min_{(x,\vi[1]) \in X \times \Vi[1]}\ \max_{(\pii[1], \pii[2]) \in \Pii[1] \times \Pii[2]} \{\La(x, \vi[1]; \pii[1], \pii[2]):= \La_1(x, \vi[1]; \pii[1], \pii[2]) + u(x)\},\vspace{-1.5mm}
\end{equation}
 $$\text{where }\La_1(x, \vi[1]; \pii[1], \pii[2]):= \pii[1][\La_2(x; \pii[2]) -\vi[1]] + \myfi[1](\vi[1]), \text{ and } \La_2(x;\pii[2]) := \pii[2] x - \fistar[2](\pii[2]).$$
We take $\Vi[1]$ and $\Piitl[1]$ to be some compact and convex sets onto which (Euclidean) projection can be efficiently computed. We require them to be  large enough to contain some important points and to have finite radii  $D_{\Vi}$ and $\Mih[1]$:
\begin{equation}\label{cst:f1-lin}
\begin{split}
&\myfi[2](\xstar) \in \Vi[1],\ \text{and}\ \tmax_{\vi[1], \bar{v}_1} \norm{\vi[1]-\bar{v}_1} \leq D_{\Vi[1]}. \\
&\{s \in \partial \myfi[1](y_1), \yi[1] \in \R^{n_1}\} \subset \Piitl[1] \subset \R^{n_1}_+, \text{ and } \tmax_{\pii[1] \in \Piitl[1]} \norm{\pii[1]} \leq \Mih[1].
\end{split}
\end{equation}
For example, if $R_1\geq \max_{x \in X} \norm{\myfi[2](x)}$, $\Vi[1]$ can be chosen to be an Euclidean ball, $\{v_1 \in \R^{n_1}| \norm{v_1} \leq R_1\}$. If $\myfi[1]$ is $\Mi[1]$-Lipschitz continuous, $\Pii[1]$ can be chosen to be $\{\pii[1] \in \R^{n_1}| \norm{\pii[1]} \leq \Mi[1]\} \cap \R^{n_1}_+$. 
\vgap

We introduce some basic properties of the nested linearization reformulation \eqref{def:2-ns-reform} and a $Q$-gap function, which will lead to the nSSD method. The next lemma is a counterpart to Lemma \ref{lm:duality}. It relates \eqref{def:2-ns-reform} to the original problem in \eqref{eq:2_prob}. 
\begin{lemma}\label{lm:ns-duality}
The following relations between $\La$ in \eqref{def:2-ns-reform} and  $f$ in \eqref{eq:2_prob} are valid.
\begin{enumerate}
\item[a)] If $\vistar:=\myfi[2](\xstar)$, then $\La(\xstar, \vistar[1]; \pii[1], \pii[2]) \leq f(\xstar)\ \forall (\pii[1] \times \pii[2]) \in \Piitl[1] \times \Pii[2]$. 
\item[b)] Given a pair $(x, \vi) \in X \times \Vi[1]$, $\La(x, \vi; \hpii[1], \hpii[2]) \geq f(x)$ if $\hpii[1] \in \partial \myfi[1](\myfi[2](x))$ and $\hpii[2] \in \partial \myfi[2](x).$
\item[c)] Thus for any $z:=(x, \vi; \pii[1], \pii[2]) \in Z$, 
$f(x) - f(\xstar) \leq \La(x, \vi; \hpii[1], \hpii[2]) - \La(\xstar, \vistar[1]; \pii[1], \pii[2])$ if $\hpii[1] \in \partial \myfi[1](\myfi[2](x))$ and $\hpii[2] \in \partial \myfi[2](x).$
\end{enumerate}
\end{lemma}
\begin{proof}
Part a) follows from the non-negativity of $\pii[1] \in \Piitl[1]$ and  $\La_2(\xstar; \pii[2]) \leq \myfi[2](\xstar) = \vistar$, i.e.,
$$\La(\xstar, \vistar[1]; \pii[1], \pii[2]) = \pii[1][\La_2(\xstar; \pii[2]) - \vistar] + \myfi[1](\vistar) + u(\xstar) \leq \myfi[1](\myfi[2](\xstar)) + u(\xstar) = f(\xstar).$$
Part b) holds because the nested linearization reformulation is always large than the nested Lagrangian reformulation, i.e., \vspace{-3mm}
\begin{align*}
\La_1(x, \vi; \hpii[1], \hpii[2])&= \hpii[1][\La_2(x; \hpii[2]) -\vi[1]] + \myfi[1](\vi[1]) = \hpii[1]\La_2(x; \hpii[2]) -  [\hpii[1] \vi[1]- \myfi[1](\vi)]\\
 &\geq \hpii[1]\La_2(x; \hpii[2]) -  \max_{\vi \in \R^{n_1}}[\hpii[1] \vi[1]- \myfi[1](\vi)] = \hpii[1]\La_2(x; \hpii[2]) - \fistar[1](\hpii[1]) =  \La_1(x; \hpii[1], \hpii[1]) = \myfi[1](\myfi[2](x)),\vspace{-3mm}\end{align*}
where $\La_1(x; \hpii[1], \hpii[1])$ is defined in \eqref{eq:2_sad} and the last inequality follows from Lemma \ref{lm:duality}.b).
\end{proof}
\vgap
For simplicity, we use the notation $Z:=X \times \Vi \times \Piitl[1] \times \Pii[2]$ and $z:=(x, \vi; \pii[1], \pii[2])$ for the rest of the section. Lemma \ref{lm:ns-duality}.c) suggests a $Q$-gap function as an alternative optimality criterion for algorithm design. Specifically, given a point $\zt:= (\xt, \vit; \piit[1], \piit[2]) \in Z$, the $Q$-gap function with respect to a reference point $z\in Z$ is given by
\begin{equation}\label{eq:ns-Q-gap}
Q(\zt,z):= \La(\xt, \vit[1]; \pii[1], \pii[2]) - \La(x, \vi[1]; \piit[1], \piit[2]).
\end{equation}
It admits a decomposition to sub-optimality criteria given by 
\begin{equation} \label{decompose_gap-p}
Q(\zt, z) = Q_2(\zt, z) +  Q_1(\zt, z) + Q^v_1(\zt, z) +  Q_0(\zt, z),
\end{equation}
\vspace{-8mm}
\begin{align}
Q_2(\zt, z) &:= \mytLa[\xt, \vit][\pii[1], \pii[2]] - \mytLa[\xt, \vit][\pii[1], \piit[2]] 
= \pii[1][\pii[2]\xt - \fistar[2](\pii[2]) ] {\boxed{-\pii[1][\piit[2] \xt - \fistar[2](\piit[2])]}}, \label{def:Q2-ns} \\
Q_1(\zt, z) &:= \mytLa[\xt, \vit][\pii[1], \piit[2]] - \mytLa[\xt, \vit][\piit[1], \piit[2]] 
= \pii[1] [\myLa[\xt][\piit[2]][2] - \vit]   {\boxed{- \piit[1]{[\myLa[\xt][\piit[2]][2] -\vit]}}},\label{def:Q1-ns}\\
Q^v_1(\zt, z) &:= \mytLa[\xt, \vit][\piit[1], \piit[2]] - \mytLa[\xt, \vi][\piit[1], \piit[2]] 
= \boxed{\myfi[1](\vit) -\piit[1]\vit} - [\myfi[1](\vi) -\piit[1]\vi], \label{def:Qv1-ns}\\
Q_0(\zt, z) &:= \mytLa[\xt, \vi][\piit[1], \piit[2]] - \mytLa[x, \vi][\piit[1], \piit[2]] 
= {\boxed{ \piit[1]\piit[2] \xt + u(\xt)}} - [\piit[1]\piit[2]x + u(x)], \label{def:Q0-ns}
\end{align}
where $Q_2$, $Q_1$, $Q^v_1$, and $Q_0$ relate to the optimality of $\piit[2]$, $\piit[1]$, $\vit[1]$ and $\xt$, respectively. 

For the simple case where both $\myfi[1]$ and $\myfi[2]$ are deterministic, each iteration of the non-smooth sequential dual (nSD) method involves prox-mappings to reduce the $Q_2$, $Q_1$, $Q^v_1$ and $Q_0$ sequentially:
\begin{align}\label{ls:SDp_idea}
\begin{split}
&\piit[2] \leftarrow \targmax_{\pii[2] \in \Pii[2]}\ \pii[2]\xtt - \fistar[2](\pii[2]);\\
&\piit[1] \leftarrow \targmax_{\pii[1] \in \Piitl[1]}\ \pii[1][\yitilt[1] - \vit[t-1]] - \tauit[1] \normsq{\pii[1]- \piitt[1]}/2, \text{ where } \yitilt[1]:= \La_2(\xtt; \piit[2]);                                                                            \\
&\vit \leftarrow \targmin_{\vi \in \Vi} \ \  \inner{\myfi[1]'(\vit[t-1]) - \piit[1]}{\vi} + \gamit\normsq{\vi - \vitt}/2;\\
&\xt \leftarrow \targmin_{x \in X}\ \ \ \ \yitilt[0]x + \etat \normsq{x - \xtt}/2,\ \text{where } \yitilt[0]:= \piit[1]\piit[2].
\end{split}
\end{align}
 There are two simplifications compared to the SD method (c.f. \eqref{ls:SD_idea}).  First, rather than momentum-extrapolated prediction terms, values from the last iterate, i.e.,  $\xtt$ and $\La_2(\xtt; \piit[2])$ (or $\xtilt$ and $\yitilt[1]$ in \eqref{ls:SD_idea} with $\thetat=0$), are used as arguments for the $\piit[2]$ and $\piit[1]$-prox mappings. Second, the $\piit[2]$-prox mapping is implemented with $\tauit[2]=0$. These simplifications are justified because the $\bigO(1/\ep^2)$ oracle complexity can not be improved by the momentum-extrapolated acceleration even for the simple deterministic nonsmooth problem. Additionally, notice that a linear approximation of $\myfi[1]$, $\inner{ \myfi[1]'(\vitt)}{\cdot}$, is utilized when performing the prox-mapping for $\vit$.

We discuss the steps required to adapt \eqref{ls:SDp_idea} to the nested stochastic setting to arrive at the nSSD method shown in Algorithm \ref{alg:ssdp_2}. First, since $\tauit[2]=0$, Lemma \ref{lm:proximal_gradient_lemma} implies that the prox-mapping for $\piit[2]$ is equivalent to a gradient evaluation, $\piit[2] = \myfi[2]'(\xtt)$. Moreover, since the conjugate duality relationship in \eqref{rel:conj_duality} implies 
$$\La_2(\xtt; \piit[2]):= \fistar[2](\piit[2]) + \piit[2] \xtt =  [\myfi[2](\xtt) - \myfi[2]'(\xtt) \xtt ] +\myfi[2]'(\xtt) \xtt =  \myfi[2](\xtt),$$
$\myfi[2](\xtt, \xi_2^1)$ provides an unbiased estimator to $\La_2(\xtt; \piit[2])$  in Line 3 of Algorithm \ref{alg:ssdp_2}. Next,  the stochastic subgradient $\myfi[1]'(\vitt, \xi_1)$ is used in place of $\myfi[1]'(\vitt)$ in Line 4 of Algorithm \ref{alg:ssdp_2}, and the unbiased estimator of $\piit[1]\piit[2]$ is constructed from an independent estimator $\myfi[2]'(\xtt, \xi_2^0)$ in Line 5 of Algorithm \ref{alg:ssdp_2}.
\begin{algorithm}[htb]
\caption{Stochastic Sequential Dual-primal (nSSD) Method for Non-Smooth Two Layer Problem}
\label{alg:ssdp_2}
\begin{algorithmic}[1]
\Require $x^0\in X$, $\vit[0]\in\Vi[1]$, $\piit[1][0] \in\Piitl[1]$.
\For{$t = 1,2,3 ... N$}
\vspace{1mm}
\parState{Call $\SOi[2]$ to obtain independent estimates  $ \{f_2(\xtt, \xi_2^1),  \myfi[2]'(\xtt, \xi_2^0)\}$. \vspace{1mm}
}
\parState{Compute $\piit[1] \leftarrow \targmax_{\pii[1] \in \Piitl[1]}\ \pii[1][\yitilt[1](\xi) - \vit[t-1]] - \tauit[1] \normsq{\pii[1]- \piitt[1]}/2$ where $\yitilt[1](\xi) :=\myfi[2](\xtt, 
\xi_2^1)$.\vspace{1mm}
}
\parState{Compute $\vit \leftarrow \targmin_{\vi \in \Vi} \ \inner{\myfi[1]'(\vit[t-1], \xi_1) - \piit[1]}{\vi} + \gamit\normsq{\vi - \vitt}/2$.\vspace{1mm}
}
\parState{Set $\yitilt[0](\xi) \leftarrow  \piit[1] \myfi[2]'(\xtt, \xi_2^0)$. Compute $\xt \leftarrow \argmin_{x \in X} \yitilt[0](\xi) x + u(x) + {\etat} \norm{x -\xtt}^2/2$.
}
\EndFor
\State Return $\xbarT = \nsumt \wt \xt / (\sumwt)$.
\end{algorithmic}
\end{algorithm}

Now we present the convergence property of the proposed nSSD method. We need to specify the Lipschitz-continuity constants, $\Mi[1]$ and $\Mi[2]$, and the variance constants, $\sig[1]$, $\sig[\myfi[2]]$ and $\sig[2]$, associated with the layer functions:
\begin{align}\label{cst:ns-2}
\begin{split}
&\Mi[1] := \tmax_{y_1 \in \R^{n_1}} \norm{\myfi[1]'(y_1)},\ \Mi[2] := \tmax_{x \in X} \norm{\myfi[2]'(x)},\\
&\E [\normsq{\myfi[1]'(y_1, \xi) -\myfi[1]'(y_1)}] \leq \sig[1]^2\ \forall y_1 \in \R^{n_1} , \\ 
&\E [\normsq{\myfi[2]'(x, \xi) -\myfi[2]'(x)}] \leq \sig[2]^2\ \text{ and } \E [\normsq{\myfi[2](x, \xi) -\myfi[2](x)}] \leq \sig[\myfi[2]]^2 \leq \Dx^2\sig[2]^2\ \forall x \in X.
\end{split}
\end{align}
We  define aggregate variance constants of the stochastic arguments in Algorithm \ref{alg:ssdp_2} as
\begin{equation}\label{cst:ns-agg_var}
\sigtl[1] := \{\max_{t \geq 0} \E[\normsq{\yitilt[1](\xi) - \E[\yitilt[1](\xi)] }]\}^{1/2},
\ \sigtl[x] := \{\max_{t \geq 0} \E[\normsq{\yitilt[0](\xi) - \E[\yitilt[0](\xi)] }] \}^{1/2}.\end{equation} 
We are now ready to state the convergence result. The proof is deferred to Subsection \ref{sb:2-ns-pf}.

\begin{theorem}\label{thm:2-ns}
Consider a two-layer problem (c.f. \eqref{eq:2_prob}) with stochastic non-smooth functions $\myfi[1]$ and $\myfi[2]$. Let their Lipschitz-continuity constants and variance constants be defined in \eqref{cst:ns-2}, and the radius constants of the tri-conjugate reformulation of $\myfi[1]$ be defined in \eqref{cst:f1-lin}.
If the solution iterates $\{\xt\}$ is generated by Algorithm \ref{alg:ssdp_2}, then the variance constants in \eqref{cst:ns-agg_var} for the aggregate stochastic estimators in the algorithm satisfy $\sigtl[1] \leq \sig[\myfi[2]]\leq \sig[2] \Dx$ and $\sigtl[x] \leq \Mih[1] \sig[2]$. Moreover, if the stepsizes are given by 
\begin{equation}\label{stp:2-ns}
\wt = 1,\ \etat= 2 \max\{4 \Mih[1] \Mi[2], \sigtl[x] \} \sqrt{t}/ \Dx,\ \tauit[1]= \sqrt{t} \sigtl[1]/\Mih[1],\ \gamit= 2 \max\{4\Mih[1], \sig[1]\}\sqrt{t}/\DVi[1],
\end{equation}
the ergodic average solution $\xbarT$ satisfies 
\begin{equation}\label{eq:2-ns-f-conv}
\E[f(\xbarT) - f(\xstar)] \leq \tfrac{1}{\sqrt{N}}\{8\Mih[1] \Mi[2]\Dx + \sigtl[x] \Dx + 5 \Mih[1]\sigtl[1] + 8 \Mih[1]\DVi[1] + \sig[1]\DVi\}.
\end{equation}
\end{theorem}

Three remarks are in order regarding the result. First, if the total iteration number $N$ is known beforehand, we can set $\Vi[1]=\R^{n_1}$ and $\gamit=2 \max\{4\Mih[1], \sig[1]\}\sqrt{N}/\DVi[1]$ with $\DVi[1]=\Mi[2]\Dx$, and keep the same $\etat$ and $\tauit[1]$ in \eqref{stp:2-ns}, then the convergence bound in \eqref{eq:2-ns-f-conv} is still valid. In fact, the bound can be simplified further to 
\begin{equation*}
\E[f(\xbarT) - f(\xstar)] \leq \bigO\{(\Mih[1]+\sig[1])(\Mi[2] + \sig[2])\Dx/\sqrt{N}\}.
\end{equation*}
Second, any stepsize choice satisfying $\etat=\Theta(\sqrt{t})$, $\gamit=\Theta(\sqrt{t})$ and $\tauit[1] = \Theta(\sqrt{t})$ will lead to a stochastic oracle complexity of $\bigO(1/\ep^2)$ .
Third, the bound in \eqref{eq:2-ns-f-conv} is order-optimal because the term $\bigO(\Mih[1]\sigfi[2]/\sqrt{N})$  is not improvable even for strongly convex problems, which we will show in the next subsection.

\subsection{Lower Complexity Bound}\label{sb:2-ns-lower}
We develop a lower stochastic oracle complexity bound for the strongly convex NSCO problem with a general nonsmooth outer-layer function. Intuitively, a structured nonsmooth $\myfi[1]$ is easier than a general nonsmooth $\myfi[1]$, so a new lower complexity bound may appear unnecessary. 
However, since the prox-mappings in the nSSD method utilize more  information than the SSD method,
the nSSD method is not covered by the abstract scheme in \eqref{alg:abs-strns}, so the result in Theorem \ref{thm:l-strns} is no longer applicable. 
Specifically, recall that for a structured nonsmooth function $\myfi[1]$,  $\piit[1]$ generated by the prox-mapping in \eqref{alg:abs-strns}  can always be written as a sub-gradient at some $\yi[1]$. In contrast,  $\pii[1]$ in the tri-conjugate reformulation plays an additional role of a Lagrange multiplier, so $\piit[1]$ contains the residue $\vi - \myfi[2](x)$  in addition to the $\myfi[1]$ sub-gradients (see Line 3 and 4 in Algorithm \ref{alg:ssdp_2}). 
This motivates us to propose a more general abstract scheme to show that the $\bigO(1/\ep^2)$ oracle-complexity for nSSD to is also unimprovable. For simplicity, we assume $\Vi[1]:=\R^{n_1}$, $\Piitl[1]=\R^{n_1}_+ \cap \{\pii[1]\in \R^{n_1}| \norm{\pii[1]} \leq \Mi[1]\}$ and take $X$ to be a ball centered at $0$.

The abstract scheme consists of the following steps. In the beginning, supplied with input $\xt[0]\in X $, $\vit[1][0] \in \R^{n_1}$ and $\piit[1][0] \in \Piitl[1]$, the affine sub-spaces are initialized to $\Xt[0]:= \sp(\xt[0])$, $V_1^0:=\{v_1^0\}$, and $\Mit = \sp(\piit[1][0])$. In the $t$\ts{th} iteration, the updates are given by (the summations below represent the Minkowski sum and $[\cdot]^+_{\Piitl[1]}$ denotes the projection onto $\Piitl[1]$)
\emph{
\begin{align}\label{alg:abs-ns}
\begin{split}
&\text{Query}\ \SOi[2] \text{ to obtain } (\myfi[2](y_2^t, \xi_2^t), \myfi[2]'(y_2^t, \xi_2^t)) \text{ for some } y_2^t \in \Xt[t-1];\\
&\Mit := \sp([\tilde{\Mit[1][]}]_{\Piitl[1]}^+) + \tilde{\Mit[1][]} \text{ where }\\
&\quad \quad    \tilde{\Mit[1][]}:= \sp\{\myfi[2](y_2^t, \xi_2^t)+ \sum_{j=1}^{t} \myfi[2]'(y_2^j, \xi_2^j) x^j - v_1^0: x^j \in \Xt[t-1]\} + \sp\{\myfi[1]'(\pii[1] + v_1^0): \pii[1] \in \Mit[1][t-1]\} +\Mit[1][t-1]; \\
& \Xt[t] := \{ \myfi[2]'(y_2^j, \xi_2^j)^\top \pii[1]^\top: j \leq t, \pii[1] \in \Pii[1]^t\}
 + \Xt[t-1]. 
\end{split}
\end{align}
}
After $N$ iterations, the scheme  outputs some $x^N \in \Xt[N]$. 

Three remarks are in order for the abstract scheme in \eqref{alg:abs-ns}. First, the inclusion of $\sp([\tilde{\Mit[1][]}]_{\Piitl[1]}^+)$ in the construction of $\Mit$ is unusual for lower complexity models (see \cite{nesterov2003introductory}), but is necessary in our case because $\Piitl[1] \subset \R^{n_1}_+$ is not rotationally invariant. Second, it can be shown recursively that $\piit[1]$ and $\vit$ in Algorithm \ref{alg:ssdp_2} satisfies $\piit[1] \in \Mit$  and $\vit \in \vit[0] + \Mit$, thus the nSSD method is a special case of the abstract scheme. Third, since  $\Mit + v^1_0$ contains all the convex combinations of $\{\myfi[2](\xt[j], \xi_2)\}_{j=1}^t$ and $\Mit$ contains $\{\myfi[1]'(\yi[1]): \yi[1] \in \Mit[1][t-1] + v^1_0\}$, the abstract scheme still covers the SCGD-type algorithms \cite{mengdi2017stochastic} which updates $x$ with some pseudo-gradient, $\myfi[1]'(y_1)\myfi[2]'(x^j, \xi^j).$ Next we state the lower complexity result for the abstract scheme and its proof is deferred to the Appendix. 

\begin{theorem}\label{thm:l-ns}
Given problem parameters $\sig[\myfi[2]]\geq 0$, $\Mi[1] \geq 0$, $\ep > 0$ and $\alpha\leq\Mi[1]^2/(4\ep)$, there exists a nested two-layer problem \eqref{eq:2_prob} consisting of a general nonsmooth $\myfi[1]$  and a stochastic linear $\myfi[2]$ such that  $\myfi[1]$ is $\Mi[1]$-Lipschitz continuous, the variance of $\myfi[2](x, \xi)$ is bounded by $\sig[\myfi[2]]$ (c.f. \eqref{cst:ns-2}), and $u(x) = \alpha \normsq{x}/2$. If the abstract scheme in \eqref{alg:abs-ns} is initialized with $\xt[0]=0$, $\piit[1][0]=0$ and some $\vit[0]$ and its output from $x^N$ satisfies $\E[f(x^N) - f(\xstar)]\leq \ep$, then $N\geq \Omega(\Mitl[1]^2\sig[\myfi[2]]^2/\ep^2).$ 
\end{theorem}
 

\subsection{Convergence Proofs}\label{sb:2-ns-pf}

We present the detailed convergence analysis for Theorem \ref{thm:2-ns} in this subsection. We begin by providing a general $Q$-gap convergence bound. 

\begin{proposition}\label{pr:ns-2}
Let $\{\zt:=(\xt, \vit; \piit[1], \piit[2])\}$ be generated by Algorithm \ref{alg:ssdp_2} when applied to a two-layer problem (c.f. \eqref{eq:2_prob}) with stochastic nonsmooth functions $\myfi[1]$ and $\myfi[2]$. Let problem parameters be defined in \eqref{cst:ns-2}, \eqref{cst:f1-lin}, \eqref{cst:ns-agg_var}. If the stepsizes satisfy 
\begin{align}\label{pr2:stp}\begin{split}
\wt\etat \geq \wtt \etatt,\ \wt \tauit[1] \geq \wtt \tauitt[1],\ \wt \gamit \geq \wtt \gamit[1][t-1], \forall t \geq 2,
\end{split}
\end{align}
the following $Q$-gap bound holds for any $z:=(\xstar, \vistar; \hpii[1], \hpii[2])$, where $\vistar:=\myfi[2](\xstar)$ and $(\hpii[1], \hpii[2])\in \Piitl[1]\times \Pii[2]$ could potentially depend on $\{\zt\}$,
\begin{equation}\label{eq:2-ns-Q-bd}
\begin{split}
 \E[\sumwt Q(\zt; z)] \leq&  (8 \Mih[1]^2 \Mi[2]^2 + {\sigtl[x]^2}) \sumt \tfrac{\wt }{\etat}  + \tfrac{\wt[N]\etat[N]\Dx^2}{2}\\
 +& {\sigtl[1]^2} \sumt \tfrac{\wt }{\tauit[1]} + 2 \wt[N]\tauit[1][N] \Mih[1]^2 + \Mih[1] \sigtl[1]\sqrt{\sumt \wt^2}  \\ 
 +& (8 \Mih[1]^2 + {\sig[1]^2})\sumt \tfrac{\wt}{\gamit} + \tfrac{\wt[N]\gamit[1][N]\DVi[1]^2}{2}.
 \end{split}
 \end{equation} 
\end{proposition}
\begin{proof}
We begin by a developing a bound for $Q_1$ and $Q^v_1$. Line 3 of Algorithm \ref{alg:ssdp_2} implies \vspace{-2mm}
$$(\piit[1] - \hpii[1])[\vit[t-1] - \myfi[2](\xtt, \xi_2^1)] +  \tfrac{\tauit[1]}{2} (\normsq{\piit[1] - \hpii[1]} + \normsq{\piit[1] - \piitt[1]} - \normsq{\piitt[1] - \pii[1]}) \leq 0.$$
Denoting $\delta_1^t:=\myfi[2](\xtt, \xi_2^1) - \myfi[2](\xtt)$, we get \vspace{-2mm}
\begin{align*}
(\piit[1] - \hpii[1])\vit[t-1] = (\piit[1] - \hpii[1])\vit[t] - (\piit[1] - \hpii[1])[\vit[t] - \vit[t-1]] \geq (\piit[1] - \hpii[1])\vit[t] - 2 \Mih[1]\norm{\vit[t] - \vit[t-1]},
\end{align*}
and \vspace{-2mm}
\begin{align*}
\E[(\piit[1] - \pii[1])\myfi[2](\xtt, \xi_2^1)] =& \E[\piit[1] \delit[1]] + \E[\pii[1]\delit[1]]   
+\E[(\piit[1] - \pii[1])\{\La_2(\xtt; \piit[2]) - \La_2(\xt; \piit[2])\} + (\piit[1] - \hpii[1])\La_2(\xt; \piit[2])] \\
\leq & \E[\pii[1]\delit[1]] + \sigtl[1]^2/\tauit[1] + 2 \Mih[1] \Mi[2]\E[\norm{\xtt - \xt}] + \E[(\piit[1] - \hpii[1])\La_2(\xt; \piit[2])],
\end{align*}
where the last inequality follows from Lemma \ref{lm:prox-noise-bd} and that $\max\{\norm{\piit[1]}, \norm{\hpii[1]}\} \leq \Mih[1].$
Moreover, Line 4 of ALgorithm \ref{alg:ssdp_2} implies that \vspace{-3mm}
$$\inner{\vit - \vistar}{\myfi[1]'(\vit[t-1], \xi_1) - \piit[1]} +  \tfrac{\gamit[1]}{2} (\normsq{\vit - \vistar} + \normsq{\vit - \vitt} - \normsq{\vitt - \vistar}) \leq 0,$$
whereby \vspace{-2mm}
\begin{align*}
\E[\inner{\vit - \vistar}{\myfi[1]'(\vit[t-1], \xi_1)}] =&\ \E[\inner{\vit - \vistar}{\myfi[1]'(\vit[t-1], \xi_1) - \myfi[1]'(\vit[t-1])}] + \E[\inner{\vitt - \vistar}{\myfi[1]'(\vit[t-1])} + \inner{\vit - \vitt}{\myfi[1]'(\vit[t])}]\\
&+\E[\inner{\vit - \vitt}{\myfi[1]'(\vitt) - \myfi[1]'(\vit)}] \\
\geq &\ \E[\myfi[1](\vit) - \myfi[1](\vistar)] - \sig[1]^2/\gamit[1] - \E[2\norm{\vit - \vitt}\Mih[1]].
\end{align*}\vspace{-2mm}
Taken together, we get 
\begin{align*}
\E[Q^v_1(\zt; z)& + Q_1(\zt; z) + \tfrac{\gamit[1]}{2} (\normsq{\vit - \vistar}  - \normsq{\vitt - \vistar}) + \tfrac{\tauit[1]}{2} (\normsq{\piit[1] - \hpii[1]}  - \normsq{\piitt[1] - \pii[1]})]\\
\leq& 8 \Mih[1]^2 /\gamit[1] + \E[\pii[1]\delit[1]] + \sigtl[1]^2/\tauit[1] + \sig[1]^2/\gamit[1].
\end{align*}
Thus, the stepsize requirements in \eqref{pr2:stp} and the use of Cauchy-Schwartz inequality (see \eqref{eq:cauchy-sm})imply the next bound for the $\wt$-weighted sum 
\begin{equation}\label{pr2:Q1}
\begin{split}
\E[Q^v_1(\zt; z) + Q_1(\zt; z)] \leq& \wt[N](\gamit[N]\DVi^2/2  + 2\tauit[N] \Mih[1]^2 )  + {\sumwt} [8 \Mih[1]^2 /\gamit[1]  + \sigtl[1]^2/\tauit[1] + \sig[1]^2/\gamit[1]] \\
& + 2 \Mih[1] \Mi[2]\E[\sumwt \norm{\xtt - \xt}] + \Mih[1] \sig[1] \sqrt{\sumt\wt^2}.
\end{split}
\end{equation}
Similarly, we can get from Line 2 of Algorithm \ref{alg:ssdp_2} that 
\begin{equation}\label{pr2:Q2}
\E[\sumwt Q_2(\zt; z)] \leq 2 \Mih[1] \Mi[2]\E[\sumwt \norm{\xtt - \xt}].
\end{equation}
It follows from Line 5 of Algorithm \ref{alg:ssdp_2} that 
\begin{equation}\label{pr2:Q0}
\E[\sumwt Q_0(\zt; z)] \leq \sumwt \sigtl[x]^2 / \etat +  \wt[N]\etat[N] \Dx^2/2 - \E[\sumwt \etat\normsq{\xtt - \xt}].
\end{equation}
Thus, the desired $Q$-gap convergence bound in \eqref{eq:2-ns-Q-bd}  can be deduced from adding up \eqref{pr2:Q2}, \eqref{pr2:Q1} and \eqref{pr2:Q0}, and applying the Young's inequality.

\end{proof}

\textbf{Proof of Theorem \ref{thm:2-ns}}
The bounds on $\sigtl[1]$ and $\sigtl[2]$ can be derived directly from the definition of $\yitilt[1](\xi)$ and $\yitilt[0](\xi)$ in Algorithm \ref{alg:ssdp_2}. Applying Lemma \ref{lm:ns-duality} with $\hpiit[1][N] \in \partial \myfi[1](\myfi[2](\xbarT))$ and $\hpiit[2][N] \in \partial \myfi[2](\xbarT)$ (see \eqref{eq:2-sm-Q-convert}), we get 
\begin{align*}
\sumwt (f(\xbarT) - f(\xstar)) &\leq \sumt \wt \La(\tfrac{\sumwt\xt}{\sumwt}, \tfrac{\sumwt\vit}{\sumwt} ; \hpiit[1][N], \hpiit[2][N]) - \sumt \wt \La(\xstar, \vistar; \piit[1], \piit[2]) \\
&\leq \sumwt Q(\zt; (\xstar, \vistar; \hpiit[1][N], \hpiit[2][N])).
\end{align*}
Thus, the desired function gap convergence bound in \eqref{eq:2-ns-f-conv} can be derived from applying Proposition \ref{pr:ns-2} and dividing both sides of the resulting inequality by $(\sumwt)$. 
\endproof

\section{Multi-layer Problem}

In this section, we extend the SSD and the nSSD methods proposed in the last two sections to the multi-layer NSCO problem:
\begin{equation}\label{eq:multi-prob}
\min_{x \in X} \{f(x) := f_1 \circ f_2\circ \ldots \circ f_k (x) + u(x)\}.
\end{equation}
We  impose a multi-layer compositional convexity assumption similar to Assumption \ref{ass:cp_cv} throughout this section.
\begin{assumption}\label{mass:cp_cv}
A nested function $f_1 \circ f_2\circ \ldots \circ f_k (x)$ in \eqref{eq:multi-prob} is said to satisfy the \textit{compositional convexity} assumption if 
\begin{itemize}
\item Every layer function $f_i: \R^{n_i}\rightarrow \R^{n_{i-1}} $ is proper closed and convex. 
\item If $\myfi$ is not affine,  $\{\pii[1] \pii[2] \ldots \pii[i-1]: \pii[j] \in \partial \myfi[j](\yi[j]), \yi[j] \in \R^{n_j}, j \leq i-1\} \subset \R^{n_{i-1}}_+$.
\end{itemize}
\end{assumption}
Notice the second assumption is an extension of the two-layer monotonicity assumption in Assumption \ref{ass:cp_cv}. For a non-affine $\myfi$, our requirement of the product of all possible subgradients of outer layer functions being component-wise non-negative is slightly weaker than the usual assumption of every outer-layer function being component-wise non-negative.   

Our development consists of three parts. In Subsection \ref{sb:mul-sm}, we propose the SSD method for the smooth multi-layer problem. In Subsection \ref{sb:mul-ns},  we propose the nSSD method  to handle the general multi-layer problem constructed from an arbitrary composition of different types of layer functions. Finally, we present the detailed convergence analysis for the smooth problem in Subsection \ref{sb:mul-pf}.

\subsection{Smooth Multi-layer Problem}\label{sb:mul-sm}
In this subsection, we study  the multi-layer problem composed of only smooth layer functions. Our development follows the same structure as that of Subsection \ref{sb:2-smo-alg}. We start by introducing a nested Lagrangian reformulation which gives rise to a $Q$-gap function. The $Q$-gap function in turn motivates a conceptual SD method, and the SSD method fleshes out the conceptual algorithm by providing concrete implementations based on $\SOi$'s. Finally, convergence guarantees to the SSD method are provided. For notation simplicity, we will use $i:j$ to represent the collection from $i$ to $j$ (inclusive), and $i:$ as a shorthand for $i:k$. For example, $\pii[i:j]:=(\pii[i], \pii[i+1], \ldots, \pii[j])$ and  $\Pii[i:j] := \Pii[i] \times \Pii[i+1] \times \ldots \times \Pii[j]$.
We provide a nested Lagrangian reformulation to \eqref{eq:multi-prob},
 \begin{equation}\label{meq:sad}
 \min_{x \in X} \max_{\pii[\ito[1]] \in\Pii[\ito[1]] } \{\nmyLa[x][\pii[\ito[1]]]:= \nmyLa[x][\pii[\ito[1]]][1] + u(x)\},
 \end{equation}
 where $\Pii:=\dom(\myfi^*)$ denotes the domain of $\pii$ and the 
 nested Lagrangian function is defined recursively as
 \begin{equation}\label{meq:con_La}
 \nmyLa[x][\pii[\ito[i]]][i] := \begin{cases}
 x &\text{if } i=k+1,\\
 \pii \nmyLa[x][\pii[\ito[\ip]]][\ip] - \fistar(\pii) &\text{if } 1 \leq i \leq k.
 \end{cases}
 \end{equation}
 Let $z:=(x; \pii[\ito[1]])$ and $Z:= X \times \Pii[\ito[1]]$ denote the collections of decision variables and of their domains, then the reformulated problem can be related to \eqref{eq:multi-prob} through a duality result. 
 \begin{lemma}\label{mlm:duality}
Let $f$ and $\La$ be defined in \eqref{eq:multi-prob} and \eqref{meq:sad}, respectively. The following relations hold for all $x \in X$. 
\begin{enumerate}
 \item[a)] Weak Duality: $f(x) \geq \nmyLa[x][\pii[\ito[1]]] \ \forall z \in Z$. 
 \item[b)] Strong Duality:  for a given $x\in X$, $f(x) = \nmyLa[x][\hpii[\ito[1]]]$ if  $\hpii[i] \in \partial \myfi[i](\myfi[\ito[i+1]](x)) \ \forall i $.
 \item[c)] There exists some $(\piistar[1], \piistar[2], \ldots, \piistar[k])$ such that $\zstar:=(\xstar; \piistar[\ito[1]])\in Z$ is a saddle point, i.e., 
 \vspace{-2mm}
    $$\La(\xstar; \pii[\ito[1]])\leq \La(\xstar;  \piistar[\ito[1]]) \leq \La(x;  \piistar[\ito[1]])\ \forall (x,  \pii[\ito[1]])\in Z.$$
    \item[d)] For any $z \in Z$, an upper bound for the optimality gap of $x$ is given by:
    \vspace{-2mm}
    \begin{equation*}
    f(x) - f(\xstar) \leq \tmax_{\piibar[\ito[1]] \in \Pii[\ito[1]]} \La(x; \piibar[\ito[1]]) - \La(\xstar; \pii[\ito[1]]).
    \end{equation*}
 \end{enumerate} 
\end{lemma}

\begin{proof}
The proof is similar to that of Lemma \ref{lm:duality}. Part b) follows from the same argument.
As for Part a), let  $z:=(x; \pii[\ito[1]])\in Z$ be given and let $\hpii[i] \in \partial \myfi[i](\myfi[\ito[i+1]](x)) \ \forall i $ such that Part b) is valid, then we have
\begin{align*}
f(x) - \myLa[x][\pik] &= \tsum_{j=1}^{k} [\myLa[x][\pii[\ifrom[j-1]],\hpii[j],\hpii[\ito[j+1]]] - \myLa[x][\pii[\ifrom[j-1]],\pii[j],\hpii[\ito[j+1]]] ] \\
&= \tsum_{j=1}^k \underbrace{\pii[\ifrom[j-1]][(\hpii[j]\myfi[ \ito[j+1]](x) - \fistar[j](\hpii[j]))-(\pii[j]\myfi[ \ito[j+1]](x) - \fistar[j](\pii[j]))]}_{A_j}. 
\end{align*}
In particular,  $A_j =0$ if  $\myfi[j]$ is affine since the singleton $\Pii[j]$ implies $\hpii[j] = \pii[j]$. If $\myfi[j]$ is not affine,  we have $A_j\geq 0$ since $\hpii[j] \in \argmax_{\pii[j] \in \Pii[j]} \pii[j]\myfi[ \ito[j+1]](x) - \fistar[j](\pii[j])$ and Assumption \ref{mass:cp_cv} implies the non-negativity of $\pii[\ifrom[j-1]]$. Thus we get $f(x) \geq \myLa[x][\pii[\ito[1]]]$. Next, the derivations of Part c) and d) are also similar to that of Lemma \ref{lm:duality}.
\end{proof}
\vgap
Part d) of the preceding lemma motivates  the use of a multi-layer $Q$-gap function as an alternative optimality criterion. Specifically, given a point $\zt:=(\xt; \piit[1:])$,  the $Q$-gap function, defined with respect to a reference point $z \in Z$, is given by 
\begin{equation}\label{meq:gap_func}
Q(\zt, z):= \myLa[\xt][\pii[\ito[1]]] - \myLa[x][\pikt].
\end{equation}
One decomposition useful for algorithm  design is given by 
\begin{equation}
Q(\zt, z)= Q_0(\zt, z) + \sumk Q_i(\zt, z), 
\end{equation}
\text{ where}
\vspace{-2mm}
\begin{equation*}
\begin{split}
Q_0(\zt, z) &:= \myLa[\pikt][\xt] - \myLa[\pikt][x]\\
&= \boxed{\pikt \xt + u(\xt)} - \pikt x - u(x),
\end{split}
\end{equation*}
\begin{equation*}
\begin{split}
Q_i(\zt, z)&:= \myLa[\pii[\ifrom], \pii, \piit[\ito[i+1]]][\xt] - \myLa[\pii[\ifrom], \piit, \piit[\ito[i+1]]][\xt] \\
&=  \pii[\ifrom] \left(\pii\myLa[\piit[\ito[i+1]]][\xt][i+1] - \fistar(\pii)  \boxed{-[\piit\myLa[\piit[\ito[i+1]]][\xt][i+1] - \fistar(\piit)]} \ \right).
\end{split}
\end{equation*}
The decomposition motivates a conceptual SD method which sequentially updates $\pii$ to reduce $Q_i$ before updating $x$ in each iteration. Specifically, initialized to $\xt[0] \in X$ and $\piit[i][0] = \grad \myfi(\uyit[i][0])\ \forall i \in [k]$, the $t$\ts{th} iteration is given by 
\begin{align}\label{mls:SD_idea}
\begin{split}
&\text{for } i = k, k-1, \ldots, 1\\
&\quad \piit \leftarrow \argmax_{\pii \in \Pii} \pii \yitilt(\xi) - \fistar(\pii) - \tauit U_i(\pii;\piitt), \text{ where } \yitilt:= \La_{i+1}(\xtt; \piit[\ito[i+1]]) + \thetat \piit[\ito[i+1]](\xtt - \xt[t-2]);\\
&\xt \leftarrow \targmin_{x \in X}\  \yitilt[0](\xi)x + \tfrac{\etat}{2} \normsq{x - \xtt},\ \text{where } \yitilt[0]:= \piit[\ito[1]].
\end{split}
\end{align}
Here $U_i$ is some Bregman distance function, $\tauit$ and $\etat$ are stepsize parameters,  and $\yitilt(\xi)$'s and $\yitilt[0](\xi)$ are some unbiased estimators to $\yitilt$'s and $\yitilt[0]$, respectively.

\begin{figure}[htb]
\centering
\includegraphics[width=0.6\textwidth]{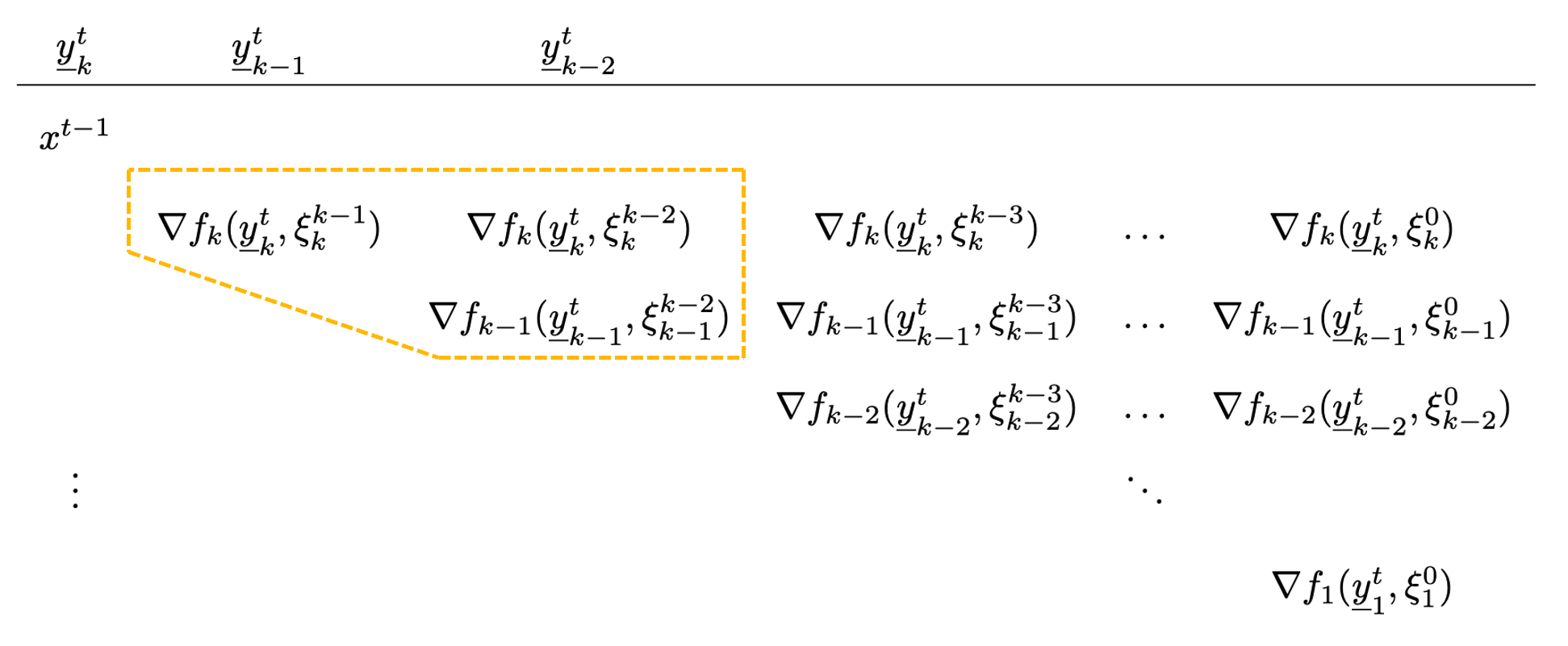}
\caption{\small Necessity of Resampling: for an unbiased estimator to   $\La_i(\xtt; \piit[i:]):= \La_{i}(\xtt; \gradfit,\ldots, \gradfit[k])$, we need conditionally independent estimators to $\gradfit, \gradfit[i+1],\ldots, \gradfit[k]$. However, for any $l > i$, $\{\gradfijt[j][l]\}_{j>l}$ are correlated with  $\uyit[l]$ such that the entire triangle of old samples $\{\gradfijt[l][j]\}_{l>i, j>i}$ (highlighted above) are correlated with $\uyit$ through $\{\uyit[j]\}_{j >i}$, i.e., $\E[\gradfijt[l][j]|\uyit, \uyit[l]]\neq \gradfit[l] \ \forall l,j>i$. Therefore,  $(k-i+1)$ new independent estimators need to be redrawn from $\SOi$'s, i.e., $\gradfijt[i][i-1], \gradfijt[i+1][i-1],\ldots, $\text{ and }$\gradfijt[k][i-2]$, for constructing an unbiased $\La_i(\xtt; \pixijt[i:][i-1])$. }
\label{mfig:resampling}
\end{figure}

Next we follow the same idea of choosing  $\Dfistar$, the dual Bregman distance function, to be $\Ui$  to obtain a concrete implementation of \eqref{mls:SD_idea} in Algorithm \ref{malg:ssd}. We use the conventions of $\La_{k+1}(x; \pii[\ito[k+1]]) = x$ and of $\grad \myfi[\ito[k+1]](\uyit[\ito[k+1]], \xi) = I$, the identity matrix, in Line 5 and 6. Moreover, as explained in Figure \ref{mfig:resampling}, the seemingly wasteful repeated independent calls to $\SOi$'s in Line 6 are essential for obtaining unbiased estimators of the nested Lagrangian function $\Laitt$.

\begin{algorithm}[htb]
\caption{ Stochastic Sequential Dual (SSD) Method for the Smooth Multi-layer Problem }
\label{malg:ssd}
\begin{algorithmic}[1]
\Require $x^0 \in X$.
\parState{Set $\uyit[k][0] := x^0$ and $x^{-1} := x^0$.}
\parState{Call $\SOi$ to obtain $\grad \myfi(\uyit[i][0], \hat{\xi}_i)$ and $\uyit[i-1][0] = \myfi(\uyit[i][0], \hat{\xi}_i)$ for $i= k, k-1, \ldots, 1$.}
\For{$t =  1,2,3 ... N$}
\For{$i=k, k-1, \ldots, 1$}
\parState{Set $\yitilt(\xi)\leftarrow \La_{i+1}[\xtt; \piit[\ito[i+1]](\xi_{\ito[i+1]}^i)] + \thetat \grad \myfi[\ito[i+1]](\uyitt[i+1:], \hxi[\ito[i+1]]) (\xtt - \xt[t-2]).$\\   Set $\uyit \leftarrow (\tauit \uyitt + \yitilt(\xi)) /(1+ \tauit).$}
\parState{Call $\SOi$ to obtain independent estimates  $ \{f_i(\uyit[i], \xi_i^j), \grad \myfi(\uyit, \xi_i^j)\}_{j=0}^{i-1} $ and $\grad \myfi(\uyit, \hxi)$. \\
Set $\La_{i}[\xtt; \piit[\ito[i]](\xi_{\ito[i]}^j)] \leftarrow \myfi(\uyit, \xi_i^j) + \grad \myfi[i](\uyit, \xi_i^j)\{\La_{i+1}[\xtt; \piit[\ito[i+1]](\xi_{\ito[i+1]}^j)] - \uyit[i]\}, \forall j = 1,\ldots, i-1$.\\
Set $\grad \myfi[\ito[i]](\uyit[\ito[i]], \hxi[\ito[i]]) \leftarrow \grad \myfi(\uyit; \hxi) \grad \myfi[\ito[i+1]](\uyit[\ito[i+1]], \hxi[\ito[i+1]]), \ \forall j = 1 \ldots i-1$. \\
 Set $\grad \myfi[\ito[i]](\uyit[\ito[i]], \xi_{\ito[i]}^0)\leftarrow \grad \myfi(\uyit; \xi_i^0)\myfi[\ito[i+1]](\uyit[\ito[i+1]], \xi_{\ito[i+1]}^0)$.}%
\EndFor
\parState{Set $\xt \leftarrow \argmin_{x \in X} \yitilt[0](\xi) x + u(x) + \tfrac{\etat}{2} \norm{x -\xtt}^2$ where $\yitilt[0](\xi):= \grad \myfi[\ito[1]](\uyit[\ito[1]], \xi_{\ito[1]}^0).$
}

\EndFor
\State Return $\xbarT := \sumwt \xp / \sumwt$.
\end{algorithmic}
\end{algorithm}
We now move on to provide the convergence result for Algorithm \ref{malg:ssd}. First, we need to define a few problem parameters. We require each layer function to be $\Li$-smooth and $\Mi$-continuous for some $\Mi \geq 1$\footnote{We assume $\Mi \geq 1$ to avoid trivialities.}, i.e., 
\begin{equation}\label{mcst:sm}
\norm{\grad \myfi(\yi)} \leq \Mi \text{ and } \norm{\grad \myfi(\yi) - \grad \myfi(\yibar)} \leq L_i \norm{\yi - \yibar} \ \forall \yi, \yibar \in \R^{n_i}. \end{equation}
For notation compactness, we will use the shorthand $\Mi[i:j]:= \prod_{l=i}^{j}\Mi[l]$ throughout this section.
We  assume the estimators from $\SOi$'s  to have bounded variances, $\sig[i] < \infty$, i.e.,
\begin{equation}\label{mcst:var}
\begin{split}
&\E[\normsq{\grad \myfi(\uyi, \xi_i) - \gradfit[i][]}] \leq \sig[i]^2,\ \E[\normsq{\myfi(\uyi, \xi_i)  - \myfi(\uyi) }] \leq \sig[\myfi]^2 \leq \sig[i]^2 \Mi[i+1:]^2 \DX^2 \ \forall \uyi \in \R^{n_i}.
\end{split}
\end{equation}
We also define a few aggregate variance constants for Algorithm \ref{malg:ssd} 
\begin{equation}\label{mcst:agg_var}
\sigtl[i] = \{\max_{t \geq 0} \E[\normsq{\yitilt(\xi) - \E[\yitilt(\xi)] }]\}^{1/2},
\ \sigtl[x] = \{\max_{t \geq 0} \E[\normsq{\yitilt[0](\xi) - \E[\yitilt[0](\xi)] }]\}^{1/2}.\end{equation}
The next theorem suggests some stepsize choices and provides the convergence rates of Algorithm \ref{malg:ssd} under both the non-strongly convex and the strongly convex settings. Its proof is deferred to Subsection \ref{sb:mul-pf}.

\begin{theorem}\label{mthm:ssd}
Let a smooth multi-layer function $f$ (c.f. \eqref{eq:multi-prob}) be given and let its problem parameters be defined in \eqref{mcst:sm} and \eqref{mcst:var}. If $\{\xt\}$ is generated by Algorithm \ref{malg:ssd} with
\begin{equation}\label{mstp:ssd-dual}
\wt = t,\ \thetat = (t-1)/t,\ \tauit=(t-1)/2, \forall i \in [k]\end{equation}
then the aggregate variance constants in \eqref{mcst:agg_var} satisfy  
\begin{equation}\label{meq:agg_var_bound}
 \sigtl[x]^2 \leq \tsum_{i=1}^{k} [\prod_{l\in[k]/\{i\}} \Mi[l]^2] \sig[i]^2 \text{ and } \sigtl[i-1]^2 \leq 2 [3 + (k-i)20^{k-i-1} ]  \tsum_{l=i}^{k} \sig[l]^2 \prod_{j\in i:k/\{l\}} (\Mi[j]^2 + \sig[j]^2) \DX^2\  \forall i \geq 2.\end{equation}
Let $\Ltl:= \sum_{i=1}^{k} \Mi[1:i-1] \Li \Mi[i+1:]^2$ denote the aggregate smoothness constant of $f$. Then the stepsize $\etat$ can be selected in accordance with the strong convexity modulus of $u(x)$ to obtain the following convergence guarantees. 
\begin{enumerate}
\item[a)] If $u(x)$ is non-strongly convex, selecting $\etat=\max\{2\Ltl/(t+1), \sigtl[x]\sqrt{t}/\Dx\}$ leads to 
\begin{equation}\label{meq:sm-ns-conv}
\E[f(\xbarT) - f(\xstar)] \leq \tfrac{\Ltl}{N(N+1)} \normsq{\xt[0] - \xstar} + \tfrac{4}{N} \tsum_{i=1}^{k-1} \Mi[1:i-1] \Li \sigtl[i]^2 + \tfrac{4}{\sqrt{N}} (\sigtl[x] \Dx + \tsum_{i=1}^{k-1} \Mi[1:i] \sigtl[i]).
\end{equation}
\item[b)] If $u(x)$ is strongly convex with a modulus $\alpha >0$, selecting $\etat=\max\{2\Ltl/(t+1), \alpha(t-1)/2\}$ leads to 
\begin{align}\label{meq:sm-str-conv}
\begin{split}
\E[f(\xbarT) - f(\xstar)] &\leq [\tfrac{\sum_{i=1}^{k-1}\Mi[1:i-1]\Li \Mi[i+1:]^2}{\alpha}  + 1][\tfrac{\log(N + 1)\tilde L}{N^2} \normsq{\xt[0] - \xstar} + \tfrac{4}{N}( \tsum_{i=1}^{k-1} \Mi[1:i-1] \Li \sigtl[i]^2 + \tfrac{\sigtl[x]^2}{\alpha})], \\
\E[\normsq{\xt[N] - \xstar}] &\leq \tfrac{\tilde{L}}{\alpha N(N+1)} \normsq{\xt[0] -\xstar} + \tfrac{4}{\alpha N}( \tsum_{i=1}^{k-1} \Mi[1:i-1] \Li \sigtl[i]^2 + \tfrac{\sigtl[x]^2}{\alpha}).
\end{split}
\end{align}
\end{enumerate}
\end{theorem}

A few remarks are in order regarding the result. In the deterministic case with $\sig[i]=0\ \forall i$, Algorithm \ref{malg:ssd} achieves the optimal oracle complexity under the non-strongly convex setting, and a restarted version of it (see Section 4.2.3 of \cite{LanBook}) achieves the optimal oracle complexity under the strongly convex setting. In the stochastic case, Algorithm \ref{malg:ssd} achieves the order-optimal stochastic oracle complexity of $\bigO(1/\ep^2)$ ($\bigO(1/\ep)$) with the parameter-free stepsizes in \eqref{mstp:ssd-dual} and any $\etat=\Theta(\sqrt{t})$ ($\etat=(t-1)\alpha/2$) under the non-strongly (strongly) convex setting. Moreover, the strong assumption of each layer function being $\Mi$-Lipschitz continuous (c.f. \eqref{mcst:sm}) is only required to attain the desired constant dependence in the deterministic case. Without it, the order-optimal oracle complexity of the SSD method is still valid for any $\etat=\Theta(\sqrt{t})$ ($\etat=(t-1)\alpha/2$). 

Additionally, we remark that the geometric dependence of $\sigtl[i]$ with respect to $k$ arises from the geometric growth of the range of the nested Lagrangian function. For example, consider the innermost function $\myfi[k]$, the diameter of $\{\La_k(x; \pii[k])\}$ is strictly larger than the range of $\myfi[k]$ in the worst case, i.e., 
$$D_k := \max_{x, y\in X }\max_{\pii[k], \pii[k]' \in \Pii[k]} \norm{\La_k(x; \pii[k]) - \La_k(y; \pii[k]')} \geq  c\Mi[k]\DX, \ c> 1.$$
Since $\La_{k-1}$ takes $\La_k$ as its argument (see \eqref{meq:con_La}), we could have $D_{k-1} \geq c^2 \Mi[k-1:]\Dx$, and $D_{i} \geq c^{k-i} \Mi[k-i:]\Dx$. This shows that  $\uyit$ and hence $\text{Var}[\gradfijt \uyit ]$ could become very large if the number of nested layer functions is large. Indeed, such a geometric dependence on $k$ appears to be an inevitable feature of the multi-layer nested function. The Lipschitz continuity constant of the nested function, $\Mi[1]\Mi[2]...\Mi[k]$, also scales geometrically with $k$. 
In practice, since $k$ is always a fixed problem parameter and it does not appear in the exponent of $\epsilon$, we expect the SSD method to converge relatively quickly.



\subsection{The nSSD Method for the General Nested Composition Problem}\label{sb:mul-ns}
In this subsection, we propose a multi-layer nSSD method to handle the general nested composition problem. Specifically, we consider the $k$-layer problem in \eqref{eq:multi-prob} constructed from an arbitrary composition of smooth, structured non-smooth and general non-smooth layer functions. We use the sets $\frS$, $\frP$ and $\frN$ to denote their respective layer indices. Our development follows the same pattern as that of the preceding subsection. We  propose a nested linearization reformulation, devise an $Q$-gap function as an alternative optimality criterion, and design the nSSD method to sequentially minimize each component of some decomposition of the $Q$-gap function.

First, we introduce the  multi-layer nested linearization reformulation. Because the general problem could involve non-smooth outer-layer functions, we need to generalize \eqref{meq:con_La} to include the tri-conjugate reformulation.  Let $\mathfrak{N}_o:=\frN \cap [k-1]$ denote the set of indices of non-smooth outer-layer functions. For simplicity, we use the notation $\vi[i:]:=\{\vi[l]\}_{ l \in \frNo \cap i:k}$ to denote the collections of auxiliary primal variables, and use $\Vi[i:]:= \prod_{l \in \frNo \cap i:}\Vii[l]$ to denote their corresponding domains (see \eqref{def:lin-reform}). 
If $i \notin \frNo$, we define $\Pii:= \dom(\fistar)$. 
If $i \in \frNo$, we choose $(\Vii, \Pii)$ to satisfy $\myfi[i:](\xstar) \in \Vii $, $\dom(\fistar) \subset \Pii$ and the monotonicity requirement in Assumption \ref{mass:cp_cv}, i.e., $\{\pii[1:l-1]: \pii[j] \in \Pii[j], j \leq l-1\} \subset \R^{n_{l-1}}_+$ if $\myfi[l]$ is not affine. Then the multi-layer nested linearization reformulation is given by 

 \begin{equation}\label{meq:sadp}
 \min_{x \in X} \min_{\vi[1:]\in \Vii[1:]} \max_{\pii[\ito[1]] \in\Pii[\ito[1]] } \{\nmyLa[x][\pii[\ito[1]]]:= \nmyLa[x][\pii[\ito[1]]][1] + u(x)\},
 \end{equation}
 where nested linearization function is defined recursively according to
 \begin{equation}\label{meq:con_Lap}
 \nmyLa[x][\pii[\ito[i]]][i] := \begin{cases}
 x &\text{if } i=k+1,\\
 \pii[i] [\nmyLa[x][\pii[\ito[\ip]]][\ip] - \vii] + \myfi(\vii) & \text{if }  i \in \frNo, \\
 \pii \nmyLa[x][\pii[\ito[\ip]]][\ip] - \fistar(\pii) &\text{otherwise.}
 \end{cases}
 \end{equation}
 Let $z:= (x, \vi[1:]; \pii[1:])$ and $Z:= X \times \Vi[1:] \times \Pii[1:]$. The next lemma relates the above reformulation to the original problem in \eqref{eq:multi-prob}. 
 \begin{lemma}\label{mlm:ns-duality}
The following relations between $\La$ in \eqref{meq:sadp} and  $f$ in \eqref{eq:multi-prob} are valid.
\begin{enumerate}
\item[a)] If $\vistar[i]=\myfi[i:](\xstar)\ \forall i \in \frNo$, then $\La(\xstar, \vistar[1:]; \pii[1:]) \leq f(\xstar)\ \forall \pii[1:] \in \Pii[1:]$. 
\item[b)] Given a pair $(x, \vi[1:]) \in X \times \Vi[1:]$, $\La(x, \vi; \hpii[1:]) \geq f(x)$ if $\hpii \in \partial \myfi[i](\myfi[i+1:](x))\ \forall i$. 
\item[c)] For any $z:=(x, \vi[1:]; \pii[1:]) \in Z$, 
$f(x) - f(\xstar) \leq \La(x, \vi[1:]; \hpii[1:]) - \La(\xstar, \vistar[1:]; \pii[1:])$ if $\hpii \in \partial \myfi[i](\myfi[i+1:](x))\ \forall i.$
\end{enumerate}
\end{lemma}
\begin{proof}
The proof is a straightforward extension to that of Lemma \ref{lm:ns-duality}. We need to show recursively that $\La_i(\xstar, \vistar[i:]; \pii[i:]) \leq \myfi(\myfi[i+1:](\xstar))\ \forall i $ for Part a), and that $\myfi[i:](x) \leq \La_i(x, \vi[i:]; \hpii[i:])\ \forall i$ for Part b).
\end{proof}
\vgap
Part c) of the preceding lemma motivates us to use the following $Q$-gap function, defined with respect to some reference point $\zbar \in Z$, as an alternative optimality criterion
$$Q(z; \zbar):= \La(x, \vi[1:]; \piibar[1:]) - \La(\xbar, \vibar[1:]; \pii[1:]).$$
A decomposition of the $Q$-gap function useful for algorithm design is given by 
\begin{equation}\label{meq:Q-ns-decom}
Q(z; \zbar) = Q_0(z; \zbar)+ \tsum_{i=1}^{k}Q_i(z;\zbar) + \sum_{i \in \frNo} Q^v_i(z; \zbar),\end{equation}
where 
\begin{align*}
Q_0(\zt, z) :=& \myLa[\pikt][\xt, \vi[1:]] - \myLa[\pikt][x, \vi[1:]]
= \boxed{\pikt \xt + u(\xt)} - \pikt x - u(x), \\
Q^v_i(\zt, z):=& \myLa[\pikt][\xt, \vvit[1:i-1], \viit, \vi[i+1:]] - \myLa[\pikt][x, \viit[1:i-1], \vii[i], \vii[i+1:]]:= \boxed{\piit[1:i-1][t] [\myfi(\viit[i]) - \piit \vvit[i] ]} - \piit[1:i-1][t] [\myfi(\vi[i]) - \piit \vi[i] ],\\
Q_i(\zt, z):=& \myLa[\pii[1:i-1], \pii[i], \piit[\ito[i+1]]][\xt,\vvit[1:]] - \myLa[\pii[\ifrom], \piit, \piit[\ito[i+1]]][\xt, \vvit[1:]] \\
=&\begin{cases}  \pii[\ifrom] \{\pii[i][\myLa[\piit[\ito[i+1]]][\xt, \vvit[i+1:]][i+1] - \vvit[i]]  \boxed{-\piit[i][t][\myLa[\piit[\ito[i+1]]][\xt, \vvit[i+1:]][i+1] - \vvit[i]]} \ \} &\text{if } i \in \frNo, \vspace{2mm} \\
\pii[\ifrom] \{\pii\myLa[\piit[\ito[i+1]]][\xt, \vvit[i+1:]][i+1] - \fistar(\pii)  \boxed{-[\piit\myLa[\piit[\ito[i+1]]][\xt, \vvit[i+1:]][i+1] - \fistar(\piit)]} \ \} &\text{otherwise}.
\end{cases}
\end{align*}
It motivates the multi-layer nSSD method shown in Algorithm \ref{malg:ssdp}, which attempts to sequentially reduce $Q_k,\  Q_{k-1},$ $\ldots,\  Q_1$, $\{Q^v_j\}_{j \in \frNo}$ and $Q_0$ in each iteration. Again we use the repeated sampling to construct unbiased estimates of the arguments to the prox-mappings according to 
\begin{equation}\label{meq:gen-repeated-sampling}
 \begin{split}
 \La_i(\xtt,& \vixijtt[i:][j]; \pixijt[i:][j])\\
 &:=\begin{cases}
    \xtt &\text{if } i=k+1\\
     \gradfijt\{\La_{i+1}[\xtt, \vixijtt[i+1:];\pixijt[i+1:]]- \uyit\} + \myfi(\uyit, \xi_i^j) &\text{else if } i \in \frS \\
     \piit {\La_{i+1}[\xtt, \vixijtt[i+1:];\pixijt[i+1:]}] - \fistar(\piit) &\text{else if } i \in \frP\\
     \piit \{\La_{i+1}[\xtt, \vixijtt[i+1:];\pixijt[i+1:]]- \viit[i][t-1]\} + \myfi(\viit[i][t-1], \xi_i^j) &\text{else if } i \in \frNo\\
     \myfi[k](\xtt, \xi_k^j) &\text{otherwise, i.e., $i=k$ and $k\in \frN$.}
 \end{cases}\\
 \pihxit[i:]& := \begin{cases}
    I &\text{if } i\geq k+1\\
    \grad \myfi(\uyit, \hxi) \pihxit[i+1:] &\text{else if } i \in \frS\\
    \piit \pihxit[i+1:] &\text{else if } i \in \frNo \cup \frP\\
    \myfi[k](\xtt, \hxi[k]) &\text{otherwise, i.e., $i=k$ and $k\in \frN$.}
 \end{cases}
 \end{split}
 \end{equation} 
Here $\{\hxi\}\cup\{\xi_i^j\}_{0 \leq j<i}$ represent independent samples drawn from $\SOi$ in each iteration. Moreover, the assignment of $\piit= \grad\myfi(\uyit)$ in Line 6 of Algorithm \ref{malg:ssdp} represents the association of $\piit$ with the primal point $\uyit$, rather than the computation of $\grad\myfi(\uyit);$ the usage of $\piit$  in the method are satisfied by calling $\SOi$ at $\uyit$ in accordance with \eqref{meq:gen-repeated-sampling}.
 It is also interesting to note that Algorithm \ref{malg:ssdp} simplifies to  Algorithm \ref{alg:ssdp_2}   if  $k=2$ and $\{1,2\}\subset \frN$, and to  Algorithm \ref{malg:ssd} if $[k] \subset \frS$.
\begin{algorithm}[htb]
\caption{ Stochastic Sequential Dual-primal (nSSD) Method for the General Multi-layer Problem }
\label{malg:ssdp}
\begin{algorithmic}[1]
\Require $x^0 \in X$, $\piit[i][0] \in \Pii\ \forall i \in \frP$, $(\viit[i][0],\piit[i][0]) \in \Vii \times \Pii\ \forall i \in \frNo$.
\parState{Set $\uyit[k][0] := x^0$ and $\uyit[i][0] = \myfi[i+1](\uyit[i+1][0], \hat{\xi}_{i+1}) \forall i \leq k-1 $. Set $\piit[i][0] = \grad \myfi(\uyit[i][0])\ \forall i \in \frS$ and $\piit[k][0] = \myfi[k]'(x^0)$ if $k\in \frN$.}
\For{$t =  1,2,3 ... N$} 
\For{$i=k, k-1, \ldots, 1$}
\parState{Set $\yit(\xi):=\La_{i+1}(\xtt, \vixijt[i+1:][i][t-1]; \pixijt[i+1:][i])$ and \\
$\yitilt(\xi):= \La_{i+1}(\xtt, \vixijt[i+1:][i][t-1]; \pixijt[i+1:][i]) + \pihxit[i+1:][t-1](\xtt - \xt[t-2])$ (c.f. \eqref{meq:gen-repeated-sampling}).}
\If{$\myfi$ is smooth}
\parState{set $\piit\leftarrow \grad \myfi(\uyit)$ where $\uyit:= (\tauit \uyitt + \yitilt(\xi))/(1 + \tauit).$}
\ElsIf{$\myfi$ is structured non-smooth}
\parState{set $\piit[i,j] \leftarrow \argmax_{\pii[i,j] \in \Pii[i,j]} \pii[i,j][\yitilt(\xi)] - \tauit \normsq{\pii[i,j]- \piitt[i,j]}/2\ \forall j \in [n_{i-1}]$.}
\ElsIf{$i=k$ and $\myfi$ is general non-smooth}
\parState{set $\piit[k] \leftarrow \myfipr[k](\xtt)$.}
\Else{\ $i \in \frNo$ ,i.e., $\myfi$ is general non-smooth}\parState{
Set $\piit[i,j] \leftarrow \argmax_{\pii[i,j] \in \Pii[i,j]} \pii[i,j] [\yit(\xi) - \viit[i][t-1]] - \tauit \normsq{\pii[i,j] - \piitt[i,j]}/2\ \forall j \in [n_{i-1}].$}
\EndIf
\EndFor
\For{$i \in \frNo$}
\parState{Set $\viit[i]\leftarrow \argmin_{\vii[i] \in \Vii[i]} \pixijt[1:i-1][0] [\myfipr(\viit[i][t-1], \hxi) - \piit]\vii + \gamiit \normsq{\vii - \viit[i][t-1]}/2$.} 
\EndFor
\parState{Set $\xt \leftarrow \argmin_{x \in X} \yitilt[0](\xi) x + u(x) + \tfrac{\etat}{2} \norm{x -\xt}^2$ where $\yitilt[0](\xi):= \pixijt[1:][0].$
}
\EndFor
\State Return $\xbarT := \sumwt \xp / \sumwt$.
\end{algorithmic}
\end{algorithm}

Now we present the convergence result for the nSSD method. First, we need to define a few problem parameters characterizing the Lipschitz-continuity and Lipschitz-smoothness of the layer functions:
\begin{align}\label{mcst:sddp}
\begin{split}
&\norm{\grad\myfi(\yi) - \grad \myfi(\yibar)} \leq \Lfi \norm{\yi - \yibar}\text{ and }\norm{\grad \myfi(\yi)} \leq \Mi\ \forall \yi , \yibar \in \R^{n_i}\ \forall i \in\frS.\\
&\norm{\pii} \leq \Mi\ \forall \pii \in \Pii \text{ and } \norm{\vii - \bar{v}_i} \leq \Mi[i+1:]\Dx\ \forall \vii, \bar{v}_i \in \Vii, \forall i \in \frNo.\\
&\norm{\pii} \leq \Mi\ \forall \pii \in \Pii\ \forall i \in \frP, \text{ and } \norm{\myfi[k]'(x)} \leq \Mi[k] \ \forall x \in X \text{ if } k \in \frN.
\end{split}
\end{align}
We also need the aggregate variance bounds for the unbiased  estimators  in Algorithm \ref{malg:ssdp}:
\begin{align}\label{mcst:gen-agg-var}
\begin{split}
&\sigtl[i] = \max\{\{\max_{t \geq 0} \E[\normsq{\yitilt(\xi) - \E[\yitilt(\xi)] }]\}^{1/2}, \{\max_{t \geq 0} \E[\normsq{\yit(\xi) - \E[\yit(\xi)] }]\}^{1/2}\},\\
& \sigtl[x] = \{\max_{t \geq 0} \E[\normsq{\yitilt[0](\xi) - \E[\yitilt[0](\xi)] }] ]\}^{1/2},\\
&\sigtl[v,i] = \{\max_{t \geq 0} \E[\normsq{\pixijt[1:i-1][0] [\myfipr(\viit[i][t-1], \hxi)] - \E[\pixijt[1:i-1][0] [\myfipr(\viit[i][t-1], \hxi) ] }]\}^{1/2}.
\end{split}
\end{align}
We are ready to state the convergence result. 

\begin{theorem}\label{mthm:ssdp}
Consider a general multi-layer function $f$ \eqref{eq:multi-prob} with
$\frS$, $\frP$ and $\frN$ denoting the indices of smooth, structured nonsmooth and general nonsmooth layer functions, respectively. Let $\frNi:=[i] \cap \frN$ denote the indices of general non-smooth functions outer to  $\myfi[i+1]$. Let the 
problem parameters be defined in \eqref{mcst:sddp} and \eqref{mcst:gen-agg-var}. Suppose $\{\xt\}$ are generated by Algorithm \ref{malg:ssdp} with the prox-update stepsizes chosen according to 
\begin{equation}\label{mstp:ssdp}
\begin{split}
&\wt := t,\ \thetat := (t-1)/t,\\\
&\tauit:=(t-1)/2,\ \etait :=2 \Mi[1:i-1] \Li \Mi[i+1:]^2/(t+1), \forall i \in \frS\\
&\tauit:= \max\{\sigtl[i]\sqrt{t},\ \Mi[i+1:]\Dx\}/\Mi,\ \etait = \Mi[1:] /\Dx,\ \forall i \in \frP\\
&\tauit:= \sigtl[i]\sqrt{t}/\Mi,\ \gamiit:= \max\{(|\frNi[1:i]| + 1) \Mi[1:i], \sigtl[i]\}\sqrt{t} / (\Mi[i+1:]\DX),\ \forall i \in \frNi[1:k-1]\\
&\etait:= \Mi[1:]\sqrt{t}/\Dx\ \forall i \in \frN.
\end{split}
\end{equation}
 Then the stepsize $\etat$ can be selected in accordance with the strong convexity modulus of $u(x)$ to obtain the following convergence guarantees. 
\begin{enumerate}
\item[a)] If $u(x)$ is non-strongly convex, selecting $\etat:=\max\{\sumk \etait, \sigtl[x]\sqrt{t}/\Dx\}$ leads to 
\begin{equation}\label{meq:mixed-conv}
\begin{split}
\E[f(\xbarT) - f(\xstar)] \leq \bigO\{&\tsum_{i \in \frS}[\tfrac{\Mi[1:i-1]\Li\Mi[i+1:]^2\Dx^2}{N^2} + \tfrac{\Mi[1:i-1]\Li\sigtl[i]^2}{N}] + \tsum_{i \in \frP}[\tfrac{\Mi[1:]\Dx}{N} ]\\
&+\tfrac{\tsum_{i \in \frNi[1:k-1]}[|\frNi| \Mi[1:]\DX + \sigtl[v,i]\Mi[i+1:]\Dx]}{\sqrt{N}} + \tsum_{i \in \frN}[\frac{\Mi[1:]\Dx}{\sqrt{N}}]\\
&+ \tsum_{i=1}^k [\tfrac{\Mi[1:i]\sigtl[i]}{\sqrt{N}}] + \tfrac{\DX\sigtl[x]}{\sqrt{N}}\}.
\end{split}
\end{equation}
\item[b)] If $u(x)$ is strongly convex with  modulus $\alpha >0$ and all outer layer functions are smooth, i.e., $\{1,2,\ldots,k-1\} \subset \frS$, then selecting $\etat:=\max\{\tsum_{i=1}^{k}\etait, \alpha(t-1)/2\}$ leads to 
\begin{align}\label{meq:mixed-str-conv}
\begin{split}
\E[f(\xbarT) - f(\xstar)] \leq \bigO\left([\tfrac{\sum_{i=1}^{k-1}\Mi[1:i-1]\Li \Mi[i+1:]^2}{\alpha}  \right.&+ 1]\{\tfrac{\log(N + 1)}{N^2}[\tsum_{i \in \frS}  \Mi[1:i-1]\Li \Mi[i+1:]^2\DX^2] \\
  &\left.+ \tfrac{1}{N}[ \tsum_{i=1}^{k-1} \Mi[1:i-1] \Li \sigtl[i]^2 +\tsum_{i\in \frP}(\Mi[1:]\Dx)+\tsum_{i \in \frN}\tfrac{\Mi[1:]\Dx}{\alpha}+\tfrac{\sigtl[x]^2}{\alpha}]  \} \right).
\end{split}
\end{align}
\end{enumerate}
\end{theorem}
\begin{proof} The analysis is a straightforward generalization to those of Theorem \ref{thm:l-strns}, \ref{thm:2-ns} and \ref{mthm:ssd}. 
\end{proof}
\vgap
We make three remarks regarding the result. Under the deterministic setting, the terms in \eqref{meq:mixed-conv} attributable to the smooth and the structured non-smooth functions, ${\Mi[1:i-1]\Li\Mi[i+1:]^2\Dx^2}/{N^2}$ and ${\Mi[1:]\Dx}/{N}$, are unimprovable. Under the stochastic setting, as argued in Theorem \ref{thm:l-strns} and \ref{thm:l-ns}, the assumption of all outer layer functions being smooth is necessary for obtaining the improved stochastic oracle complexity of $\bigO(1/\ep)$ in \eqref{meq:mixed-str-conv}. Moreover, the specific stepsize choices listed above are required only to achieve the desired constant dependence. In practice, the following parameter-independent choices can lead to the same order-optimal stochastic oracle complexity of $\bigO(1/\ep^2)$ in the non-strongly convex case:
$$\wt =t, \tauit:= (t-1)/2\ \forall i \in \frS, \tauit:=\Theta(\sqrt{t})\ \forall i \in \frN \cup \frP,\ \gamiit:= \Theta(\sqrt{t})\ \forall i \in \frNo,\ \etat:= \Theta(\sqrt{t});$$
and the oracle complexity of $\bigO(1/\ep)$ in the strongly convex case if all outer-layer functions are smooth:
$$\wt =t, \tauit:= (t-1)/2\ \forall i \in \frS, \tauit:=\Theta(1)\ \forall i \in  \frP,\ \etat:= (t-1)\alpha/2.$$


\subsection{Convergence Analysis}\label{sb:mul-pf}
We present in this subsection the detailed analysis of the SSD method applied to the multi-layer smooth problem. We follow the same pattern as Subsection \ref{sb:2-smo-pf} by first presenting a general convergence result of the $Q$-gap function and then specializing it to provide the convergence rates in Theorem \ref{mthm:ssd}. 

 \begin{proposition}\label{mpr:sm}
 Consider a smooth multi-layer problem of the form \eqref{eq:multi-prob}. Let $\{\zt:=(\xt; \piit[1:])\}$ be generated according to \eqref{mls:SD_idea} with $\Ui = \Dfistar\ \forall i \in [k]$. 
 Assume its Lipschitz-continuity and smoothness constants are defined in \eqref{mcst:sm}, and the stochastic arguments $\yitilt(\xi)\ 0\leq i \leq k-1$  satisfy the aggregate variance bounds in \eqref{mcst:agg_var}. 
 Let $z=(\xstar;\pii[1] = \grad \myfi[1](\yi[1]) \ldots \pii[k] = \grad \myfi[k](\yi[k]))$, for any $y_i \in \R^{n_i}$, be a reference point which could potentially depends on $\{\zt\}$ and let the following requirements be satisfied for all $t \geq 1$ with some non-negative weight $\wt$:
 \begin{align}
  \begin{split}
  &\wt[t] = \thetat[t+1]\wt[t+1],\ \wt(\tauit + 1) \geq \wt[t+1] \tauit[i][t+1] \ \forall i \in [k], \\
&\etat \geq \tsum_{i=1}^k \tfrac{\thetat[t+1] \Mi[1:i-1] \Li \Mi[i+1:]^2}{\tauit[i][t+1]},\  
 \etat[N] \geq \tsum_{i=1}^k \tfrac{ \Mi[1:i-1] \Li \Mi[i+1:]^2}{\tauit[i][N]+1},
  \end{split}
  \end{align} 
  Assuming $\wt=0$, then the following $Q$-gap bound is valid
  \begin{align}\label{meq:sm-Q-conv}
  \begin{split}
  \E[&\sumwt Q(\zt; z)] + {\wt[N](\etat[N] + \alpha)} \normsq{\xt[N] - \xstar}/2 \leq \E\{\sumt [\wt\etat - \wt[t-1](\etat[t-1] + \alpha)] \normsq{\xt[t-1] - \xstar}/2  \} \\
  & + \tsum_{i=1}^k \wt[1]\tauit[i][1] \Dfistar(\pii; \piit[i][0]) + \sumt \wt \sigtl[x]^2 /(\etat + \alpha) + \tsum_{i=1}^{k-1}\{\Mi[1:i-1]\Li \sigtl[i]^2 [ \sumwt   / (\tauit + 1)]\} \\
  & + \tsum_{i=1}^{k-1} \E\{\sumt \wt \pii[1:i-1](\pii - \pii[i,\E])[\yitilt[i]- \yitilt[i](\xi)]\},
  \end{split}
  \end{align}
  where $\pii[i, \E] = \grad \myfi(\yi)$ for some $\yi\in\R^{n_i}$ and $\pii[i, \E]$ is independent of $\yitilt[i]- \yitilt[i](\xi)$ conditioned on $\yitilt[i]$.

\end{proposition}

\begin{proof}
First, let us develop a convergence bound for $Q_i\ \forall i \geq 1$ (c.f. \eqref{meq:gap_func}). The $\piit$ update in \eqref{mls:SD_idea} implies a (vector) three-point inequality given by 
$$(\pii - \piit) \yitilt(\xi) + \fistar(\piit) - \fistar(\pii) + (\tauit + 1)\Dfistar(\pii; \piit) + \tauit \Dfistar(\piit; \piitt) \leq \tauit \Dfistar(\pii; \piitt).$$
Multiplying both sides with the non-negative weight $\pii[p]:=\pii[1:i-1]$ leads to 
\begin{equation}\label{pfeq:pr3-1}
\pii[p][(\pii - \piit) \yitilt + \fistar(\piit) - \fistar(\pii)] + (\tauit + 1)\pii[p]\Dfistar(\pii; \piit) + \tauit \pii[p]\Dfistar(\piit; \piitt) \leq \tauit \pii[p]\Dfistar(\pii; \piitt) + \pii[p](\pii - \piit) \delit,\end{equation}
where $\delit:= \yitilt - \yitilt(\xi).$
Take $\piit[i, \E] \leftarrow \argmax_{\pii \in \Pii} \pii \yitilt - \fistar(\pii) - \tauit U_i(\pii;\piitt)$ such that it is conditionally independent of $\delit$. Since $\pii[p]\piit[i, \E] = \argmax_{\pii \in \Pii} \pii[p]\pii \yitilt - \pii[p]\fistar(\pii) - \tauit \pii[p] \Dfistar(\pii;\piitt)$, we can apply to argument in \eqref{eq:smo-noise-bd}, with the $1/(\norm{\piip}\Li)$ strong convexity of $\pii[p] \Dfistar(\pii;\piitt)$ with respect to $\pii[p]\pii$ (see Lemma \ref{lm:smooth_breg}), to obtain $\norm{\pii[p](\piit - \piit[i, \E])} \leq \norm{\pii[p]}\Li\norm{\delit}/(\tauit +1).$ Therefore 
\begin{equation}\label{pfeq:pr3-2}
\begin{split}
\E[\pii[p](\pii - \piit) \delit] &= \E[\pii[p](\pii - \piit[i, \E]) \delit] + \E[\pii[p](\piit[i, \E] - \piit) \delit]\\
&\leq \E[\pii[p](\pii - \piit[i, \E]) \delit] + \Li\E[\norm{\pii[p]} \normsq{\delit}]/(\tauit +1) \\
&\leq \E[\pii[p](\pii - \piit[i, \E]) \delit] + \Mi[1:i-1] \Li \sigtl[i]^2 /(\tauit + 1).
\end{split}
\end{equation}
Next consider the left-hand-side of \eqref{pfeq:pr3-1}, we have 
$$\pii[p](\pii - \piit) \yitilt = \pii[p](\pii - \piit) \Laitt[i+1][t] - \underbrace{\{\pii[p](\pii - \piit) \piit[i+1:](\xt - \xtt) - \thetat \pii[p](\pii - \piit) \piit[i+1:](\xtt - \xt[t-2]) \}}_{A_t}.$$
Also, an argument similar to that of \eqref{pfpr1:Q1_bound} implies 
\begin{equation}\label{pfeq:pr3-3}
\begin{split}
\sumwt& \{A_t - (\tauit + 1)\pii[p]\Dfistar(\pii; \piit) - \tauit \pii[p]\Dfistar(\piit; \piitt) + \tauit \Dfistar(\pii; \piitt)\} \\
&\leq \tauit[i][1] \Dfistar(\pii; \piit[i][0]) + \tsum_{t=1}^{N-1} \tfrac{\wt\theta_{t+1}}{2\tauit[i][t+1]} (\Mi[1:i-1] \Li \Mi[i+1:]^2) \normsq{\xt - \xtt}\\
&\quad  +  \tfrac{\wt[N]}{2(\tauit[i][N] +1)} (\Mi[1:i-1] \Li \Mi[i+1:]^2) \normsq{\xt[N] - \xt[N-1]}.
\end{split}
\end{equation}
Applying \eqref{pfeq:pr3-2} and \eqref{pfeq:pr3-3} to an $\wt$-weighted sum of \eqref{pfeq:pr3-1}, we obtain 
\begin{equation}
\begin{split}
\E[\sumwt Q_i(\zt; z)] \leq& \tsum_{t=1}^{N-1} \tfrac{\wt\theta_{t+1}}{2\tauit[i][t+1]} (\Mi[1:i-1] \Li \Mi[i+1:]^2) \normsq{\xt - \xtt} +  \tfrac{\wt[N]}{2(\tauit[i][N] +1)} (\Mi[1:i-1] \Li \Mi[i+1:]^2) \normsq{\xt[N] - \xt[N-1]} \\
&\tauit[i][1] \Dfistar(\pii; \piit[i][0]) + \E[\sumwt \pii[p](\pii - \piit[i, \E]) \delit] + \sumwt \Mi[1:i-1] \Li \sigtl[i]^2 /(\tauit + 1).
\end{split}
\end{equation}
The convergence bound for $Q_0$ is similar to \eqref{pfpr1:Q0_bound}. Therefore, adding up the bounds for $Q_k$, $Q_{k-1}$, ..., and $Q_0$ and noting the stepsize requirement for $\etat$, we obtain the desired convergence bound for $Q$ in \eqref{meq:sm-Q-conv}.

\end{proof}

\textbf{Proof of Theorem \ref{mthm:ssd}}
Clearly, Algorithm \ref{malg:ssd} is a concrete implementation of \eqref{mls:SD_idea}. We need to show three bounds: the aggregate variance bounds in  \eqref{meq:agg_var_bound}, the convergence bounds for the non-strongly convex problem in \eqref{meq:sm-ns-conv} and for the strongly convex problem in \eqref{meq:sm-str-conv}.

We first develop the aggregate variance bounds. As argued in Figure \ref{mfig:resampling}, the sampled $\yitilt(\xi)$'s in Algorithm \ref{malg:ssd} are unbiased estimators to $\yitilt$'s in \eqref{mls:SD_idea}. 
We establish upper bounds to their variances (c.f. \eqref{meq:agg_var_bound}) if the stepsize in \eqref{mstp:ssd-dual} is utilized. For simplicity, we use the notation $\Var[y(\xi)]:= \E[\inner{y(\xi) - \E[y(\xi)]}{y(\xi) - \E[y(\xi)]}]$ if $y(\xi)$ is a random vector. The upper bound for $\sigtl[x]$ is straightforward, but that for $\sigtl[i]^2 := \max_{t \geq 0} \Var[\yitilt(\xi)]$ requires more effort. Specifically, the definition of $\yitilt(\xi)$ in Algorithm \ref{malg:ssd} implies that 
\begin{equation}\label{mpf:sigtl-decom}
\begin{split}
\Var[\yitilt[i-1](\xi)] &= \Var[\Laijtt[i][i-1] + \thetat \grad \myfi[i:](\uyitt[i:], \hat{\xi})(\xtt - \xt[t-2])] \\
&= \Var[\Laijtt[i][i]] + \thetat^2 \Var[\grad \myfi[i:](\uyitt[i:], \hat{\xi})(\xtt - \xt[t-2])] \\
& \leq \Var[\Laijtt[i][i]] + \tsum_{l=i}^k [\prod_{j\in i:k/\{l\}} (\Mi[j]^2 + \sig[j]^2)] \sig[l]^2 \DX^2,
\end{split}
\end{equation}
where the second equality follows from the conditional independence between the terms. Regarding the first term above, its definition in Line 6 of Algorithm \ref{malg:ssd} implies the following inequality for all $j < i$,
\begin{equation}\label{mpf:var-ltil}
\begin{split}
\Var&[\Laijtt[i][j]] = \Var[\myfi(\uyit, \xi_i^j) + \grad \myfi[i](\uyit, \xi_i^j)\{\La_{i+1}[\xtt; \piit[\ito[i+1]](\xi_{\ito[i+1]}^j)] - \uyit[i+1]\}] \\
&\leq 2 \Var[\myfi(\uyit, \xi_i^j)] + 2\Var[\grad \myfi[i](\uyit, \xi_i^j)\{\La_{i+1}[\xtt; \piit[\ito[i+1]](\xi_{\ito[i+1]}^j)] - \uyit[i+1]\}] \\
&\leq 2 \sig[i]^2 \Mi[i+1:]^2 \DX^2 + 2 (\sig[i]^2 + \Mi^2) \Var[\Laijtt[i+1]] + 2 \sig[i]^2 \E[\normsq{\Laijtt - \uyit}].
\end{split}
\end{equation}

We provide bounds to $\Var[\Laijtt[i+1]]$ and $\E[\normsq{\Laijtt - \uyit}]$ above. It is useful to define a few constants.
\begin{align}\label{mpf:aux-cst}
\begin{split}
&\sig[\La_i]^2:=\max_{t, t' \geq 0}  \max_{j, j' < i} \Var[\Laijtt].  \\
&\hat{D}^2_i := \max\{\max_{t, t' \geq 0} \max_{j, j' < i}\E [\max_{x, x'} \normsq{\hLaijtt[i][j][] - \hLaijtt[i][j']['][t']}], \Mi[i:]^2 \DX^2\} \text{ where }\\
&\hLaijtt[i][j][]:= \Laijtt +  \gradfijht[i:] (x - \xtt).
\end{split}
\end{align}
Similar to \eqref{eq:uyit-decom}, the stepsize choice in \eqref{mstp:ssd-dual} implies that $\uyit= 2[\tsum_{l=1}^{t-1} l \hLaijtt[i+1][i][l][l] + t \hLaijtt[i+1][i][t-1][t]] / [t(t+1)]$. Thus the Jensen's inequality allows to simplify  \eqref{mpf:var-ltil} can be simplified into a recursive relation 
\begin{equation}\label{mpf:sig-Li-recur}
\sig[\La_i]^2 \leq  
 \begin{cases} 0 &\text{if }i = k+1, \\
 2 \sig[i]^2 \Mi[i+1:]^2 \DX^2 + 2 (\sig[i]^2 + \Mi^2)\sigLai[i+1]^2 + 2 \sig[i]^2 \hat{D}^2_{i+1} &\text{otherwise.}\end{cases}
\end{equation}

Now let's focus on $\hat{D}_i$. We can derive from the recursive definition of $\hat{\La}_i$ in \eqref{mpf:aux-cst} the following recursive bound 
\begin{equation*}
\begin{split}
\hat{D}_i^2 \leq \begin{cases} \DX^2 &\text{if } i= k+1 \\
20(\Mi^2 + \sig[i]^2)\hat{D}_{i+1}^2 + 20 \prod_{l=i}^k (\Mi^2 + \sig[i]^2) \DX^2 &\text{otherwise.}    
\end{cases}
\end{split}
\end{equation*}
Thus 
$\hat{D}_i^2 \leq [1+(k-i + 1) 20^{k-i}] \prod_{l=i}^k (\Mi^2 + \sig[i]^2) \DX^2.$
Substituting it into \eqref{mpf:sig-Li-recur} leads to 
\begin{equation*}
\begin{split}
\sig[\La_i]^2 &\leq \tsum_{l=i}^{k} \{\prod_{j=i}^{l-1}[2 (\Mi[j]^2 + \sig[j]^2)]\} \{2\sig[l]^2\} \{
[2+ (k-l) 20^{k-l-1}] \prod_{j=l+1}^k (\Mi[j]^2 + \sig[j]^2) \DX^2\}\\
&\leq 2 [2 + (k-i)20^{k-i-1} ]  \tsum_{l=i}^{k} \sig[l]^2 \prod_{j\in i:k/\{l\}} (\Mi[j]^2 + \sig[j]^2) \DX^2. 
\end{split}
\end{equation*}
The desired bound on $\sigtl[i]$ then follows from substituting the preceding bound into \eqref{mpf:sigtl-decom}.

Next, with the specific choices of $\etat$ in the theorem statement, the requirements of Proposition \ref{mpr:sm} are satisfied, so the convergence rates in \eqref{meq:sm-ns-conv} and \eqref{meq:sm-str-conv} can be derived in a similar fashion as those of Theorem \ref{thm:2-sm} and \ref{thm:2-sm-str}.
\endproof

\section{Applications}

In this section, we demonstrate the practical use of the SSD and nSSD methods by applying them to two concrete problems. Even though these algorithms are order-optimal, we can tailor them to specific problem structures to improve the constant dependence of the oracle complexity. 

\subsection{Risk Averse Optimization}

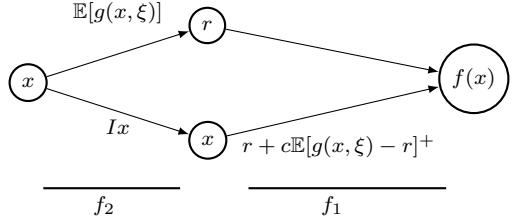
\begin{figure}[!htb]
    \centering
    \begin{minipage}{.5\textwidth}
        \centering
          \begin{tikzpicture}[auto,node distance=8mm,>=latex,font=\small]
            \draw [thick] (0.2,-1.4) -- (2, -1.4);
            \draw [thick] (2.9,-1.4) -- (5.5, -1.4);
            \tikzstyle{round}=[thick,draw=black,circle]

            \node[round] (x) {$x$};
            \node[round,above right=.4 and 2 of x] (y) {$r$};
            \node[round,below right=.4 and 2 of x] (x1) {$x$};
            \node[round,  below right=.2 and 3 of y] (fx) {$f(x)$};

            \draw[->] (x) -- node[above=0.3] {$\E[g(x, \xi)]$} ++ (y);
            \draw[->] (x) -- node[below] {$I x$} ++ (x1);
            \draw[->] (y) -- (fx);
            \draw[->] (x1) --node[below right =.2 and -1.3 ]{$r + c \E[g(x, \xi) - r]^+$} ++ (fx);
            \node[below] at (1, -1.4) {$f_2$};
            \node[below] at (4, -1.4) {$f_1$};
        \end{tikzpicture}
        \caption{Two Layer Formulation for \eqref{apeq:devi}.}
        \label{fig:semi_orig}
    \end{minipage}%

\end{figure}



First, we consider a risk-averse two-stage stochastic program given by 
\begin{equation}\label{apeq:devi}
\tmin_{x \in X} \{\rho(Z(x)) := \E[g(x, \xi)] + c \E[g(x, \xi) - \E[g(x, \xi)]]^+\},
\end{equation}
where the random variable $Z(x):=g(x, \xi)$ denotes the cost incurred by the decision $x$ under the scenario $\xi$, the risk measure $\rho$ is the mean-upper-semideviation of order one \cite{shapiro2014lectures} and the trade-off parameter $c\in[0,1]$ characterizes the optimizer's risk aversiveness.  We take $g(x,\xi)$ to be non-smooth for generality. For example, in the two-stage LP, $g(x,\xi)$, the minimum total cost incurred by the first-stage decision $x$ under the scenario $\xi$, is  piecewise linear.

To apply our nSSD method in Algorithm \ref{alg:ssdp_2}, we can formulate \eqref{apeq:devi} as a two-layer problem shown in Figure \ref{fig:semi_orig}. 
We assume $g(x, \xi)$ to be Lipschtiz continuous and to have bounded variances:
\begin{equation*}
\norm{\E[g'(x, \xi)\mathbbm{1}_{\{g(x, \xi) \geq r\}}]} \leq \Mi[g],\ \E[\normsq{g'(x, \xi) - \E[g'(x, \xi)]}] \leq \sig[g']^2,\ \E[g(x,\xi) - \E[g(x,\xi)]]^2\leq \sig[g]^2\ \forall x \in X\ \forall r \in \R. 
\end{equation*}
To handle the  non-smooth outer-layer function, we need to define a tri-conjugate reformulation to $\myfi[1]$.  We set the domains $\Piitl[1]$ and $\Vi[1]$ to be
\begin{equation}\label{def:domain-risk-averse}
\Vi[1]:= \mathbb{B}^1(\bar r; \Dx\Mg) \times X, \ \Piitl[1]:=[0,1] \times \mathbb{B}^n(0; \Mg), 
\end{equation}
where the notation $\mathbb{B}^m(x;l)$ denotes the $m$-dimensional ball centered at $x$ with a radius $l$, and $\bar r$ denotes any possible value of $E[g(x, \xi)]$ for any $x\in X$. It can be checked that they satisfy all requirements in \eqref{cst:f1-lin} except for $\Piitl[1, 2:] =  \mathbb{B}^n(0; \Mg) \not\subset  \R^{n}_+$. However, the results in Section \ref{sec:2-ns} are still valid because the input corresponding to $\pii[1,2:]$ \vspace{1mm} is a linear function $Ix$. Thus a direct application of Theorem \ref{thm:2-ns} leads to an oracle complexity of
$\bigO\{(\Mg^4 \Dx^4 + \Mg^2 \sig[g]^2 + \sig[g']^2 \Mg^2 \Dx^2)/\ep^2\}.$

To improve the constant dependence, we need to look into the structure of $\pii[1]$ and $\vi[1]$. Their domains are separable Cartesian products of the domains of the first coordinate and those of the next $n$ coordinates (see \eqref{def:domain-risk-averse}). Moreover, in the prox-mapping for $\vi$, the first coordinate of the argument, $1-c\mathbbm{1}_{\{g(\xtt,\xi_1)\geq \vit[t-1]\}}$, differs greatly from the next $n$ coordinates of the argument, $cg'(\xtt, \xi_1) \mathbbm{1}_{\{g(\xtt,\xi_1)\geq \vit[t-1]\}}$. In the prox-mapping for  $\pii[1]$, only the first coordinate of the argument, $g(\xtt, \xi_2^1)$, is stochastic. These separable structures motivate us to modify Line 3 and 4 of Algorithm \ref{alg:ssdp_2} to use different stepsizes for updating the first coordinates:
\begin{align*}
&\piit[1,1] \leftarrow \targmax_{\pii[1,1] \in \Piitl[1,1]}\ \pii[1,1][\yitilt[1,1](\xi) - \viit[1,1][t-1]] - \tauit[1]^r \normsq{\pii[1,1]- \piitt[1,1]}/2,\ \\
&\piit[1,2:] \leftarrow \targmax_{\pii[1,2:] \in \Piitl[1,2:]}\ \pii[1,2:][\yitilt[1,2:](\xi) - \viit[1,2:][t-1]] - \tauit[1]^\grad \normsq{\pii[1,2:]- \piitt[1,2:]}/2,\\
&\viit[1,1] \leftarrow \targmin_{\vi[1,1] \in \Vi[1,1]} \ \inner{\myfi[1,1]'(\vit[t-1], \xi_1) - \piit[1,1]}{\vi[1,1]} + \gamit^r\normsq{\vi[1,1] - \viitt[1,1]}/2,\\
&\viit[1,2:] \leftarrow \targmin_{\vi[1,2:] \in \Vi[1,2:]} \ \inner{\myfi[1,2:]'(\vit[t-1], \xi_1) - \piit[1,2:]}{\vi[1,2:]} + \gamit^\grad \normsq{\vi[1,2:] - \viitt[1,2:]}/2.
\end{align*}
Specifically, if we set $\tauit[1]^\grad = 0$, $\tauit[1]^r= \sig[g]\sqrt{t}$, $\gamit[1]^\grad = \max\{\Mg, \sig[g']\}\sqrt{t}/\Dx$, and $\gamit^r= \sqrt{t}/(\Dx \Mg),$ the oracle complexity can be improved to 
$\bigO(\{\Mg^2 \DX^2 + \sig[g]^2 + \sig[g']^2\Dx^2\}/\ep^2).$
Comparing it to the $\bigO(\{\Mg^2 \DX^2 + \sig[g']^2\Dx^2\}/\ep^2)$ oracle complexity for solving the risk-neutral two-stage program,  the only extra cost is $\bigO(\sig[g]^2/\ep^2)$, which arises from estimating the function value when computing the risk measure $\rho$.
 

\subsection{Stochastic Composite Optimization}

\begin{figure}[!htb]
    \centering
    \begin{minipage}[b]{.5\textwidth}
        \centering
          \begin{tikzpicture}[auto,node distance=8mm,>=latex,font=\small]
            \draw [thick] (0.2,-1.4) -- (2, -1.4);
            \draw [thick] (2.9,-1.4) -- (4.8, -1.4);
            \tikzstyle{round}=[thick,draw=black,circle]

            \node[round] (x) {$x$};
            \node[round,above right=.4 and 2 of x] (z) {$r$};
            \node[round,below right=.4 and 2 of x] (x1) {$x$};
            \node[round,  below right=.2 and 2 of y] (fx) {$f(x)$};

            \draw[->] (x) -- node[above=0.3] {$A x$} ++ (z);
            \draw[->] (x) -- node[below] {$I x$} ++ (x1);
            \draw[->] (z) -- (fx);
            \draw[->] (x1) --node[below right =.2 and -.5 ]{$F(r) + g(x)$} ++ (fx);
            \node[below] at (1, -1.4) {$f_3$};
            \node[below] at (4, -1.4) {$f_2$};
        \end{tikzpicture}

        \caption{Two Layer Formulation for \eqref{apeq:composite_pro}.}
        \label{fig:composite}
    \end{minipage}%
    \begin{minipage}[b]{0.5\textwidth}
        \centering
          \begin{tikzpicture}[auto,node distance=8mm,>=latex,font=\small]
            
            \tikzstyle{round}=[thick,draw=black,circle]
            \tikzstyle{rect}=[thick,draw=black,rectangle]

            \node[round] (x) {$x$};

            \node[below right=0 and 3 of x] (fd) {$\vdots$};
            \node[rect,above=0.5 of fd] (f2) {$f^{(2)}(x) \quad$};
            \node[rect, above=0.5 of f2] (f1) {$f^{(1)}(x) \quad$};
            \node[rect, below=0.5 of fd] (fm) {$f^{(m)}(x) \quad$};
            \draw[->] (x)--(2.6, 1.7);
            \draw[->] (x)--(2.6, 0.65);
            \draw[->] (x)--(2.56, -1.8);

            \node[round, right=6 of x] (fx) {$f(x)$};
            \draw[->] (4.2, 1.7)--(fx);
            \draw[->] (4.2, .65)--(fx);
            \draw[->] (4.2, -1.8)--node[above=0.4] {max} ++(fx);

            \draw[thick] (2.6, -2.5) --  (3.1, -2.5);
            \node[below] at (2.85, -2.5) {$f_3$};
            \draw[thick] (3.2, -2.5) -- (4, -2.5);
            \node[below] at (3.6, -2.5) {$f_2$};
            \draw[thick] (4.5, -2.5)  -- (6, -2.5);
            \node[below] at (5.3, -2.5) {$f_1$};

        \end{tikzpicture}
        \caption{Three Layer Formulation for \eqref{apeq:muli_comp}.}
        \label{fig:composite_epdd}
    \end{minipage}
\end{figure}

\begin{algorithm}[htb]
\caption{SSD Algorithm for Composite Optimization}
\label{apalg:composite}
\begin{algorithmic}[1]
\Require $x_{-1} = x_{0}\in X$ and $\piiO[F] \in \Pii[F]$. 
\parState{Set $\uyiO[g] := x_{0}$ and call $\SOi[]$ to obtain estimate $A(\xi^{2}_{3,0})$.}
\For{$t =  1,2,3 ... N$}
\parState{Call $\SOi[]$ to obtain estimates $A(\xi^{2}_{3,t})$ and ${A^{\transpose}(\xi^{0}_{3,t})}$.}
\parState{Let $\xtilt:= \xtt + \thetat (\xtt - \xt[t-2])$.\\
Let $\uyit[g]:= (\tauit[g] \uyitt[g] + \xtilt)/(1 + \tauit[g])$ and call $\SOi[]$ to obtain $\piit[g](\xi^{0}_{2,t}):= g'(\uyit[g],\xi^{0}_{2,t})$.}

\parState{Let $\yitilt[F](\xi): = A(\xi^{2}_{3,t})\xtt + A(\hxi[3,t-1])(\xtt -\xt[t-2]) $. \\
Compute $\piip[F]:= \argmin_{\pii[F] \in \Pii[F]} -\inner{\pii[F]}{\yitilt[F](\xi)} + F^*(\pii[F]) + \tauit[F] \normsq{\pii[F] - \piitt[F]}/2$}

\parState{Set $\xt := \argmin_{x \in X}  \inner{\piit[g](\xi^{0}_{2,t}) + A^{\transpose}(\xi^{0}_{3,t})\piit[F]}{x}+  \etat \normsq{x - \xt}/2$.}
\EndFor
\State Return $\xbarT := \sumwt \xp / \sumwt$.
\end{algorithmic}
\end{algorithm}

Next we consider a stochastic composite optimization problem that arises frequently in machine learning and data analysis \cite{chen2014optimal}:
\begin{equation}\label{apeq:composite_pro}
\min_{x \in X} \{f(x):= F(Ax) + g(x)\equiv \max_{\pii[F] \in \Pii[F]} \inner{\pii[F]}{ Ax} - F^*(\pii[F])+ g(x)\},\end{equation}
where $F$ is a structured non-smooth function, for example, the total variation loss function, and $g$ is a stochastic smooth function, for example, the data fidelity loss function.  The dimension of $A$ is usually large, so we assume existence of stochastic oracles to return unbiased estimators $A(\xi)$, $A\transpose(\xi)$, and $g'(x, \xi)$ for $A$, $A{\transpose}$ and $g'(x)$, respectively. Additionally, we assume their variances to be uniformly bounded by $\sigma^2_{A}$, $\sigma^2_{A{\transpose}}$
 and $\sig[g']^2$ respectively. 

Clearly, \eqref{apeq:composite_pro} can be formulated as a two layer problem shown in Figure \ref{fig:composite}. Treating the outer $\myfi[2]$ as a non-smooth layer function, a direct application of the nSSD method in Algorithm \ref{alg:ssdp_2} then leads to an order-optimal oracle complexity of $\bigO(1/\ep^2)$. To further improve the constant dependence, we observe a separable structure of $\myfi[2]$. Since $\myfi[2](x, r):= F(r) + g(x)$, the dual variable $\pii[2]$ and the Fenchel conjugate can be decomposed to components attributable to $F$ and $g$ respectively:
$$\pii[2]:=[\pii[F]| \pii[g]], \quad \fistar[2](\pii[2]):=\max_{x, r} \inner{x}{\pii[g]} + \inner{r}{\pii[F]} - F(r) - g(x) = F^*(\pii[F]) + g^*(\pii[g]). $$
The gap-function $Q_2$ is also decomposable 
$$Q_{2}(\zt, z):= \{(\pii[F]- \piit[F])Ax - [F^*(\pii[F])- F^*(\piit[F])]\} + \{(\pii[g]- \piit[g])x - [g^*(\pii[g])- g^*(\piit[g])]\}.$$ 
So we can perform the prox-mappings for $\pii[F]$ and for $\pii[g]$ separately, in accordance with their own structural properties. This leads us to Algorithm \ref{apalg:composite}.
 With appropriately chosen stepsizes, it is easy to derive from Theorem \ref{mthm:ssdp}  an oracle complexity of 
$$\bigO\{ \tfrac{\sqrt{L_g} \norm{x - \xstar}} {\sqrt{\ep}} + \tfrac{\norm{A} \Dx \Mi[F]}{\ep} + \tfrac{(\sigma_A^2 + \sigma_{A{\transpose}}^2) \Dx^2 \Mi[F]^2}{\ep^2} + \tfrac{\sig[g']^2 \Dx^2}{\ep^2}\},$$
where $\Lg$ is the Lipschitz-smoothness constant of $g$ and $\Mi[F]$ is an upper bound for $\norm{\pii[F]}\ \forall  \pii[F] \in \Pii[F]$.

Comparing it to the accelerated primal dual (APD) algorithm \cite{chen2014optimal} designed specifically for problem~\eqref{apeq:composite_pro}, Algorithm \ref{apalg:composite} achieves the same oracle complexity.  However, our general approach allows an easy extension to handle more complicated problems in which $f$ in \eqref{apeq:composite_pro} is but one sub-component. In particular, if $f\upperi:= g\upperi(x) + F\upperi(A\upperi x)$, a minimax sub-problem which arises frequently from either constrained optimization or multi-objective optimization is given by 
\begin{equation}\label{apeq:muli_comp}
\min_{x \in X}\{f(x):= \max \{f\upperi[1](x), f\upperi[2](x),\ldots, f\upperi[m](x)\}\equiv \max_{\pii[1] \in \Delta^+_m} \sum \pii[1]\upperi f\upperi(x)\}\footnote{$\Delta^m_+:= \{\pii[1] \in \R^m_{+}| \sum_{i=1}^{m} \pii[1]\upperi = 1\}$ is the probability simplex.}.\end{equation}
Clearly, \eqref{apeq:muli_comp} admits a three-layer formulation shown in Figure \ref{fig:composite_epdd}, so we only need to add an additional prox-mapping for $\pii[1]$ to Algorithm \ref{apalg:composite}. Moreover, if  the variances for $g^{(i)}(\cdot, \xi_i)$, $\grad g^{(i)}(\cdot, \xi_i)$, $A^{(i)}(\xi_i)$, ${A^{(i)}}\transpose(\xi_i)$ are uniformly bounded by $\sigma_g^2$, $\sig[g']^2$, $\sigma^2_{A}$ and $\sigma^2_{A{\transpose}}$, and if $\norm{A^{(i)}}$, $\Mi[F^{(i)}]$ and  $\norm{\grad g^{(i)}(\cdot)}$ are uniformly bounded by $\norm{A}$, $\Mi[F]$ and $M_g$, then a straightforward application of Theorem \ref{mthm:ssdp} implies an oracle complexity of 
$$\bigO\{ \tfrac{\sqrt{L_g} \norm{x - \xstar}} {\sqrt{\ep}} + \tfrac{{\sqrt{m}}\norm{A} \Dx \Mi[F]}{\ep} + \tfrac{{{m}}(\sigma_A^2 + \sigma_{A{\transpose}}^2) \Dx^2 \Mi[F]^2}{\ep^2} + \tfrac{{m}{(\sig[g']^2 \Dx^2 + {\sigma_g^2})}}{\ep^2} + {\tfrac{\sqrt{m} M_g}{\ep}}\} .$$
If the entropy Bregman distance function is selected as the prox-function for the $\pii[1]$-prox mapping similar to \cite{zhang2019efficient}, the above complexity can be improved to be nearly independent of the number of sub-components, $$\bigO\{ \tfrac{\sqrt{L_g} \norm{x - \xstar}} {\sqrt{\ep}} + \tfrac{{\sqrt{\log (m)}}\norm{A} \Dx \Mi[F]}{\ep} + \tfrac{{{\log(m)}}(\sigma_A^2 + \sigma_{A{\transpose}}^2) \Dx^2 \Mi[F]^2}{\ep^2} + \tfrac{{\log (m)}{(\sig[g']^2 \Dx^2 + {\sigma_g^2})}}{\ep^2} + {\tfrac{\sqrt{\log(m)} M_g}{\ep}}\} .$$.

\section{Conclusion}

To sum up, this paper studies the order of stochastic oracle complexity for the convex NSCO problem by proposing order-optimal algorithms under a mild compositional convexity assumption. Our complexity results reveal that the convex NSCO problem has the same order of oracle complexity as those without the nested composition in all but the strongly convex and outer-non-smooth problem. The proposed SSD/nSSD method is general since it can handle an arbitrary multi-layer composition of smooth, structured non-smooth and general non-smooth layer functions. Moreover, we introduced two motivating applications to show the method is flexible enough to allow modifications to exploit their special problem structures.




\renewcommand\refname{Reference}

\bibliographystyle{siam}
\bibliography{Reference.bib}
\section{Appendix}
\subsection{Technical lemmas for vector-valued functions}
\textbf{Proof of Lemma \ref{lm:proximal_gradient_lemma}}
This is a direct consequence of the conjugate duality relationship \cite{beck2017first}, 
$$\pii \in \partial g_i(\uyi[]) \Leftrightarrow \pii \in \argmax_{\piibar \in \dom(g_i^*)} \piibar \uyi[] - g_i^*(\piibar).$$
 In particular,  since $\uyitt[] \in \partial g^*(\piitt[])$, the $i$\ts{th} row of $\piit[]$ satisfies 
\begin{align*}
&\piit \in \argmin_{\pii \in \dom(g_i^*)} -\inner{\pii}{\yit[]} + g_i^*(\pii) + \taut D_{g_i^*}( \pii; \piitt), \\
\stackrel{}{\Longleftrightarrow}\ \ & \piit \in \argmin_{\pii \in \dom(g_i^*)} - \inner{\pii}{\yit[]} - \taut \inner{\pii}{ \uyitt[]} + (1 + \taut) g_i^*(\pii),\\
\Longleftrightarrow\ \ & \piit \in \argmin_{\pii \in \dom(g_i^*)} - \inner{\pii}{ \underbrace{\tfrac{\yit[] + \taut \uyitt[]}{1 + \taut}}_{\uyit[]}} + g_i^*(\pii) \ \stackrel{}{\Longleftrightarrow} \ \piit = g_i'(\uyit[]) \in \partial g_i(\uyit[]).
\end{align*}
Therefore, we have $\pi^t = g'(\uyit[]).$

\endproof

\begin{lemma}\label{lm:smooth_breg}
Let a closed convex and proper vector-vectored function $g$ be defined on $R^n$. Let it be $L_g$-smooth , i.e., 
$$\norm{\grad g(y) - \grad g(\ybar)}\footnote{$l_2$ operator norm.} \leq \Lg \norm{y - \ybar}, \forall y, \ybar \in \R^n.$$   
Let $\gstar$ and $\Dgstar$ denote its (component-wise) conjugate function and (component-wise) conjugate Bregman's distance function. Then given an m-dimensional non-negative weight vector $w$, we have 
 \begin{equation}\label{msq:1}
 \norm{w}w\transpose \Dgstar(\pibar; \pi) \geq \frac{\norm{w\transpose (\pibar - \pi) }^2}{2 \Lg }
 \end{equation}
if $\pibar = \grad g(\ybar)$ and $\pi = \grad g(y)$ for some $ \ybar, y$.
\end{lemma}

\begin{proof}
First, if $w=0$, \eqref{msq:1} is clearly true. 
Now, assume  $w\neq 0$. The definition of operator norm implies
$$\norm{u\transpose(g'(y) - g'(\ybar))} \leq \Lg \norm{y-\ybar} \ \forall u \text{ with } \norm{u}=1.$$
So the one-dimensional $\guw(y) := \uw\transpose g(y)$ with $\uw := \frac{w}{\norm{w}}$ is $\Lg$-Lipschitz smooth and  its Fenchel conjugate $\guwstar$ is $1/\Lg$ strongly convex. More specifically, since $\guw'(y)=\uw\transpose g'(y)$, we have $$\guwstar(\uw\transpose \pi) - \guwstar(\uw\transpose \pibar) - \uw \transpose (\pi - \pibar)\ybar \geq \frac{1}{2 \Lg} \normsq{\uw \transpose (\pi - \pibar)}, \text{ if } \pibar = g'(\ybar) \text{ , i.e., } \uw\transpose \pibar = \guw'(\ybar). $$
Thus the key to showing \eqref{msq:1} is to relate $\guwstar(\uw\transpose \pi) $ to $\uw \transpose \gstar(\pi)$.

Those two quantities are quite different in general. For $\guwstar(\uw\transpose \pibar):=\max_{y} \uw\transpose\pibar y - \guw(y)$, we can choose only one overall  maximizer, $y^*$, but for $\uw\transpose g^*(\pibar)= \tsum_{j}u_{w,j} \max_{y_j} \pibar_j y_j - g_j(y_j)$, different maximizers $\bar{y}_j^*$'s can be selected for different $\pibar_j$'s. So we have $\guwstar(\uw\transpose \pi) \leq \uw \transpose \gstar(\pi)\ \forall \pi$.  
However, if $\tilde{\pi}$ is associated with some primal solution $\tilde y$, i.e.,  $\tilde\pi =  g'(\tilde y)$, all those $\tilde y_j$'s are the same. Let $\tilde \pi=g'(\tilde y)$ such that $\tilde{\pi}_j = g_j'(\tilde y)$. The conjugate duality implies that $\gstar(\tilde \pi)= \tilde \pi \tilde y - g(\tilde y)$, 
 so $$\guwstar(\uw \transpose \tilde \pi) := \max_{\ybar} \uw\transpose\tilde \pi \ybar - \guw(\ybar) \geq  \uw \transpose (\tilde \pi \tilde y - g(\tilde y)) = \uw \transpose \gstar(\tilde \pi).$$ 
Therefore, $\guwstar(\uw\transpose \tilde \pi) =\uw \transpose \gstar(\tilde \pi)$ holds if $\tilde \pi$ is associated with some $\tilde{y}$. 

Now returning to the $\pi$ and $\pibar$ in the lemma statement, since $\pibar=g'(\ybar)$ and $\pi = g'(y)$, we have
\begin{align*}
    \uw \transpose \Dgstar(\pibar, \pi) &= \uw\transpose\gstar( \pi) - \uw\transpose\gstar(\pibar) - \uw \transpose (\pi - \pibar)\ybar\\
    &= \guwstar(\uw \pi) - \guwstar(\uw \pibar) - \uw \transpose (\pi - \pibar)\ybar\geq \frac{1}{2 \Lg} \normsq{\uw \transpose (\pi - \pibar)}.
\end{align*}
Then \eqref{msq:1} follows from multiplying both sides of the above inequality by $\norm{w}^2$.
\end{proof}
\vgap
The following lemma provides a bound on the error incurred from utilizing a stochastic argument during prox update.
\begin{lemma}\label{lm:prox-noise-bd}
  Let $\Pi \subset \R^m$ be a non-empty closed and  convex domain and let function $u(\pi)$ be $\mu$-strongly convex. Let $\hat \pi$ be generated via a prox-mapping with the argument $g + \delta$,
$\hat \pi \leftarrow \argmin_{\pi \in \Pi} \inner{\pi}{g + \delta} + u(\pi),$
  where  $\delta$ denote a noise term with $\E[\delta] = 0$ and $\E[\normsq{\delta}] \leq \sig^2.$ Then $|\E[\inner{\hat \pi}{\delta}]| \leq \sig^2/\mu$.
\end{lemma}
\begin{proof}
Let $\pi(y):=\argmin_{\pi \in \Pi} \inner{\pi}{y} + u(\pi)$. The $\mu$-strong convexity of $u(\pi)$ implies that $\pi(y)$ is an $1/\mu$-Lipschitz continuous function of $y$ (see \cite{Nestrov2004Smooth}). Define a auxiliary point
 $\piibar[] := \argmin_{\pi \in \Pi} \inner{\pi}{g} + u(\pi)$ which is independent of $\delta$, i.e., $\E[\inner{\delta}{\piibar[]}] = 0$. The $1/\mu$-Lipschitz continuity of $\pi(y)$ implies  
 \begin{equation}\label{eq:smo-noise-bd}
  \norm{\piibar[] - \hat \pi} \leq \norm{\delta}/\mu.\end{equation} 
 Thus we get 
 \begin{align*}
 |\E \inner{\hat \pi}{\delta}| &\leq  |\E \inner{\hat \pi - \piibar[]}{\delta}| + |\E \inner{\bar \pi}{\delta}| \\
 &\leq \E[ \norm{\bar \pi - \hat \pi}\norm{\delta} ] 
 \leq \E [\normsq{\delta}/\mu]
 \leq \sigma^2 / \mu.
 \end{align*}

\end{proof}

\subsection{Lower Complexity Bounds}
We present in this subsection the detailed analysis for lower complexity results in Subsection \ref{sb:2-smo-l} and  \ref{sb:2-ns-lower} . First, we state an one-dimensional hard problem which will be useful for proving both results.  Parameterized by $(\alpha, \sigma, \nu, \beta)$ with $\sigma \geq \nu$, the problem is given by

 \begin{minipage}{0.4\textwidth}
\centering
  \begin{tikzpicture}[scale=0.8]

   \begin{axis}[
        xmin=-1, xmax=1,
        ymin=-1.3, ymax=2.5,
        xtick distance=10, ytick distance=10, axis x line=middle, axis y line=middle]

      \addplot [domain=-1:-0.5, samples=100, thick, color=blue!50]
        {-1 };

      \addplot [domain=-.5:1, samples=100, thick, color=blue!50]
        {2*x};


      \draw [dashed, color=red!50, line width=0.3mm,opacity=0.8] (axis cs:{-0.5,0}) -- (axis cs:{-0.5,-1});
      \draw [dashed, color=red!50, line width=0.3mm, opacity=0.8] (axis cs:{-0,-1}) -- (axis cs:{-0.5,-1});
      \node[color=black, font=\footnotesize] at (axis cs: 0.1,2) {$\myfi[1]$};
      \node[color=black, font=\scriptsize] at (axis cs: -.5, .2) {$-\nu$};
      \node[color=black, font=\scriptsize] at (axis cs: 0.1, -1) {$-\nu\beta$};
    \end{axis} 
\end{tikzpicture}

\end{minipage}
\begin{minipage}{.6\textwidth}
\begin{align}\label{eq:hard_prob}
\begin{split}
f(x)& := \myfi[1](\E[\myfi[2](x, \xi_2)]) + \alpha \normsq{x}/2, \ X := [-2\nu, 2\nu],xXX \text{ where }\\
&\myfi[1](y_1) :=\beta  \max\{y_1, - \nu\}, \\
&\myfi[2](x,\xi_2):= x + \xi_2,\\
& \text{\small with iid r.v. } \xi_2 :=\begin{cases}-\nu & \text{w.p. } 1-q \\ \nu(1-q)/q &\text{w.p. } q, \end{cases} \\
& \text{\small and }q:= \nu^2 / \sig^2.
\end{split}
\end{align}
\end{minipage}
Observe that $\myfi[1]$ is $\beta$-Lipschitz continuous and the variance of $\myfi[2](x, \xi_2)$ satisfies
$$\E\normsq{\myfi[2](x, \xi_2) - \E[\myfi[2](x, \xi)]} = \E[\normsq{\xi_2}] \leq \sigma^2.$$
The key to our construction is to show the reachable subspace of $x$ being restricted to $0$, i.e., $\Xt= \{0\}$, if a certain condition is met for all generated stochastic estimators. The next technical lemma characterizes the probability and the optimality gap of that scenario. 

\begin{lemma}\label{lm:lower}
The following results are valid for \eqref{eq:hard_prob}. \\
a) $f(0) - f(\xstar) \geq \min\{\beta\nu, \beta^2/\alpha\}/2.$\\
b) If $\xi_2^j$ denotes the $j$\ts{th} query to $\SOi[2]$ and $N <\sigma^2/(4\nu^2)$,  then $\probP\{\myfi[2](0, \xi_2^l) =-\nu\ \forall l \leq N\}> 1/2.$
\end{lemma}

\begin{proof}
Part a) can be derived from the first order optimality condition. Since $f(x)=\beta \max\{x, - \nu \} + \alpha \normsq{x}/2$, the optimal solution and the optimal objective value are 
\begin{equation*}
    \xstar = \begin{cases}-\beta/\alpha & \text{if } \alpha > \beta/\nu \\
    -\nu &\text{if } \alpha \in \tfrac{\beta}{\nu} [0, 1], \end{cases} \Rightarrow f(\xstar) \leq \begin{cases}-\beta^2/(2\alpha) & \text{if } \alpha > \beta/\nu \\
    -\beta\nu/2   &\text{if } \alpha \in \tfrac{\beta}{\nu} [0, 1]. \end{cases}
\end{equation*}
So the inequality in part a) represents an uniform lower bound on $f(0) - f(\xstar).$ 
Part b) follows from the algebraic fact that  $(1 - p)^t >  3/4 - tp$ if $t p <  1/4$:
\begin{align*}
\probP\{\myfi[2](0, \xi_2^l) =-\nu\ \forall l \leq N\} = (1-q)^N > 3/4 - N q > 1/2.
\end{align*}
\end{proof}

Now we are ready to prove the lower bound results. \\

\textbf{Proof for Theorem \ref{thm:l-strns}}: Consider applying the abstract scheme in \eqref{alg:abs-strns} to the problem in \eqref{eq:hard_prob}. A structured non-smooth formulation to $\myfi[1]$ is given by 
$$\myfi[1](\yi[1])= \tmax_{\pii[1] \in [0,\beta]}  \pii[1]\yi[1] -  \nu(\beta-\pii[1]),$$
where $\Pii[1] = [0, \beta]$ and $\fistar[1](\pii[1]) = \nu(\beta - \pii[1]).$ 
Choosing $\yit[1][0] = - \nu$, we have $\Xt[0]= \Piit[0] = \{0\}$ and $\Yit[0]= \{-\nu\}.$ Now assume  $\Xt[t-1]= \Piit[t-1] = \{0\}$ and $\Yit[t-1]= \{-\nu\}$, then $\myfi[2](0, \xi_2^t) = -\nu$ (c.f. \eqref{alg:abs-strns}) implies $\Yit = \{-\nu\}$ and $\Piit= \{0\}$ since 
\begin{align*}
\piit[1]&=\targmax_{\pii[1] \in \Pii[1]} \inner{\pii[1]}{y_1^t} - \fistar[1](\pii[1]) - \tauit[1] \normsq{\pii[1] - \piitt[1]}/2 \\
&= \targmax_{\pii[1] \in [0, \beta]} \pii[1](\yit[1] + \nu) - \beta \nu - \tauit[1] \normsq{\pii[1] - \piitt[1]}/2
= \targmin_{\pii[1] \in [0, \beta]} \tauit[1] \normsq{\pii[1]}/2 = 0.
\end{align*}
It then follows $\Xt = \{0\}$. Such an argument can be applied recursively to show $\Xt[N]=\{0\}$ if the event $B^N :=\{\myfi[2](0, \xi_2^l) =-\nu\ \forall l \leq N\}$ occurs.

Now selecting $\beta =\Mi[1]$, $\nu = 4\ep/\Mi[1]$, $\alpha = \bar{\alpha}$ and $\sigma= \sig[\myfi[2]]$, the hard problem in \eqref{eq:hard_prob} satisfies the hard problem requirements in the theorem statement. Moreover, with $N < \sig[\myfi[2]]^2\Mi[1]^2/(4\ep^2)$, Lemma \ref{lm:lower} implies $\probP(B^N)> 1/2$ such that 
\begin{equation}
\E[f(x^N) - f(\xstar)] \geq \probP(B^N) \E[f(0) - f(\xstar)]> \min\{\beta \nu, \beta^2/\alpha\} /4 = \ep.
\end{equation}
Thus it takes at least $\Omega(\Mi[1]^2\sig[\myfi[2]]^2/\ep^2)$ $\SOi[2]$ queries to obtain an $\ep$-optimal solution.
\vgap

\textbf{Proof for Theorem \ref{thm:l-ns}}
The analysis is similar to that of Theorem \ref{thm:l-strns}. With $\vit[0]=-\nu$, we need to show $\Mit[1][N]=\Xt[N]=\{0\}$ if $\myfi[2](0, \xi^j_2) = - \nu\ \forall j \leq N$.
\endproof
\vgap
We remark that the lower complexity bound of $\Omega(\Mi[1]^2\sig[f_2]^2/\ep^2)$ is applicable beyond  the first-order schemes like \eqref{alg:abs-strns} and \eqref{alg:abs-ns}. In fact, it is not hard to use
 the hard instance in \eqref{eq:hard_prob} to show that at least $\Omega(\Mi[1]^2\sig[f_2]^2/\ep^2)$ samples are required by any SAA-type method to find an $\ep$-optimal solution for the strongly convex NSCO problem with either a structured non-smooth or a general non-smooth outer-layer function.
\end{document}